\documentclass[11pt, a4paper, reqno]{amsart}

  \usepackage{verbatim, ifpdf}
  \usepackage{graphicx, color}

 \usepackage{amssymb, amsmath}
 \usepackage{enumitem}
 \usepackage[a4paper,margin=30mm]{geometry} 

 \usepackage{amsthm}
 \usepackage{listings}
 \usepackage{bbm}
 \usepackage{mathtools}

 \usepackage{tikz-cd}

 \usepackage{soul}

 \usepackage{stmaryrd}

 \usepackage{enumitem}
 \usepackage{framed}

 \usepackage{array, tabularx, multirow}

 \usepackage{tikz}
\usepackage{adjustbox}

\usepackage[table]{xcolor}

\usepackage{appendix}
\usepackage{todonotes}

 \usetikzlibrary{decorations.markings}

  \IfFileExists{hyperref.sty}{%
    \ifpdf
      \usepackage[pdftex]{hyperref}
    \else
      \usepackage[dvips]{hyperref}
    \fi
    \definecolor{maroon}{RGB}{128,0,0}
    \hypersetup{pdfborder=0 0 0,colorlinks=true,plainpages=false,
                pageanchor=true,linkcolor=black,urlcolor=maroon, citecolor=colorForCitations}%
  }{%
    \providecommand{\href}[1]{\texttt{##1}}
    \providecommand{\href}[2]{##2\footnote{~See \texttt{##1}.}}%
  }

\newtheorem{theorem}{Theorem}[section]
\newtheorem{proposition}[theorem]{Proposition}
\newtheorem{lemma}[theorem]{Lemma}

\newtheorem{corollary}[theorem]{Corollary}

\theoremstyle{remark}
\newtheorem{remark}[theorem]{Remark}
\theoremstyle{definition}

\newtheorem{definition}[theorem]{Definition}
\newtheorem{example}[theorem]{Example}

\newcommand{\m}{\mathsf{m}}
\newcommand{\len}{\ell}
\newcommand{\find}{{\mathbf{k_0}}}
\newcommand{\psind}{{\boldsymbol{\rho_0}}}

\newcommand{\ZZ}{\mathbb{Z}}
\newcommand{\QQ}{\mathbb{Q}}
\newcommand{\RR}{\mathbb{R}}
\newcommand{\CC}{\mathbb{C}}

\newcommand{\cP}{\mathcal{P}}
\newcommand{\cPmin}{{\mathcal{P}_{\mathrm{min}\, a}}}
\newcommand{\cPminmax}{{\mathcal{P}_{\mathrm{min}\, a}^{\mathrm{max}\, \len}}}
\newcommand{\vdim}{\mathrm{vdim}}
\newcommand{\Gdec}{{\overrightarrow{\Gamma}}}
\newcommand{\src}{\mathrm{src}}
\newcommand{\dst}{\mathrm{dst}}
\newcommand{\Lie}{\mathrm{Lie}}
\newcommand{\cM}{\overline{\mathcal{M}}}
\newcommand{\ev}{\mathrm{ev}}
\newcommand{\dec}[1]{\overrightarrow{#1}}
\newcommand{\val}{\mathrm{val}}
\newcommand{\tP}{{\widetilde{P}}}
\newcommand{\cO}{{\mathcal{O}}}
\newcommand{\cN}{{\mathcal{N}}}
\newcommand{\PP}{{\mathbb{P}}}
\newcommand{\Th}{{\mathrm{Th}}}
\newcommand{\GKM}{{\mathrm{GKM}}}
\newcommand{\Gfat}{{\boldsymbol{\mathcal{G}}}}
\newcommand{\G}{{\mathcal{G}}}
\newcommand{\C}{\mathcal{C}}
\newcommand{\Cmin}{\mathcal{C}^\text{min}}
\newcommand{\cQ}{\mathcal{Q}}

\definecolor{tableRow}{HTML}{FAFAFA}
\definecolor{tableRowGray}{HTML}{F3F3F3}
\definecolor{tableHead}{HTML}{DEF6AA}
\definecolor{colorForCitations}{HTML}{00643A}

\begin{document}

\title[Hamiltonian GKM spaces via GW theory]{Hamiltonian GKM spaces via Gromov--Witten theory}

\author[D. Holmes]{Daniel Holmes}
\address{Institute of Science and Technology Austria (ISTA)}
\email{daniel.holmes@ist.ac.at}
\urladdr{https://www.daniel-holmes.at}



\begin{abstract}
We use equivariant Gromov--Witten theory to prove the generalised symplectic Mukai inequality for Hamiltonian GKM spaces with \emph{large pseudoindex} or a \emph{special vertex}, including smooth Schubert varieties.
We also classify all Hamiltonian GKM graphs when the pseudoindex is at least the complex dimension, obtaining a Hamiltonian GKM analogue of the classical Kobayashi--Ochiai classification of complex projective space and complex quadrics.

Of independent interest are our connectivity theorem for a low-magnitude subgraph of the GKM graph (reminiscent of Campana and Koll{\'a}r--Miyaoka--Mori's rational chain connectedness of Fano manifolds), and a local Mukai-type theorem that applies at each vertex.
We obtain corollaries on reflexive Delzant polytopes that are purely combinatorial in nature.
\end{abstract}

\maketitle

\tableofcontents

\section{Introduction}
\label{sec : intro}

In this paper, we apply equivariant Gromov--Witten  theory to obtain new structural results on closed symplectic manifolds with sufficiently nice (\emph{Hamiltonian GKM}) torus actions.
The axiomatic behaviour of Gromov--Witten invariants serves as a substitute for techniques from algebraic geometry, allowing us to obtain results in Hamiltonian geometry which have their classical origins in algebraic geometry.

For example, the \emph{Mukai conjecture} \cite{Mukai1988Problems} states that any Fano manifold $X$ of complex dimension $n$ and Picard number $\mathbf{p}(X)$ satisfies
\[
\mathbf{p}(X)(\find(X) - 1) \le n.
\]
Here, $\find(X)$ is the \emph{Fano index} of $X$, defined to be the largest integer dividing $c_1(TX)$ in $H^2(X;\ZZ)$.
Similarly, the \emph{generalised Mukai conjecture} \cite{MR1998396} asserts
\[
\mathbf{p}(X)(\psind(X) - 1) \le n,
\]
where $\psind(X)$ is the \emph{pseudoindex} of $X$, defined as the minimum degree of the anticanonical line bundle $-K_X$ on any rational curve $C$ in $X$.
In both cases, equality is conjectured to hold precisely for $(\CC\PP^{\find(X) - 1})^{\mathbf{p}(X)}$ and $(\CC\PP^{\psind(X) - 1})^{\mathbf{p}(X)}$, respectively.

Both conjectures were generalised to positive monotone\footnote{Here, \emph{positive monotone} means $[\omega]=c_1(TX)$.} symplectic manifolds $(X,\omega)$ in \cite{MR4671382}.
Here, the Picard rank $\mathbf{p}(X)$ is replaced by the second Betti number $b_2(X)$ (as they agree for Fano manifolds), the \emph{index} $\find(X,\omega)$ is defined as before (up to torsion), and the \emph{pseudoindex} $\psind(X,\omega)$ is redefined\footnote{
    $\find(X,\omega)$ and $\psind(X,\omega)$ depend on $\omega$ through $c_1(TX)$, which is well-defined by picking any $\omega$-compatible almost complex structure on $X$, cf.\ \cite[\S2.7]{MR3674984}.
} following \cite{MR4029085} as the infimum of $\int_C c_1(TX)$ where $C$ ranges over all embedded symplectic $S^2$ in $X$.
The \emph{symplectic Mukai inequality} and its generalisation are
\[
b_2(X)(\find(X,\omega) - 1) \le n
\hspace{7mm}
\text{and}
\hspace{7mm}
b_2(X)(\psind(X,\omega) - 1)\le n,
\]
respectively, where $2n$ is the dimension of $X$.
While important special cases of the (generalised) Mukai conjecture have been proved (e.g., \cite{MR1063222, MR2228683, MR2761927, MR2660578, MR2725271, MR2997473, MR3592522, MR4790927, MR5048523, Gagliardi_Hofscheier_Pearson_2025, Pearson_2026}), its general form remains open --- up to human verification of the AI generated proof claimed by
\cite{OAI:The-Generalized-Mukai-Conjecture-September-24-2026}, which appeared during the final stage of preparing this article and which the author could not fully follow.
In the symplectic case, much less is known beyond the initial paper \cite{MR4671382}.

\textbf{Spaces.} We consider here the following spaces.
Let $(X,\omega)$ be a closed connected symplectic manifold of dimension $2n$ and let $T=(S^1)^r$ be a compact torus with an effective Hamiltonian action on $(X,\omega)$.
Then $(X,\omega,T)$ is a \emph{Hamiltonian GKM space} (named after \cite{GKM98}) if there are only finitely many $T$-fixed points and the isotropy weights at each fixed point are pairwise linearly independent.
Examples of Hamiltonian GKM spaces include closed toric symplectic manifolds, homogeneous spaces \cite{Guillemin_Holm_Zara_2006}, smooth Schubert varieties, certain Hessenberg varieties \cite{Goldin_Tymoczko_2023}, and even non-equivariantly-K{\"a}hler examples like those constructed by Tolman and Woodward \cite{Tolman_1998, Woodward_1998} (cf.\ \cite{MR4088417}).

\textbf{Tools.} A key tool to study a GKM space is its \emph{GKM graph} $\Gfat=(\G,\alpha)$, where $\G$ is the simple $n$-regular graph whose vertices are the $T$-fixed points of $X$ and whose edges are the $T$-stable {2-spheres} of $X$.
The \emph{axial function} $\alpha$ decorates each edge with the $T$-weight associated to its 2-sphere (see \S\ref{sec : gkm} for details).
Like in the toric case, $\Gfat$ encodes much of the topology of $(X,\omega,T)$, including its Betti numbers, Chern numbers, equivariant cohomology, and $K$-theory \cite{GKM98, GZ01, Guillemin_Zara_2002, Rosu_03}.
In particular, we can access the quantities appearing in the symplectic Mukai inequality through $\Gfat$.
On the other hand, $\Gfat$ does not determine the homotopy type of $X$ in dimension greater than 6 \cite{MR4545925}.

\textbf{Methods.} Our main technique is to compute equivariant Gromov--Witten invariants in terms of $\Gfat$ \cite{LS17, Hirschi_2023, Holmes_Muratore_2025b}.
Any non-vanishing invariant in curve class $\beta\in H_2(X;\ZZ)$ produces an upper bound on $\int_\beta c_1(TX)$, giving rise to our technical results listed in \S\ref{sec : intro technical}.

\textbf{Symplectic vs.\ algebraic.} Although we work in symplectic generality, all our results  hold equally well for algebraic GKM spaces $(X,T_\CC)$ after replacing $\psind(X,\omega)$ by $\psind(X)$
(which implicitly depends on the complex structure of $X$)
and using $b_2(X) =\mathbf{p}(X)$ in this case.
In particular, whenever we prove a special case of the generalised symplectic Mukai inequality, the same argument proves the analogous case of the classical (algebraic) generalised Mukai conjecture.
We give details on this in \S\ref{sec : alg Mukai}.

\subsection{Main results}

First, we prove a large-pseudoindex special case of the generalised symplectic Mukai inequality, whose analogue in algebraic geometry is due to Wi{\'s}niewski \cite{MR1063222}.
Note that although the Mukai conjecture was classically formulated for Fano manifolds \cite{Mukai1988Problems, MR1998396} and monotone symplectic manifolds \cite{MR4671382}, our results do not require these assumptions.

\begin{theorem}[Theorem \ref{thm : Mukai large index}]\label{thm : intro large index}
    Let $(X,\omega,T)$ be a Hamiltonian GKM space of dimension $2n>0$ with GKM graph $\Gfat$.
    \begin{enumerate}
        \item If $\psind(X,\omega) > (n+2)/2$ then $b_2(X)=1$.
        \item If $\psind(X,\omega) = (n+2)/2$ then either
        $b_2(X)=1$; or $b_2(X)=2$ and $\Gfat$ is isomorphic to a projection of the\footnote{Here, we mean the standard GKM action of Example~\ref{ex: GKM graph of CPn} on $\CC\PP^{n/2}$, see also Example~\ref{ex: CPn proj criterion}.} GKM graph of $(\CC\PP^{n/2})^2$.
    \end{enumerate}
\end{theorem}

Second, we introduce the notion of a \emph{special vertex} on $\Gfat$ (see Definition~\ref{def : special vertex}) and prove the generalised symplectic Mukai inequality for any Hamiltonian GKM space that has one.
Informally, a special vertex is one whose \emph{local} cone of curves already is the \emph{full} cone of curves inside $H_2(X;\RR)$.

\begin{theorem}[Corollary \ref{cor : key vertex}]\label{thm : intro Mukai}
    Let $(X,\omega,T)$ be a Hamiltonian GKM space of dimension $2n>0$ with GKM graph $\Gfat$.
    If $\Gfat$ has a special vertex, then
    \[
    b_2(X)(\psind(X,\omega) -1) \le n.
    \]
\end{theorem}

While not every Hamiltonian GKM space has a special vertex (see \S\ref{sec : no special vertex}), every homogeneous space and every smooth Schubert variety has one (see \S\ref{sec : examples with special vertex}).
This gives a new proof of the generalised Mukai inequality for homogeneous spaces, which was previously proved by \cite{MR2415341, MR4671382}.
For smooth Schubert varieties, however, the following result appears to be new to the best of the author's knowledge (up to verification of \cite{OAI:The-Generalized-Mukai-Conjecture-September-24-2026}).

\begin{corollary}[Corollary \ref{cor : mukai for schubert}]\label{cor : intro Mukai Schubert}
    Let $X_{[w]}\subset G/P$ be a smooth Schubert variety, where $G$ is a complex semisimple connected Lie group and $P\subset G$ a parabolic.
    Then the generalised Mukai inequality holds for $X_{[w]}$.
\end{corollary}

\subsection{Corollaries and polytopes}
Before presenting our technical results in \S\ref{sec : intro technical}, we highlight some corollaries that are of independent interest.
For each edge $e$ of the GKM graph, let $C_e\subset X$ be its $T$-invariant $2$-sphere, which is a symplectic submanifold by \cite[Lemma~5.5.7]{MR3674984}.
Following the notation of \cite{Godinho_Sabatini_2014}, the \emph{magnitude of $e$} is
\[
\m(e):=\int_{C_e}c_1(TX).
\]
The following is immediate from Theorem~\ref{thm: connectivity restated} and reminiscent of Campana and Koll{\'a}r--Miyaoka--Mori's rational chain connectedness of Fano manifolds \cite{MR1191735, MR1189503}, cf. \cite[Proof of Corollary~II.2.14]{MR1440180}.
It is different from the recent connectivity result in \cite{GKZ2026}, which concerns the removal of any single vertex or edge.

\begin{corollary}\label{cor : connectedness}
Let $\Gfat$ be the GKM graph of a Hamiltonian GKM space of dimension $2n$.
Then $\G$ remains connected after removing all edges of magnitude greater than $n+1$.
\end{corollary}

Next, we obtain the upper bounds $\find(X,\omega),\psind(X,\omega) \le n+1$ and present Hamiltonian GKM analogues of the classical Kobayashi--Ochiai characterisation of projective space and complex quadrics \cite{MR316745, MR1929792, MR1922103, MR3098425, MR3592468}.
These correspond precisely to the extremal cases of (pseudo)index $n+1$ and $n$, respectively.

The upper bound is well-known in many geometric settings, including \cite{MR316745}, \cite{MR554387}, \cite[Corollary~0.3]{MR1929792} for Fano manifolds, \cite[Corollary~7.17]{MR595007} for the index of compact K{\"a}hler manifolds with positive Ricci curvature,
and \cite[Theorem~5.1]{MR776689} for the index of compact almost complex $S^1$-manifolds with quasi-ample line bundle and isolated fixed points.
For Hamiltonian $S^1$-manifolds with isolated fixed points (which includes Hamiltonian GKM spaces), the bound $\find(X,\omega)\le n+1$ was proved in \cite[Corollary~1.3]{MR3665359_Sabatini_2017}.
We present a different proof for Hamiltonian GKM spaces using Theorem~\ref{thm: unique min path restated} in Remark~\ref{rem : deduce find le n+1}.
Similarly, the bound on $\psind(X,\omega)$ and the statement for $\psind(X,\omega)=n+1$ were proved for Hamiltonian $S^1$-manifolds with isolated fixed points in \cite{MR4029085} assuming unimodality of the even-index Betti numbers and existence of a toric 1-skeleton.
In the Hamiltonian GKM case, our result extends the latter as we do not require unimodality.

\begin{corollary}\label{cor : index}
    Let $(X,\omega,T)$ be a Hamiltonian GKM space of dimension $2n>0$ with GKM graph $\Gfat$.
    Then
    \[
    \psind(X,\omega)\le n+1.
    \]
    Moreover\footnote{In the following, we mean the GKM graph of the standard $(\CC\PP^n, \omega_{FS}, T)$ as in Example~\ref{ex: GKM graph of CPn}, see also Example~\ref{ex: CPn proj criterion}.
    The complex quadric $\cQ_n$ and its GKM graph are recalled in \S\ref{sec : quadrics}.},
    \begin{align*}
        \psind(X,\omega)=n + 1
        &
        \Longrightarrow
        \begin{cases}
            \Gfat \text{ is isomorphic to a projection of the GKM graph}
            \\
            \text{ of $\CC\PP^n$ and $X$ is homotopy equivalent to }\CC\PP^n.
        \end{cases}
        \\
        \psind(X,\omega)=n &\Longrightarrow 
        \begin{cases}
            \Gfat\text{ is isomorphic to a projection of the}
            \\
            \text{GKM graph of the complex quadric }\cQ_n.
        \end{cases}
    \end{align*}
    In the two cases above, $X$ is complex cobordant to $\CC\PP^n$ and $\cQ_n$, respectively.
\end{corollary}

The upper bound and its equality case are proved as Theorem~\ref{thm: index n+1 with proof} using Corollary~\ref{cor : connectedness},  Theorem~\ref{thm: local Mukai} and \cite[Theorem~1.7]{MR4029085}.
The case $\psind(X,\omega)=n$ is proved as Theorem~\ref{thm : quadrics} in \S\ref{sec : quadrics}.
The complex cobordism statement follows immediately since the GKM graph determines all Chern numbers of $X$.

Note that the cases $\psind(X,\omega)\in\{n,n+1\}$ of Corollary~\ref{cor : index} are proved for almost complex $S^1$-manifolds with isolated fixed points in \cite[Theorems~5.7, 6.1, 6.17]{MR776689} using the index instead of the pseudoindex, and assuming that $\chi(X)=\chi(\CC\PP^n)$ and $\chi(X)=\chi(\cQ_n)$, respectively.
However, our proofs only determine $\chi(X)$ at their final steps and are therefore not easily shortened by appealing to Hattori's results.

Finally, we translate Corollary~\ref{cor : connectedness} and Theorem~\ref{thm: local Mukai} to statements about \emph{reflexive Delzant polytopes} $\Delta$ (see \S\ref{sec : polytopes}) using the correspondence between reflexive Delzant polytopes and closed symplectic toric manifolds $(X,\omega)$ with $[\omega]=c_1(TX)$ (cf.\ \cite[\S3.3]{MR3695881}).

\begin{corollary}[Corollary \ref{cor : polytopes}]\label{cor : intro polytopes}
    Let $\Delta\subset\RR^n$ be a reflexive Delzant polytope. Then
    \begin{enumerate}
        \item The $1$-skeleton of $\Delta$ remains connected after removing all edges of affine length greater than $n+1$.
        \item Let $v$ be any vertex of $\Delta$ and let $L_{\min}(v)$ be the minimum affine length of any edge of $\Delta$ incident to $v$.
        Then there exist at least $L_{\min}(v) - 1$ edges of affine length $L_{\min}(v)$ at $v$.
    \end{enumerate}
\end{corollary}

We provide background on reflexive Delzant polytopes and deduce Corollary~\ref{cor : intro polytopes} in \S\ref{sec : polytopes}.

\subsection{Overview of technical results}
\label{sec : intro technical}

Tables~\ref{tab: magnitude bounds} and \ref{tab: line chern bounds} give an overview of our technical results using the following notation.
Given a GKM graph $\Gfat=(\G,\alpha)$, let $E(\G)$ and $V(\G)$ be its set of edges and vertices, respectively (see Definition~\ref{def : gkm graph}).
For $v\in V(\G)$, let $E(\G)_v$ be the set of edges incident to $v$.

\begin{definition}\label{def : paths and subpaths}
    A \emph{path $P$ of length $l>0$} on the GKM graph $\Gfat$ is a sequence $v_0\dots v_l\in V(\G)$ such that $v_{i-1}v_{i}\in E(\G)$ for each $i\in\{1,\dots,l\}$.
We say that $P$ is a path \emph{from $v_0$ to $v_l$}.
The \emph{length of $P$} is $\len(P):=l$.
The \emph{total magnitude of $P$} is
$$\m(P):=\sum_{i=1}^{\len(P)}\m(v_{i-1}v_i)$$
and the \emph{class (or curve class) of $P$} is 
$$[P]:=\sum_{i=1}^{\len(P)} [C_{v_{i-1}v_i}]\in H_2(X;\ZZ).$$

A \emph{subpath of $P$} is a path of the shape $v_j\dots v_k$ for any $0\le j<k\le l$.
In particular, $P$ is a subpath of itself, and any edge $v_{i-1}v_i$ of $P$ is a subpath of $P$.
\end{definition}

\begin{definition}\label{def : partial order on H2}
    Let $<$ be the strict partial order on $H_2(X;\ZZ)$ generated by the following relations: for any $B,C\in H_2(X;\ZZ)$ and $e\in E(\G)$ with
\[
    [C_e] + B = C
\]
we have
\[
    B < C.
\]
Note that $<$ is well-defined (irreflexive and asymmetric) 
because each $C_e$ is a symplectic submanifold of $(X,\omega)$, so $\omega(B)<\omega(C)$ whenever $B<C$.
For any set $S\subset H_2(X;\ZZ)$ and $C\in S$, we call $C$ \emph{minimal among $S$} if there is no $B\in S$ with $B < C$.
As usual, we write $B\le C$ for \emph{"$B<C$ or $B=C$"}.
For a path $P$ from $v$ to $w$, we say that \emph{$[P]$ is minimal among paths from $v$ to $w$} if $[P]$ is minimal among $S:=\{[Q] : Q\text{ is a path from } v\text{ to }w\}$.

\end{definition}

We are now ready to summarise our magnitude bounds in Table~\ref{tab: magnitude bounds}.

\begin{table}[h]
    \centering
    \small
    \setlength{\tabcolsep}{4pt}
    \renewcommand{\arraystretch}{1.25}
    \renewcommand{\tabularxcolumn}[1]{m{#1}}
    \begin{tabularx}{\textwidth}{|
        >{\hsize=0.5\hsize\centering\arraybackslash}X|
        >{\hsize=1.65\hsize\centering\arraybackslash}X|
        >{\hsize=0.85\hsize\centering\arraybackslash}X|
    }
        \hline
        \textbf{Theorem} \cellcolor{tableHead}&
        \textbf{Object} \cellcolor{tableHead}&
        \textbf{Conclusion on magnitudes}\cellcolor{tableHead}
        \\
        \hline\hline\cellcolor{tableRow}
            Theorem~\ref{thm: unique min path restated} &\cellcolor{tableRow}
            Let $v\neq w\in V(\G)$ and let
            $P$ be a path from $v$ to $w$
            such that no other path from $v$ to $w$ has class $[P]$, and $[P]$ is minimal among all paths from $v$ to $w$.
        &\cellcolor{tableRow}
            $\m(P)\le n+1$
        \\
        \hline\cellcolor{tableRowGray}
        Theorem~\ref{thm: connectivity restated} &\cellcolor{tableRowGray}
            Let $v\neq w\in V(\G)$. Then there exists a path $P$ from $v$ to $w$ such that for each subpath $Q$ of $P$:
        &\cellcolor{tableRowGray}
            $\m(Q) \le (n-1)\len(Q) + 2$
        \\
        \hline\cellcolor{tableRow}
        Theorem \ref{thm: local Mukai}
        (\emph{local Mukai})
        &\cellcolor{tableRow}
            Let $v\in V(\G)$ and $e\in E(\G)_v$ such that $[C_e]$ is minimal among curve classes of edges at $v$.
        &\cellcolor{tableRow}
            There exist at least $\m(e)-1$ edges at $v$ with curve class $[C_e]$.
        \\
        \hline\cellcolor{tableRowGray}
        Theorem~\ref{thm: unique minimal edge restated} &\cellcolor{tableRowGray}
            Let $e\in E(\G)$ be the unique edge of $\G$ with curve class $[C_e]$ and assume $[C_e]$ is minimal among the curve classes of all edges.
        &\cellcolor{tableRowGray}
            $\m(e)\le 3-n$.
        \\
        \hline
    \end{tabularx}
    \vspace{1mm}
    \caption{Overview of magnitude bound results.}
    \label{tab: magnitude bounds}
\end{table}

Our results also bound the integers $a_1,\dots,a_{n-1}$ appearing in the normal bundle splittings $N_{C_e/X}\cong \mathcal{O}_{C_e}(a_1)\oplus\cdots\oplus\mathcal{O}_{C_e}(a_{n-1})$ for certain $e\in E(\G)$.
Given $e$, let $\cN(e,a)$ be the number of indices $i\in \{1,\dots,n-1\}$ for which $a_i\ge a$.
The resulting bounds are presented in Table~\ref{tab: line chern bounds}, where the considered objects refer row by row to those of Table~\ref{tab: magnitude bounds}.
In fact, as explained in \S\ref{sec : connection and splittings}, any compatible connection $\nabla$ on $\Gfat$ induces integers $a_1,\dots,a_{n-1}$ for each $e\in E(\G)$, and one can analogously define $\cN^\nabla(e,a)$.
The bounds in Table~\ref{tab: line chern bounds} then hold for any compatible connection $\nabla$.

\begin{table}[h]
    \centering
    \small
    \setlength{\tabcolsep}{4pt}
    \renewcommand{\arraystretch}{1.25}
    \renewcommand{\tabularxcolumn}[1]{m{#1}}
    \begin{tabularx}{\textwidth}{|
        >{\hsize=0.67\hsize\centering\arraybackslash}X|
        >{\hsize=2.15\hsize\centering\arraybackslash}X|
        >{\hsize=0.59\hsize\centering\arraybackslash}X|
        >{\hsize=0.59\hsize\centering\arraybackslash}X|
    }
        \hline
        \cellcolor{tableHead}&\cellcolor{tableHead}&
        \multicolumn{2}{c|}{\textbf{Conclusion on $\mathbf{\cN^\nabla(e,a)}$}\cellcolor{tableHead}}
        \\
        \cline{3-4}
       \multirow{-2}{=}{\centering\arraybackslash\textbf{Theorem}} \cellcolor{tableHead}&
        \multirow{-2}{=}{\centering\arraybackslash\textbf{Object (continuing from the corresponding row of Table~\ref{tab: magnitude bounds})}} \cellcolor{tableHead}& $\GKM_2$ \cellcolor{tableHead}& $\GKM_3$\cellcolor{tableHead}
        \\
        \hline\hline
        Theorem~\ref{thm: unique min path restated}\cellcolor{tableRow}
        &
        For each edge $e$ of $P$:\cellcolor{tableRow}
        &
            $\cN^\nabla(e,2)\le k$\cellcolor{tableRow}
        &
            $\cN^\nabla(e,2)=0$\cellcolor{tableRow}
        \\
        \hline
        Theorem~\ref{thm: connectivity restated} \cellcolor{tableRowGray}&\cellcolor{tableRowGray}
            For each edge $e$ of $P$:
        &\cellcolor{tableRowGray}
            $\cN^\nabla(e,2)\le k$
        &\cellcolor{tableRowGray}
            $\cN^\nabla(e,2)=0$
        \\
        \hline
        Theorem \ref{thm: magnitude at vertex restated} \cellcolor{tableRow}&\cellcolor{tableRow}
            If no other edge at $v$ has curve class $[C_e]$ then:
        &\cellcolor{tableRow}
            $\cN^\nabla(e,1)\le k$
        &\cellcolor{tableRow}
            $\cN^\nabla(e,1)=0$
        \\
        \hline
        Theorem~\ref{thm: unique minimal edge restated} \cellcolor{tableRowGray}&\cellcolor{tableRowGray}
            For the given $e$, we have:
        &\cellcolor{tableRowGray}
            $\cN^\nabla(e,0)\le k$
        &\cellcolor{tableRowGray}
            $\cN^\nabla(e,0)=0$
        \\
        \hline
    \end{tabularx}
    \vspace{1mm}
    \caption{Overview of bounds on $\cN^\nabla(e,a)$,  corresponding to the cases of Table~\ref{tab: magnitude bounds}. Results distinguish between the generic $\GKM_2$ case and the non-generic $\GKM_3$ case (see Definition~\ref{def: k-indep} and Remark~\ref{rem : gkm2 vs gkm3}). Here, $k$ is the complexity of $(X,\omega,T)$ (see Definition~\ref{def : complexity}).}
    \label{tab: line chern bounds}
\end{table}

In \cite{KSZ_moment_maps_in_a_box}, the integers $a_1,\dots,a_{n-1}$ are called \emph{line Chern numbers} and quantitative upper and lower bounds for them in terms of $n$ and $\chi(X)$ are proved in the monotone $\GKM_3$ case.
By contrast, our upper bounds in the $\GKM_3$ case are much stronger and do not require monotonicity, but only apply to certain special edges.

\subsection{Algebraic versions of our results}
\label{sec : alg Mukai}
All results in this paper apply equally well in the algebraic category.
An \emph{algebraic GKM space} $(X,T_\CC)$ is a smooth projective variety $X$ with an algebraic action of the complex torus $T_\CC\cong(\CC^\times)^r$ satisfying the same properties as a Hamiltonian GKM space (see Definition~\ref{def : gkm space}).
The only difference is that no choice of symplectic form $\omega$ is required.
Indeed, when the symplectic version of the theory picks any $T$-invariant $\omega$-compatible almost complex structure $J$ on $X$, the algebraic version simply picks the canonical integrable complex structure of the complex variety $X$.
Note also that restricting the action to $(S^1)^r\subset T_\CC$ produces a Hamiltonian GKM space with the same GKM graph (cf.\ \cite[Remark~2.3]{Holmes_Muratore_2025b}).
Therefore, all arguments using Hamiltonian geometry still apply and it only remains to check the following two details.
\begin{enumerate}
    \item \emph{All arguments work with $\psind(X,\omega)$ replaced by $\psind(X)$.}
    This holds because we only ever access the pseudoindex through the GKM graph, via Lemma~\ref{lem : find and psind from G}.
    Its proof holds completely analogously for $\psind(X)$ as each $C_e\subset X$ is a rational curve in the algebraic setup.

    \item \emph{The Picard number $\mathbf{p}(X)$ agrees with $b_2(X)$.}
    This holds since $X$ is rational by the Białynicki-Birula decomposition \cite{BialynickiBirula1973}
    and hence $H^{2,0}(X)= H^0(X,\Omega_X^2)=0$ which implies $\mathrm{NS}(X)\otimes\QQ\cong H^2(X;\QQ)\cap H^{1,1}(X)\cong H^2(X;\QQ)$ via the Lefschetz theorem on $(1,1)$-classes and Hodge theory.
\end{enumerate}

\subsection{Strengthening the pseudoindex}
The only properties of the pseudoindex $\psind(X,\omega)$ that we ever use are that $\m(e)\ge \psind(X,\omega)$ for each $e\in E(\G)$, and that $\find(X,\omega)\mid \psind(X,\omega)$ whenever $\psind(X,\omega)$ is finite (see Lemma~\ref{lem : find and psind from G}).
The same properties hold for the \emph{equivariant pseudoindex} of \cite[\S1]{MR4029085}, which in our setting is
\[
\psind(\Gfat) := \min\{\m(e) : e\in E(\G)\}.
\]
Therefore, all results in this paper hold after replacing $\psind(X,\omega)$ by $\psind(\Gfat)$ throughout, which leads to potentially stronger results since $\psind(\Gfat)\ge \psind(X,\omega)$.
While the definition of $\psind(X,\omega)$ is intrinsic to $(X,\omega)$, the advantage of using $\psind(\Gfat)$ is that it is readily computable from the GKM graph $\Gfat$ (see \S\ref{sec : magnitudes background}).
This allows our results to be applied even when only $\Gfat$ is known.

\subsection{Structure of the paper}

In \S\ref{sec : gkm}, we recall the required elements of GKM theory.
In \S\ref{sec : GW}, we set up Gromov--Witten theory as a black box and deduce the main lemmas from the localisation formula.
In \S\ref{sec : curve class combinatorics}, we prepare the curve class combinatorics that we will need to identify non-vanishing Gromov--Witten invariants.
In \S\ref{sec : proofs of technical results}, we prove our technical results from \S\ref{sec : intro technical} using the previous sections.
In \S\ref{sec : pseudoindex n+1 and projective space}, we treat the pseudoindex $n+1$ case of projective space.

The main results on the generalised symplectic Mukai inequality are proved in \S\ref{sec: mukai}, while the pseudoindex $n$ case of quadrics is treated in \S\ref{sec : quadrics}.
Finally, our corollaries on reflexive Delzant polytopes are deduced in \S\ref{sec : polytopes}.

\subsection{Acknowledgements}
The author would like to thank Silvia Sabatini for stimulating discussions and previous results that inspired this work, as well as detailed comments that improved this draft; Benjamin Becker, Oliver Goertsches, Panagiotis Konstantis, and Grigory Solomadin for fruitful conversations at the University of Marburg in January 2026; Tam{\'a}s Hausel for supporting this project during the author's PhD studies at ISTA; and Rahul Pandharipande for hosting the author at ETH Z{\"u}rich from November 2025 to January 2026.
Further thanks are due to Mischa Elkner and Jakub L{\"o}wit for helpful comments on aspects of Lie theory; and Giosu{\`e} Muratore for helpful comments and encouragement to write down the case of complex quadrics.
Drawings of GKM graphs were produced using \textsf{GKMtools.jl} \cite{Holmes_Muratore_2025a}.
This work is supported by the ERC grant ViaFiPoS, No 101199663.
From September 2026, the author is a recipient of a DOC Fellowship of the Austrian Academy of Sciences.

\subsection{AI disclosure}
This paper is entirely human-written and the mathematical content is due to the author. The only exception is \S\ref{sec : no special vertex}, where the author used \textsf{ChatGPT} and \textsf{GKMtools.jl} to identify the example of $\mathrm{Bl}_{\PP^2\sqcup\PP^2}\PP^5$, which was subsequently verified and written up by the author.
The \LaTeX\ formatting and some illustrations were optimised using \textsf{Claude Code} with \textsf{Sonnet 4.5/4.6/5.5}, and
\textsf{ChatGPT 6} was used as one layer of the final proofreading.

\section{GKM theory}
\label{sec : gkm}

In this section, we recall the required background on GKM theory.
We refer to \cite{GZ01, Guillemin_Zara_2002}, \cite[\S2]{GKZ22}, and \cite[\S2]{Charton_Kessler_2023} for more details.

\subsection{Hamiltonian GKM spaces and their GKM graphs}

Given a compact real torus $T=(S^1)^r$, we write $\Lie(T)$ for its Lie algebra, $\Lie(T)^\vee$ for its dual, and 
$\CC_\alpha$ for a complex $T$-representation on $\CC$ with infinitesimal weight $\alpha\in\Lie(T)^\vee$.

\label{sec : gkm spaces and their graphs}

\begin{definition}\label{def : gkm space}
    A \emph{Hamiltonian GKM space} $(X,\omega,T)$ is a closed connected symplectic manifold $(X,\omega)$ of dimension $2n$ with an effective Hamiltonian action by the torus $T = (S^1)^r$, such that:
    \begin{enumerate}
        \item The fixed point set $X^T$ is finite.
        \item At each fixed point $p\in X^T$, the tangent weights $\alpha_1,\dots,\alpha_n\in\Lie(T)^\vee$ in the splitting $T_pX = \bigoplus_{i=1}^n\CC_{\alpha_i}$ are pairwise linearly independent.
        Note that the tangent weights are determined by the choice of any $T$-invariant $\omega$-compatible almost complex structure $J$ on $X$. Since the set of such $J$ is contractible, all choices of $J$ lead to the same weights.
    \end{enumerate}
\end{definition}

For a given Hamiltonian torus action on $(X,\omega)$ with $\dim X > 0$, being GKM is equivalent to $\# X^T <\infty$ and $\dim X_1 =2$, where
\[
X_1 := \{x\in X : \dim(T\cdot x) \le 1\}
\]
is the \emph{1-skeleton of $(X,T)$}
\cite[Proposition~1.1]{MR2064318_Guillemin_Holm_04}.
Under this condition, $X_1$ is a union of transversely intersecting $2$-spheres.

The effectiveness assumption means that no element of $T$ other than the identity acts trivially on $X$.
We make this assumption without loss of generality, to simplify Definition~\ref{def : complexity} and the statements of Theorems~\ref{thm : intro large index} and Corollary~\ref{cor : index}.
Throughout, let $(X,\omega,T)$ be a Hamiltonian GKM space of dimension $2n$.

\textit{Graph-theoretic notation:}
Given a simple undirected graph $\G$, we denote its vertex set by $V(\G)$ and its edge set by $E(\G)$.
We write $E(\G)^\pm$ for the set of edges of $\G$ with a choice of vertex, called \emph{oriented edges}.
    That is, the forgetful map
$$\pi\colon E(\G)^\pm\rightarrow E(\G)$$
is two-to-one.
Given $\epsilon\in E(\G)^\pm$, we write $\src(\epsilon)$ for its chosen vertex (the \emph{source}) and $\dst(\epsilon)$ for its other vertex (the \emph{destination}).
For vertices $v,w$ defining an edge, we write $\{v,w\}\in E(\G)$ for the unoriented edge and $vw\in E(\G)^\pm$ for the oriented edge with $\src(vw)=v$ and $\dst(vw)=w$.
We define
\[
E(\G)^\pm_v=\{\epsilon\in E(\G)^\pm : \src(\epsilon)=v\}\,\,\,\mathrm{and}\,\,\,
E(\G)_v=\{e\in E(\G) : v\in e\}.
\]

\begin{definition}\label{def : gkm graph}
    The \emph{GKM graph} $\Gfat=(\G,\alpha)$ of $(X,\omega,T)$ is given by:
    \begin{enumerate}
        \item The undirected simple graph $\G$ whose vertices $V(\G)$ are the fixed points $X^T$ and whose edges are given by the $T$-invariant 2-spheres in $X$, which are symplectic by \cite[Lemma~5.5.7]{MR3674984}. Each $T$-invariant 2-sphere contains precisely two distinct fixed points, which are the vertices joined by the corresponding edge.

        \textit{Notation:} For $e\in E(\G)$, we write $C_e\subset X$ for the associated $T$-invariant 2-sphere and $[C_e]\in H_2(X;\ZZ)$ for its \emph{class} (or \emph{curve class}).
    \item The \emph{axial function}
\[
\alpha\colon E(\G)^\pm\longrightarrow \Lie(T)^\vee
\]
sends $\epsilon\in E(\G)^\pm$ to the $T$-weight of $T_{\src(\epsilon)}C_{\pi(\epsilon)}$.
    \end{enumerate}
\end{definition}

The 1-skeleton $X_1$ of $X$ is precisely $\cup_{e\in E(\G)} C_e$, and there are only finitely many $T$-stable $2$-spheres in $X$ because the isotropy weights at each fixed point are pairwise linearly independent.
The axial function satisfies $\alpha(\overline{\epsilon})=-\alpha(\epsilon)$, where $\overline{\epsilon}$ is the other preimage of $\pi(\epsilon)$ under $\pi$.
By the Hamiltonian assumption, $\G$ is simple \cite[Remark~2.7]{Charton_Kessler_2023} and connected \cite[Remark~2.9]{Charton_Kessler_2023}.
By \cite[\S1.5]{Guillemin_Zara_2002}, each vertex of $\G$ has valency $n$.

\begin{example}
\label{ex: GKM graph of CPn}
    On $\CC\PP^n$ with Fubini--Study form $\omega_{FS}$, the action of $T=(S^1)^n$ given by
    \[
    (\lambda_1\dots,\lambda_n)\cdot [x_0:x_1:\cdots:x_n] = [x_0:\lambda_1x_1:\cdots:\lambda_n x_n]
    \]
    is Hamiltonian (see \cite[Example~5.5.2]{MR3674984}) and GKM.
    The fixed points are
    \[
    (\CC\PP^n)^T=\{p_i : i\in\{0,1,\dots,n\}\}\hspace{5mm}\text{where}\hspace{5mm}
    p_i :=[0:\cdots:0:\underbrace{1}_{i^\text{th}}:0:\cdots:0]\in \CC\PP^n.
    \]
    The $T$-invariant 2-spheres are
    \[
    \{C_{\{i,j\}} : 0\le i < j \le n\}
    \hspace{5mm}\text{where}\hspace{5mm}
    C_{\{i,j\}} := \{[x_0:\cdots:x_n ] : x_k=0\,\,\,\forall\,k\neq i,j\}.
    \]
    Thus, the GKM graph is $\Gfat_{\CC\PP^n}=(\G,\alpha)$ where $\G$ is the complete graph on the vertex set $\{p_0,\dots,p_n\}$ and the axial function is
    \[
    \alpha(p_ip_j) = t_j - t_i
    \]
    where $t_1,\dots,t_n\in \Lie(T)^\vee$ correspond to $\lambda_1,\dots,\lambda_n$, respectively, and $t_0:=0$.
\end{example}

\begin{definition}\label{def: k-indep}
    A GKM space $(X,\omega,T)$ is called $\GKM_k$ or \emph{$k$-independent} for $k\in\ZZ_{\ge 1}$ if for any $v\in V(\G)$ and any distinct $\epsilon_1,\dots,\epsilon_k\in E(\G)^\pm _v$, the elements $\alpha(\epsilon_1),\dots,\alpha(\epsilon_k)\in \Lie(T)^\vee$ are linearly independent.
\end{definition}

\begin{remark}\label{rem : gkm2 vs gkm3}
    By definition, every GKM space is $\GKM_2$, but there are examples which are not $\GKM_3$.
\end{remark}

\begin{definition}\label{def : iso and projection of GKM graphs}
Let $\Gfat_1=(\G_1,\alpha_1)$ and $\Gfat_2=(\G_2,\alpha_2)$ be the GKM graphs of Hamiltonian GKM spaces $(X_1,\omega_1,T_1)$ and $(X_2,\omega_2,T_2)$, respectively.
\begin{enumerate}
    \item $\Gfat_2$ is \emph{isomorphic to $\Gfat_1$} if there exists a graph isomorphism $\varphi\colon \G_1\rightarrow \G_2$ and a linear isomorphism $\psi\colon \Lie(T_1)^\vee \rightarrow \Lie(T_2)^\vee$ such that each $vw\in E(\G_1)^\pm$ satisfies
    \[
    \alpha_2(\varphi(v)\varphi(w))
    =
    \psi(\alpha_1(vw)).
    \]
    \item $\Gfat_2$ is a \emph{projection of $\Gfat_1$} if $\G_1=\G_2$ and there is a linear surjection $\psi\colon \Lie(T_1)^\vee \rightarrow \Lie(T_2)^\vee$ such that each $vw\in E(\G_1)^\pm$ satisfies
    \[
    \alpha_2(vw)
    =
    \psi(\alpha_1(vw)).
    \]
\end{enumerate}
\end{definition}

\begin{remark}
    By \cite[Lemma~2.1]{Charton_Kessler_2023}, Definition~\ref{def : iso and projection of GKM graphs} (1) and (2) are equivalent to the same conditions for $\Gfat_1$ and $\Gfat_2$ as \emph{abstract GKM graphs} in the sense of \cite[\S2.4--2.5]{Charton_Kessler_2023}.
    If $X_1=X_2$ and $T_2\subset T_1$ is a subtorus whose action is still GKM, then $\Gfat_2$ is a projection of $\Gfat_1$ and $\psi\colon \Lie(T_1)^\vee\rightarrow \Lie(T_2)^\vee$ is induced by $T_2\hookrightarrow T_1$ \cite[Lemma~2.19]{Charton_Kessler_2023}.
\end{remark}

\begin{example}\label{ex: CPn proj criterion}
    Recall the GKM graph $\Gfat_{\CC\PP^n}$ from Example~\ref{ex: GKM graph of CPn}.
    One easily checks that a GKM graph $\Gfat=(\G,\alpha)$ of a Hamiltonian GKM space $(X,\omega,T)$ is isomorphic to a projection of $\Gfat_{\CC\PP^n}$ if and only if $\G$ is the complete graph on $n+1$ vertices, and any distinct $u,v,w\in V(\G)$ satisfy
    \[
    \alpha(uv) + \alpha(vw) = \alpha(uw).
    \]
    Note that we used effectiveness of the $T$-action on $X$ here, which implies that the values of the axial function at each vertex of $\G$ span $\Lie(T)^\vee$.
\end{example}

\subsection{Complexity of the action}

\begin{definition}\label{def : complexity}
    The \emph{complexity} $k$ of $(X,\omega,T)$ is given by $k:= n - r$, where $2n$ is the dimension of $X$ and $r$ is the rank of $T$.
\end{definition}

\begin{remark}
    We have $k\ge 0$  since the $T$-action on $X$ is effective and Hamiltonian by Definition~\ref{def : gkm space}.
    The case $k=0$ corresponds to toric symplectic manifolds.
\end{remark}

\subsection{Magnitudes}\label{sec : magnitudes background}

We adopt the following terminology from \cite{Godinho_Sabatini_2014}.

\begin{definition}
    The \emph{magnitude map} $\m\colon E(\G)\rightarrow \ZZ$ is defined by
    \[
    \m(e) = \int_{C_e} c_1(TX).
    \]
\end{definition}
By equivariant localisation (cf.\ \cite[Remark~4.4]{Godinho_Sabatini_2014}), these numbers are easily computed from $\Gfat$ via
\begin{equation}\label{eqn : magnitude from G}
\m(\{v,w\}) = \frac{\sum_{\epsilon_i\in E(\G)_v^\pm} \alpha(\epsilon_i) - \sum_{\epsilon_i'\in E(\G)_w^\pm} \alpha(\epsilon_i')}{\alpha(vw)}.
\end{equation}
Note that swapping $v$ and $w$ only changes the sign of the numerator and denominator, so the right hand side does not depend on a choice of orientation for $\{v,w\}$.

We also adopt the terms \emph{positive} from \cite[\S1.1]{Charton_Kessler_2023} and \emph{almost positive} from \cite[Definition~5.8]{Holmes_Muratore_2025b}.

\begin{definition}\label{def : pos and almost positive}
    A Hamiltonian GKM space is \emph{positive} if $\m(e)>0$ for all $e\in E(\G)$.
    It is \emph{almost positive} if there exists $e_0\in E(\G)$ with $\m(e_0)=0$ and $\m(e)>0$ for all $e\neq e_0$.
\end{definition}

\begin{lemma}\label{lem : exists positive magnitude edge}
    Let $(X,\omega,T)$ be a Hamiltonian GKM space of positive dimension.
    Then there exists at least one edge $e\in E(\G)$ with $\m(e)>0$.
\end{lemma}
\begin{proof}
    Assume for contradiction that $\m(e)\le 0$ for all $e\in E(\G)$ and let $\nu\colon X\rightarrow \Lie(T)^\vee$ be the moment map of $(X,\omega,T)$.
    By \cite[\S2.2]{Guillemin_Guillemin_Ohsawa_Ginzburg_Karshon_2014} there exists some $\lambda_\epsilon > 0$ for each $\epsilon=vw\in E(\G)^\pm$ such that
    \begin{equation}\label{eqn : alpha positive multiple of nu}
    \nu(w)-\nu(v)=\lambda_{\epsilon}\alpha(vw).
    \end{equation}
    By the convexity theorem \cite{Ati82, MR664117} and connectedness of $X$, there exists $\xi\in\Lie(T)$ such that $\langle\xi,\alpha(\epsilon)\rangle\neq 0$ for all $e\in E(\G)$ and such that $\langle\xi,\nu\rangle$ has a unique minimum and maximum on $X$.
    They are automatically $T$-fixed points, corresponding to vertices $v_{\min},v_{\max}\in V(\G)$.
    Note that $v_{\min}\neq v_{\max}$ as $\dim(X)>0$.
    Moreover, there exists a path $v_0\dots v_\len$ on $\G$ with $v_0=v_{\min}$ and $v_\len=v_{\max}$, such that $\langle\xi,\alpha(v_{i-1}v_i)\rangle > 0$ for each $i\in\{1,\dots,\len\}$.
    This follows from the uniqueness of $v_{\min},v_{\max}$ and equation \eqref{eqn : alpha positive multiple of nu}.

    For every $\epsilon=vw\in E(\G)^\pm$ with $\langle\xi,\alpha(\epsilon)\rangle>0$, equation \eqref{eqn : magnitude from G} gives
    \begin{equation}\label{eqn : c1 change for m le 0}
        \left\langle\xi, \sum_{\epsilon_i\in E(\G)^\pm_v} \alpha(\epsilon_i)\right\rangle
        \le
        \left\langle\xi, \sum_{\epsilon_i'\in E(\G)^\pm_w} \alpha(\epsilon_i')\right\rangle.
    \end{equation}
    
    Thus, equation~\eqref{eqn : c1 change for m le 0} applied to each edge of the path implies
    \[
        \left\langle\xi, \sum_{\epsilon_i\in E(\G)^\pm_{v_{\min}}} \alpha(\epsilon_i)\right\rangle
        \le
        \left\langle\xi, \sum_{\epsilon_i'\in E(\G)^\pm_{v_{\max}}} \alpha(\epsilon_i')\right\rangle,
    \]
    while equation~\eqref{eqn : alpha positive multiple of nu} implies
    \[
        \left\langle\xi, \sum_{\epsilon_i\in E(\G)^\pm_{v_{\min}}} \alpha(\epsilon_i)\right\rangle
        > 0
    \]
    and
    \[
        \left\langle\xi, \sum_{\epsilon_i\in E(\G)^\pm_{v_{\max}}} \alpha(\epsilon_i)\right\rangle
        < 0,
    \]
    a contradiction.
    Hence, some $e\in E(\G)$ must satisfy $\m(e)>0$.
\end{proof}

\begin{remark}
    The above proof is similar to \cite[Lemma~5.2]{MR3695881}, which shows that $c_1(TX) = \lambda[\omega]$ with $\lambda\in \RR$ implies $\lambda > 0$ for any Hamiltonian $S^1$-space $X$.
\end{remark}

\subsection{Index, pseudoindex, and irreducible edges}\label{sec : indices}

From now on, we assume that $(X,\omega,T)$ is a Hamiltonian GKM space of positive dimension, which only excludes the trivial example of a point.
The following is well-defined since $c_1(TX)\neq 0$ by Lemma~\ref{lem : exists positive magnitude edge}, and since $H^2(X;\ZZ)$ is torsion-free as $X$ is Hamiltonian with isolated fixed points\footnote{This follows since the moment map evaluated at generic $\xi\in \Lie(T)$ is a Morse function whose critical points all have even index, see \cite{Ati82, MR766741} and cf.\ \cite[Theorem~4.5(i)]{Li_2017}. 
Thus, $H^*(X;\ZZ)$ and $H_*(X;\ZZ)$ are torsion-free and supported in even degree.}\label{footnote : torsion}.

\begin{definition}
    The \emph{index} $\find(X,\omega)\in\ZZ_{\ge 1}$ is the greatest integer dividing $c_1(TX)$ in $H^2(X;\ZZ)$.
\end{definition}

\begin{definition}
    Following \cite{MR4029085}, the \emph{pseudoindex} $\psind(X,\omega)\in\ZZ\cup\{-\infty\}$ is the infimum of $\int_{C}c_1(TX)$ where $C$ ranges over all embedded symplectic $2$-spheres in $X$.
\end{definition}
Both $\find(X,\omega)$ and $\psind(X,\omega)$ depend on $\omega$ through $c_1(TX)$, which is well-defined in terms of any $\omega$-compatible almost complex structure $J$ on $X$ \cite[\S2.7]{MR3674984}.

\begin{lemma}\label{lem : find and psind from G}
    Let $(X,\omega,T)$ be a Hamiltonian GKM space.
    Then
    \[
        \find(X,\omega) = \gcd\{\m(e) : e\in E(\G)\}
    \]
    and
    \[
        \psind(X,\omega) \le \min\{\m(e) : e\in E(\G)\}.
    \]
    Moreover, $H_2(X;\ZZ)$ is a free $\ZZ$-module spanned by $\{[C_e] : e\in E(\G)\}$, and $\find(X,\omega)$ divides $\psind(X,\omega)$ whenever it is finite.
\end{lemma}
\begin{proof}
    By \cite[Remark~3.13]{MR3665359_Sabatini_2017}, $\find(X,\omega)$ agrees with the \emph{minimum Chern number}, which is the non-negative integer $N$ such that $\langle c_1(TX), H_2(X;\ZZ)\rangle=N\ZZ$.
    Moreover, evaluating the moment map $\nu$ of $(X,\omega,T)$ on any generic $\xi\in \Lie(T)$ (where generic means  $\langle\xi,\alpha(\epsilon)\rangle\neq 0$ for all $\epsilon\in E(\G)^\pm$) produces a Morse function with critical points $V(\G)$, where the index of $v\in V(\G)$ is
    \[
    2\cdot \#\{\epsilon\in E(\G)_v^\pm : \langle\xi,\alpha(\epsilon)\rangle<0\},
    \]
    cf.\ \cite[Theorem~1.3.2]{GZ01}.
    In particular, $H_2(X;\ZZ)$ is generated freely over $\ZZ$ by the classes of downward flows of vertices of index 2, which are all of the form $C_e$ for some $e\in E(\G)$.
    Thus, $H_2(X;\ZZ)=\text{span}_\ZZ\{[C_e] : e\in E(\G)\}$ is free and the claim on $\find(X,\omega)$ follows.

    The claim on $\psind(X,\omega)$ holds because each $C_e$ is a symplectic 2-sphere embedded in $X$.

    The final divisibility statement holds by definition, since $c_1(TX)/\find(X,\omega)\in H^2(X;\ZZ)$.
\end{proof}

\begin{definition}\label{def : irreducible}
    An edge $e\in E(\G)$ is \emph{irreducible} if for any $e_1,\dots,e_s\in E(\G)$, the equation
    \[
    [C_e] = \sum_{i=1}^s [C_{e_i}]
    \]
    in $H_2(X;\ZZ)$ implies $s=1$.
\end{definition}

Note that every curve class $[C_e]$ is a sum of curve classes of irreducible edges.

\begin{remark}\label{rem : exists irred edge with positive magnitude}
    Any edge $e\in E(\G)$ whose class $[C_e]$ is minimal among $\{[C_e] : e\in E(\G)\}$ (see Definition~\ref{def : partial order on H2}) is irreducible.
Thus, there exists at least one irreducible edge.
In fact, there exists an irreducible edge $e$ with $\m(e)>0$, since otherwise we would have $\m(e)\le0$ for all $e\in E(\G)$, contradicting Lemma~\ref{lem : exists positive magnitude edge}.
\end{remark}

\begin{remark}\label{rem : deduce find le n+1}
Any irreducible edge is a path of length one that satisfies the hypothesis of Theorem~\ref{thm: unique min path restated} (see Table~\ref{tab: magnitude bounds}).
Hence, this theorem will imply $\m(e)\le n+1$ for any irreducible edge $e$.
Combining this with Remark~\ref{rem : exists irred edge with positive magnitude}, there exists an irreducible edge $e$ with $1\le \m(e)\le n+1$.
By Lemma~\ref{lem : find and psind from G}, this gives the new proof of the inequality $\find(X,\omega)\le n+1$, as claimed in the paragraph before Corollary~\ref{cor : index}.
\end{remark}

\subsection{Connection and normal splittings}\label{sec : connection and splittings}

For each edge $e\in E(\G)$, the Grothendieck--Birkhoff theorem \cite{MR87176, Klyachko1990} gives the existence of an equivariant splitting
\begin{equation}\label{eqn: splitting}
    N_{C_e/X}\cong \cO_{C_e}(a_1)\oplus\cdots\oplus \cO_{C_e}(a_{n-1})
\end{equation}
where $N_{C_e/X}$ is the normal bundle of $C_e$ in $X$ (cf.\ \cite[equation~(1.1)]{GZ01}).
The following definition is a combinatorial abstraction of this splitting \cite[\S2.1]{Guillemin_Zara_2002}.

\begin{definition}\label{def : connection}
    A \emph{compatible connection} $\nabla$ on $\Gfat$ consists of a bijection
    \[
    \nabla_\epsilon\colon E(\G)^\pm _v\longrightarrow E(\G)^\pm_w
    \]
    for each $\epsilon=vw\in E(\G)^\pm$ satisfying:
    \begin{enumerate}
        \item $\nabla_{\overline{\epsilon}} = \nabla_{\epsilon}^{-1}$.
        \item For each $\epsilon_i\in E(\G)^\pm_v$ there exists $a_i\in \ZZ$ with
    \begin{equation}\label{eqn : connection def a}
    \alpha(\epsilon_i') = \alpha(\epsilon_i) - a_i\alpha(\epsilon)
    \end{equation}
    where $\epsilon_i'=\nabla_\epsilon(\epsilon_i)$.
    \end{enumerate}
\end{definition}

\begin{remark}\label{rem : connection notation}
    Since every GKM graph is $\GKM_2$ (see Remark~\ref{rem : gkm2 vs gkm3}), $\nabla_\epsilon(\epsilon)=\overline{\epsilon}$ holds automatically for all $\epsilon\in E(\G)^\pm$, and the corresponding $a_i$ in \eqref{eqn : connection def a} is $2$.
    Thus, for $\epsilon=vw$, we usually write $E(\G)_v^\pm = \{\epsilon_0,\dots,\epsilon_{n-1}\}$ with $\epsilon_0=\epsilon$ and $a_0=2$.
    The integers $a_1,\dots,a_{n-1}$ correspond to the \emph{normal directions}.
    Moreover, the integers associated to $\epsilon$ and $\overline{\epsilon}$ are the same, so a connection $\nabla$ determines integers $2,a_1,\dots,a_{n-1}$ for every unoriented edge $e\in E(\G)$.
\end{remark}

The splitting~\eqref{eqn: splitting} always induces a compatible connection on $\Gfat$ \cite[Theorem~1.1.2]{GZ01}.
However, when $\Gfat$ is not $\GKM_3$, there can be multiple compatible connections on $\Gfat$.

\begin{remark}\label{rem : vertex span is constant}
    By equation~\eqref{eqn : connection def a} and the existence of a compatible connection for $\Gfat$, the vector space $V_v:=\mathrm{span}_\RR\left\langle\alpha(\epsilon) : \epsilon\in E(\G)^\pm_v\right\rangle$ is the same for each $v\in V(\G)$.
    By the effectiveness assumption in Definition~\ref{def : gkm space}, we have $\dim V_v=r$, the rank of the torus $T$.
\end{remark}

By \eqref{eqn : magnitude from G}, the integers $a_1,\dots,a_{n-1}$ in \eqref{eqn : connection def a} always satisfy
\begin{equation}\label{eqn: sum of a is m}
\m(e)=2+a_1+\cdots+a_{n-1}.
\end{equation}

\begin{definition}\label{def : N(e,a)}
    Let $\nabla$ be a compatible connection on $\Gfat$.
    For any $e\in E(\G),a\in\ZZ$, define
\[
\cN^\nabla(e,a):=\left|\{i\in\{1,\dots,n-1\} : a_i \ge a\}\right|
\]
where $a_1,\dots,a_{n-1}$ are defined by $\nabla$ and $e$ via \eqref{eqn : connection def a} and Remark~\ref{rem : connection notation}.
\end{definition}

\subsection{Equivariant cohomology}

The \emph{$T$-equivariant cohomology} of $X$, denoted $H_T^*(X;\QQ)$, is defined as the ordinary cohomology with $\QQ$ coefficients of the Borel construction \cite{Hsiang_1975, MR4655919}.
It is an algebra over the coefficient ring
\[
H_T^*:=H_T^*(\mathrm{point};\QQ)=\QQ[t_1,\dots,t_r]
\]
where $t_1,\dots,t_r$ are the \emph{equivariant parameters}, which are homogeneous of degree two.
For each vertex $v\in V(\G)$, let $p_v\in X^T$ be the associated fixed point.
The $T$-equivariant inclusion map $\{p_v\}\hookrightarrow X$ induces the \emph{restriction map} of $H_T^*$-algebras
\[
\iota_v\colon H_T^*(X;\QQ) \longrightarrow H_T^*.
\]
Taking direct sum over all fixed points, we get
\[
\iota\colon H_T^*(X;\QQ)\longrightarrow \bigoplus_{v\in V(\G)} H_T^* \cong \bigoplus_{v\in V(\G)}\QQ[t_1,\dots,t_r].
\]
By the GKM theorem (named after \cite[Theorem~1.2.2]{GKM98}, see also \cite[Theorem~1.7.3]{GZ01} and \cite{MR766741} for equivariant formality), $\iota$ is injective and its image is described as follows.

Identifying $t_1,\dots,t_r$
with a linear basis of $\Lie(T)^\vee$,
the GKM theorem states that
\[
H_T^*(X;\QQ)\cong \left\{
    (f_v)_{v\in V(\G)}\in \bigoplus_{v\in V(\G)}\QQ[t_1,\dots,t_r]\,\,\,
    :\,\,\,
    \alpha(\epsilon)\mid f_{\src(\epsilon)}-f_{\dst(\epsilon)}\,\,\,\forall \,\,\, \epsilon\in E(\G)^\pm
\right\}.
\]

\textit{Important classes:}
For each $v\in V(\G)$ we have a \emph{Thom class}
\[
\Th_v\in H_T^{2n}(X;\QQ)
\]
defined by the property that for each $w\in V(\G)$ its restriction is
\[
(\Th_v)_w = \begin{cases}
    0 & w\neq v
    \\
    \prod_{\epsilon\in E(\G)^\pm_v}\alpha(\epsilon)
    &
    w = v.
\end{cases}
\]
Similarly, each edge $e=\{v,w\}\in E(\G)$ has a Thom class $\Th_e\in H_T^{2n-2}(X;\QQ)$ defined by
    \[
        (\Th_e)_s = \begin{cases}
            0 & s\notin\{v,w\}
            \\
            \prod_{\epsilon\in E(\G)^\pm_v\setminus\{vw\}} \alpha(\epsilon)
            & s = v
            \\
            \prod_{\epsilon\in E(\G)^\pm_w\setminus\{wv\}} \alpha(\epsilon)
            & s = w
        \end{cases}
    \]
    for all $s\in V(\G)$.
    One easily uses the GKM theorem to verify that $\Th_v,\Th_e$ are indeed cohomology classes.
    Moreover, for each $e=\{v,w\}\in E(\G)$ we have
    \begin{equation}\label{eqn : thom edge and vertices}
        \Th_v - \Th_w = \alpha(vw) \Th_e.
    \end{equation}

\section{Equivariant Gromov--Witten theory}
\label{sec : GW}

\subsection{Equivariant Gromov--Witten invariants}

Let $(X,\omega,T)$ be a Hamiltonian GKM space, $\beta\in H_2(X;\ZZ)$, and $m,g\in\ZZ_{\ge 0}$.
The \emph{equivariant Gromov--Witten invariants} of $(X,\omega,T)$ in curve class $\beta$ and genus $g$ with $m$ marked points are the $H_T^*$-multilinear map
\[
GW_{g,m}^{X;\beta}\colon H_T^*(X;\QQ)^{\oplus m}
\longrightarrow H_T^*\cong \QQ[t_1,\dots,t_r].
\]
Its construction can be summarised as
\[
GW_{g,m}^{X;\beta}(x_1,\dots,x_m) :=
\int_{\left[\cM_{g,m}(X;\beta)\right]^\mathrm{vir}} \prod_{i=1}^m \ev_i^* (x_i)
\]
where $\cM_{g,m}(X;\beta)$ is the moduli space of stable maps from (connected) nodal curves of genus $g$ with $m$ marked points to $X$ in class $\beta$, $\left[\cM_{g,m}(X;\beta)\right]^\mathrm{vir}$ is the virtual fundamental class, and $\ev_i\colon\cM_{g,m}(X;\beta)\rightarrow X$ is the evaluation map at the $i$-th marked point.
We refer to \cite{Cox_Katz_1999} for details of the construction in the algebraic case, and to \cite{Hirschi_2023} and references therein for the general Hamiltonian symplectic setting in which we work.

By construction, $GW_{g,m}^{X;\beta}$ is homogeneous of degree
\[
- \vdim_\RR(\cM_{g,m}(X;\beta)) = 
-2\left( (1-g)(n-3) + m + \int_\beta c_1(TX) \right).
\]
This gives rise to the following simple key observation.

\begin{lemma}\label{lem:  GW key idea}
    Let $(X,\omega,T)$ be a Hamiltonian GKM space and $x_i\in H_T^{2d_i}(X;\QQ)$ for each $i\in \{1,\dots,m\}$. If
    \[
    GW_{g,m}^{X;\beta}(x_1,\dots,x_m)\neq 0
    \]
    then
    \[
    \left(\sum_{i=1}^m d_i\right) - (1-g) (n-3) - m
    \ge 
    \int_\beta c_1(TX).
    \]
\end{lemma}
\begin{proof}
    This follows immediately from
    \[
    \deg_\RR GW_{g,m}^{X;\beta}(x_1,\dots,x_m)
    =
    \sum_{i=1}^m 2d_i - \vdim_\RR(\cM_{g,m}(X;\beta))
    \ge 0
    \]
    since $H_T^*\cong\QQ[t_1,\dots,t_r]$ is supported in non-negative degree.
\end{proof}

We will obtain the magnitude bounds announced in the introduction by finding suitable non-vanishing equivariant Gromov--Witten invariants and applying Lemma~\ref{lem:  GW key idea}.

\begin{remark}
    Compactness of $X$ (or at least of $\cM_{g,m}(X;\beta)$) is an essential ingredient in our arguments.
    If $\cM_{g,m}(X;\beta)$ is not compact, one may instead \emph{define} $GW_{g,m}^{X;\beta}$ by the equivariant localisation formula (see \S\ref{sec: localisation formula} and cf.\ \cite[\S3.5.2]{LS17}).
    Then, however, the values of $GW_{g,m}^{X;\beta}$ will generally lie in $\QQ(t_1,\dots,t_r)$ rather than $\QQ[t_1,\dots,t_r]$, and the argument in Lemma~\ref{lem:  GW key idea} no longer holds as $\QQ(t_1,\dots,t_r)$ has non-zero elements in negative degree.
\end{remark}

\subsection{Localisation formula}
\label{sec: localisation formula}

Our tool to identify non-vanishing Gromov--Witten invariants is the localisation formula (Theorem~\ref{thm : localisation formula}).
The idea of applying equivariant localisation \cite{Atiyah_Bott_1984} to the moduli space of stable maps goes back to \cite{Kontsevich_torus_enumeration, Manin_1995, Graber_Pandharipande_1999}.
As a result, the computation of Gromov--Witten invariants is reduced to combinatorics of the GKM graph $\Gfat=(\G,\alpha)$ of $(X,\omega,T)$ and, in positive genus, Hodge integrals.

For our purposes, it is sufficient to restrict to genus zero.
Fix $m\in\ZZ_{\ge 0}$ and $\beta\in H_2(X;\ZZ)$.
Let $\Gamma_{0,m}(X;\beta)$ be the set of decorated trees $\Gdec=(\Gamma, \dec{f},\dec{d},\dec{s})$ consisting of:
\begin{itemize}
    \item a tree $\Gamma$ (i.e., a simple undirected connected acyclic graph),
    \item a map of graphs $\dec{f}\colon \Gamma\rightarrow \G$ sending vertices to vertices and edges to edges, preserving their incidence relations.
    \item a \emph{degree map} $\dec{d}\colon E(\Gamma)\rightarrow \ZZ_{\ge 1}$ such that
    \begin{equation}\label{eqn : degree condition}
    \sum_{e\in E(\Gamma)} \dec{d}(e) \left[C_{\dec{f}(e)}\right] = \beta,
    \end{equation}
    \item a \emph{marking map} $\dec{s}\colon \{1,\dots,m\}\rightarrow V(\Gamma)$.
\end{itemize}

Note that $\Gamma_{0,m}(X;\beta)$ is finite since \eqref{eqn : degree condition} implies
$|E(\Gamma)|\le\omega(\beta) / \min_{e\in E(\G)}\omega(C_e)$.

\begin{example}
Recall the GKM space $(\CC\PP^3,\omega_{FS},T)$ from Example~\ref{ex: GKM graph of CPn}.
The underlying graph of its GKM graph is the complete graph on $\{p_0,\dots,p_3\}$.
With $H_2(\CC\PP^3;\ZZ)=\ZZ\langle H \rangle$ where $H$ is a curve class, the decorated tree $\Gdec$ on the left below belongs to $\Gamma_{0,3}(X;5H)$.
Each vertex $v$ is labelled by $\dec{f}(v)$, each edge $e$ by $\dec{d}(e)$, and the red flags illustrate $\dec{s}$.

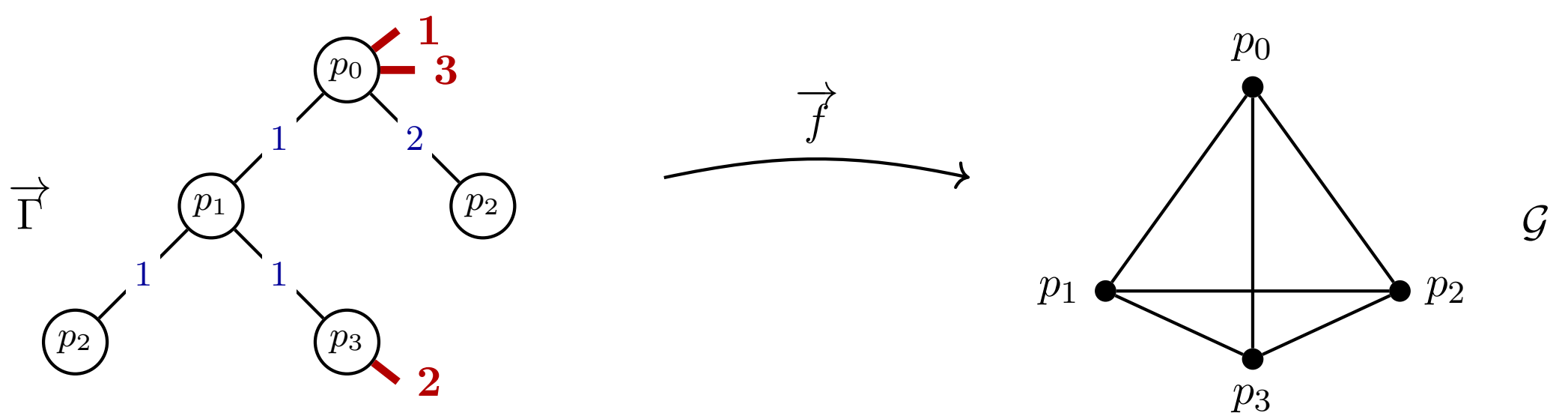
\begin{center}
\begin{tikzpicture}[
    vertex/.style={circle, draw, fill=black, inner sep=0pt, minimum size=5pt},
    labeledvertex/.style={circle, draw, thick, fill=white, inner sep=2pt, font=\footnotesize},
    every node/.style={font=\normalsize}
]

\begin{scope}[shift={(-5,0.15)}]
    \node at (-2.8,0) {$\Gdec$};

    \node[labeledvertex] (g1) at (0,1.2) {$p_0$};
    \node[labeledvertex] (g2) at (-1.2,0) {$p_1$};
    \node[labeledvertex] (g3) at (1.2,0) {$p_2$};
    \node[labeledvertex] (g4) at (0,-1.2) {$p_3$};
    \node[labeledvertex] (g5) at (-2.4,-1.2) {$p_2$};

    \draw[thick] (g1) -- node[midway, fill=white, inner sep=2pt, font=\footnotesize, text=blue!60!black] {$1$} (g2);
    \draw[thick] (g1) -- node[midway, fill=white, inner sep=2pt, font=\footnotesize, text=blue!60!black] {$2$} (g3);
    \draw[thick] (g2) -- node[midway, fill=white, inner sep=2pt, font=\footnotesize, text=blue!60!black] {$1$} (g4);
    \draw[thick] (g2) -- node[midway, fill=white, inner sep=2pt, font=\footnotesize, text=blue!60!black] {$1$} (g5);

    \draw[line width=2pt, red!70!black] (g1) -- ++(0.45,0.35) node[right, font=\bfseries\normalsize, red!70!black] {$\mathbf{1}$};
    \draw[line width=2pt, red!70!black] (g1) -- ++(0.6,0) node[right, font=\bfseries\normalsize, red!70!black] {$\mathbf{3}$};
    \draw[line width=2pt, red!70!black] (g4) -- ++(0.45,-0.35) node[right, font=\bfseries\normalsize, red!70!black] {$\mathbf{2}$};
\end{scope}

\draw[->, thick] (-2.2,0.4) to[bend left=12] node[above] {$\dec{f}$} (0.5,0.4);

\begin{scope}[shift={(3,0)}]
    \node at (2.5,0) {$\G$};

    \node[vertex, label={[font=\normalsize]above:$p_0$}] (v0) at (0,1.2) {};
    \node[vertex, label={[font=\normalsize]left:$p_1$}] (v1) at (-1.3,-0.6) {};
    \node[vertex, label={[font=\normalsize]right:$p_2$}] (v2) at (1.3,-0.6) {};
    \node[vertex, label={[font=\normalsize]below:$p_3$}] (v3) at (0,-1.2) {};

    \draw[thick] (v0) -- (v1);
    \draw[thick] (v0) -- (v2);
    \draw[thick] (v0) -- (v3);
    \draw[thick] (v1) -- (v2);
    \draw[thick] (v1) -- (v3);
    \draw[thick] (v2) -- (v3);
\end{scope}

\end{tikzpicture}
\end{center}
\end{example}

\begin{definition}
    Given $\Gdec\in \Gamma_{0,m}(X;\beta)$ and $x_1,\dots,x_m\in H_T^*(X;\QQ)$, the \emph{contribution $GW_{\Gdec}(x_1,\dots,x_m)\in\QQ(t_1,\dots,t_r)$ of $\Gdec$}
    is defined by
    \begin{align}\label{eqn: localisation contrib}
        G&W_{\Gdec}(x_1,\dots,x_m) := 
        \frac{1}{|\mathrm{Aut}(\Gdec)|}
        \prod_{e\in E(\Gamma)}\frac{h(\dec{f}(e), \dec{d}(e))}{\dec{d}(e)}
        \\\nonumber
        &\prod_{v\in V(\Gamma)} \alpha(\dec{f}(v))^{\val(v)-1}
        \prod_{i\in \dec{s}^{-1}(v)} (x_i)_{\dec{f}(v)}
        \prod_{\epsilon\in E(\Gamma)^\pm_v}\alpha_{(\epsilon,v)}^{-1}
        \left(\sum_{\epsilon\in E(\Gamma)^\pm_v} \alpha_{(\epsilon,v)}^{-1}\right)^{|\dec{s}^{-1}(v)| + \val(v) - 3}.
    \end{align}
    Here, 
    \begin{itemize}
        \item $\alpha_{(\epsilon,v)} := \alpha(\dec{f}(\epsilon)) / \dec{d}(\pi(\epsilon))$, where $\pi\colon E(\Gamma)^\pm \rightarrow E(\Gamma)$ forgets the chosen vertex,
        \item $\alpha(\dec{f}(v)):= \prod_{\epsilon\in E(\G)^\pm_{\dec{f}(v)}}\alpha(\epsilon)$,
        \item $\mathrm{Aut}(\Gdec)$ is the group of automorphisms of $\Gamma$ preserving the decorations $\dec{f},\dec{d},\dec{s}$.
        \item For any $e=\{v,w\}\in E(\G), d\in\ZZ_{>0}$, pick an orientation $\epsilon=vw\in E(\G)^\pm$ and let $E(\G)^\pm_v=\{\epsilon,\epsilon_1,\dots,\epsilon_{n-1}\}$.
        Let $\nabla$ be the compatible connection on $\Gfat$ induced by the splitting of $N_{C_e/X}$ in equation \eqref{eqn: splitting}.
        Let $\epsilon_i' := \nabla_\epsilon(\epsilon_i)$ and let $a_1,\dots,a_{n-1}\in\ZZ$ be determined by $\nabla$ as in equation \eqref{eqn : connection def a}.
        Then
        \begin{equation}\label{eqn: h factor}
            h(e,d):= \frac{(-1)^d d^{2d}}{(d!)^2 \alpha(\epsilon)^{2d}}
            \prod_{i=1}^{n-1} b\left(\frac{\alpha(\epsilon)}{d},\alpha(\epsilon_i),da_i\right)
        \end{equation}
        where
        \begin{equation}\label{eqn: b factor}
        b(u,w,a) := \begin{cases}
            \prod_{j=0}^a (w-ju)^{-1} & a\ge 0
            \\
            \prod_{j=1}^{-a-1}(w+ju) & a<0.
        \end{cases}
        \end{equation}
        Using $\alpha(\overline{\epsilon})=-\alpha(\epsilon)$ and equation~\eqref{eqn : connection def a}, it is easy to verify that $h(e,d)$ does not depend on the choice of $\epsilon$ for $e$.
    \end{itemize}
\end{definition}

The following formula was proven in \cite{LS17} in the complex algebraic case and in \cite{Hirschi_2023} in the Hamiltonian symplectic case.
Even though it holds for any genus (with a more involved definition of $GW_\Gdec$), the genus zero version is enough for us.

\begin{theorem}\label{thm : localisation formula}
    Let $(X,\omega,T)$ be a Hamiltonian GKM space and $x_1,\dots,x_m\in H_T^*(X;\QQ)$. Then
    \[
    GW_{0,m}^{X;\beta}(x_1,\dots,x_m) = 
    \sum_{\Gdec\in \Gamma_{0,m}(X;\beta)} GW_\Gdec(x_1,\dots,x_m).
    \]
\end{theorem}

\begin{remark}\label{rem: connection independence}
    One can use any compatible connection $\nabla$ on $\Gfat$ to define $h(e,d)$ as in equation~\eqref{eqn: h factor}.
    By \cite[\S3.5, Theorem~3.6]{Holmes_Muratore_2025b}, any such $\nabla$ leads to the same value of $h(e,d)$ and therefore $GW_\Gdec(x_1,\dots,x_m)$.
    This observation is necessary to compute Gromov--Witten invariants in terms of $\Gfat$, since the splitting of $N_{C_e/X}$ in equation~\eqref{eqn: splitting} is not determined by $\Gfat$ \cite[Example~3.12]{KSZ_moment_maps_in_a_box}.
\end{remark}

\subsection{Analysis of the localisation formula}
We record the following two properties of $GW_{\Gdec}$ for later use.
As before, let $(X,\omega,T)$ be a Hamiltonian GKM space with GKM graph $\Gfat=(\G,\alpha)$.

\begin{lemma}\label{lem: PD at marked point}
    Let $\Gdec\in\Gamma_{0,m}(X;\beta)$ and $m>0$.
    Then for any $v\in V(\G)$, we have
    \[
    GW_{\Gdec}(x_1,\dots,x_{i-1}, \Th_v, x_{i+1}, \dots, x_m) = 0
    \]
    unless $\dec{f}(\dec{s}(i)) = v$.
\end{lemma}
\begin{proof}
    This follows immediately from \eqref{eqn: localisation contrib} since the factor
    \[
    (x_i)_{\dec{f}(\dec{s}(i))} = (\Th_v)_{\dec{f}(\dec{s}(i))}
    \]
    is zero unless $\dec{f}(\dec{s}(i))=v$.
\end{proof}

The next lemma concerns decorated trees $\Gdec$ of the following shape.

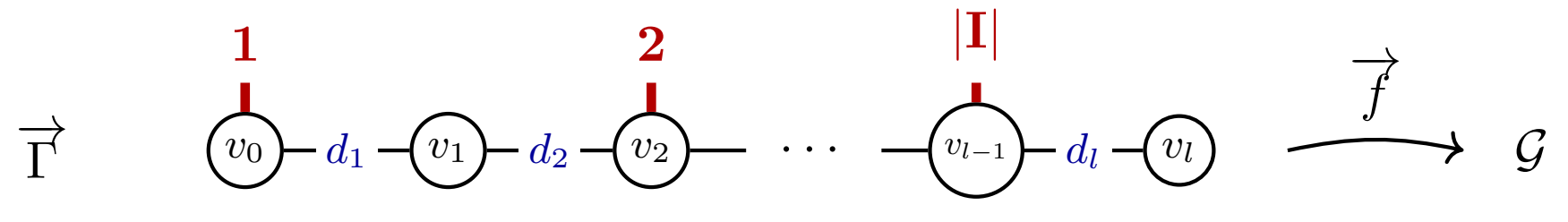
\begin{center}
\begin{tikzpicture}[
    vertex/.style={circle, draw, fill=black, inner sep=0pt, minimum size=5pt},
    labeledvertex/.style={circle, draw, thick, fill=white, inner sep=2pt, font=\footnotesize},
    smalllabeledvertex/.style={circle, draw, thick, fill=white, inner sep=1pt, font=\scriptsize},
    every node/.style={font=\normalsize}
]

\node at (-1.5,0) {$\Gdec$};

\node[labeledvertex] (g0) at (0,0) {$v_0$};
\node[labeledvertex] (g1) at (1.5,0) {$v_1$};
\node[labeledvertex] (g2) at (3,0) {$v_2$};
\node at (4.2,0) {$\cdots$};
\node[labeledvertex] (glm1) at (5.4,0) {\scalebox{0.8}{$v_{l-1}$}};
\node[labeledvertex] (gl) at (6.9,0) {$v_l$};

\draw[thick] (g0) -- node[midway, fill=white, inner sep=2pt, font=\footnotesize, text=blue!60!black] {$d_1$} (g1);
\draw[thick] (g1) -- node[midway, fill=white, inner sep=2pt, font=\footnotesize, text=blue!60!black] {$d_2$} (g2);
\draw[thick] (g2) -- (3.7,0);
\draw[thick] (4.7,0) -- (glm1);
\draw[thick] (glm1) -- node[midway, fill=white, inner sep=2pt, font=\footnotesize, text=blue!60!black] {$d_l$} (gl);

\draw[line width=2pt, red!70!black] (g0) -- ++(0,0.5) node[above, font=\bfseries\normalsize, red!70!black] {$\mathbf{1}$};
\draw[line width=2pt, red!70!black] (g2) -- ++(0,0.5) node[above, font=\bfseries\normalsize, red!70!black] {$\mathbf{2}$};
\draw[line width=2pt, red!70!black] (glm1) -- ++(0,0.5) node[above, font=\bfseries\normalsize, red!70!black] {$\mathbf{|I|}$};

\draw[->, thick] (7.7,0) to[bend left=12] node[above] {$\dec{f}$} (9,0);
\node at (9.5,0) {$\G$};

\end{tikzpicture}
\end{center}

\begin{lemma}\label{lem: path Gamma nonzero}
    Let $P=v_0\dots v_l$ be a path on $\G$ and let $I\subseteq \{0,\dots,l\}$.
    Let $\Gdec$ be a decorated tree with $|I|$ marked points such that
    \begin{enumerate}
        \item $\Gamma$ is a path of length $l$,
        \item  $\dec{f}$ sends the $i$-th vertex of $\Gamma$ to $v_i$, and
        \item if $I=\{i_1<\dots<i_{|I|}\}$
        then $\dec{s}$ sends $j\in \{1,\dots,|I|\}$ to the $i_j$-th vertex of $\Gamma$.
    \end{enumerate}
    Then
    \[
    GW_{\Gdec}(\Th_{v_{i_1}},\dots,\Th_{v_{{i_{|I|}}}})
    \neq 0.
    \]
\end{lemma}

\begin{proof}
    It suffices to verify that each factor in \eqref{eqn: localisation contrib} is nonzero.
    By 2-independence of GKM graphs, we have $\alpha(\epsilon)\neq 0$ for all $\epsilon\in E(\G)^\pm$ and $\alpha(\epsilon_1),\alpha(\epsilon_2)$ are linearly independent whenever $\epsilon_1\neq\epsilon_2$ and $\src(\epsilon_1)=\src(\epsilon_2)$.
    Therefore, all factors of the shapes
    \[
    h(\dec{f}(e),\dec{d}(e))
    \,\,\,\,,\,\,\,\,
    \alpha(\dec{f}(v))^{\val(v)-1}
    \,\,\,\,,\,\,\,\,
    \alpha_{(\epsilon,v)}^{-1}
    \]
    are nonzero.
    Moreover, the factors
    \[
    (x_j)_{\dec{f}(\dec{s}(j))}  = (\Th_{v_{i_j}})_{v_{i_j}} = \alpha(v_{i_j})
    \]
    are also nonzero.
    It remains to show that
    \[
    \left(\sum_{\epsilon\in E(\Gamma)^\pm_v} \alpha_{(\epsilon,v)}^{-1}\right)^{|\dec{s}^{-1}(v)| + \val(v) - 3}
    \]
    is nonzero for each $v\in V(\Gamma)$.
    As $\Gamma$ is a path, the sum in brackets has either one or two terms.
    \begin{enumerate}
        \item If $\val(v)=1$ then $E(\Gamma)^\pm_v = \{\epsilon\}$ and the sum is $\alpha_{(\epsilon,v)}^{-1}$, which is nonzero by definition.

        \item If $\val(v)=2$ then $E(\Gamma)^\pm_v=\{\epsilon_1,\epsilon_2\}$ and the sum is
        \[
        \alpha_{(\epsilon_1, v)}^{-1} + \alpha_{(\epsilon_2, v)}^{-1}
        = \frac{\alpha_{(\epsilon_1, v)} + \alpha_{(\epsilon_2, v)}}{\alpha_{(\epsilon_1, v)}\alpha_{(\epsilon_2, v)}}
        \]
        which is nonzero because of 2-independence.
    \end{enumerate}
    Thus, $GW_{\Gdec}(\Th_{v_{i_1}},\dots,\Th_{v_{i_{|I|}}})$ is nonzero, as desired.
\end{proof}

Finally, we consider the following three decorated trees $\Gdec$ and analyse when their contribution to a specific invariant is polynomial in the equivariant parameters.
From left to right, they correspond to cases 1., 2., 3. in the following proposition.

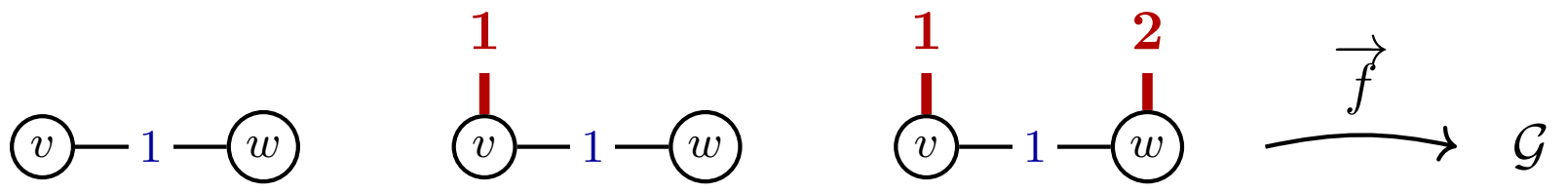
\begin{center}
    \begin{tikzpicture}[
    vertex/.style={circle, draw, fill=black, inner sep=0pt, minimum size=5pt},
    labeledvertex/.style={circle, draw, thick, fill=white, inner sep=2pt, font=\footnotesize},
    smalllabeledvertex/.style={circle, draw, thick, fill=white, inner sep=1pt, font=\scriptsize},
    every node/.style={font=\normalsize}
]

\node[labeledvertex] (a1) at (0.5,0) {$v$};
\node[labeledvertex] (a2) at (2,0) {$w$};
\draw[thick] (a1) -- node[midway, fill=white, inner sep=2pt, font=\footnotesize, text=blue!60!black] {$1$} (a2);

\node[labeledvertex] (b1) at (3.5,0) {$v$};
\node[labeledvertex] (b2) at (5,0) {$w$};
\draw[thick] (b1) -- node[midway, fill=white, inner sep=2pt, font=\footnotesize, text=blue!60!black] {$1$} (b2);
\draw[line width=2pt, red!70!black] (b1) -- ++(0,0.5) node[above, font=\bfseries\normalsize, red!70!black] {$\mathbf{1}$};

\node[labeledvertex] (c1) at (6.5,0) {$v$};
\node[labeledvertex] (c2) at (8,0) {$w$};
\draw[thick] (c1) -- node[midway, fill=white, inner sep=2pt, font=\footnotesize, text=blue!60!black] {$1$} (c2);
\draw[line width=2pt, red!70!black] (c1) -- ++(0,0.5) node[above, font=\bfseries\normalsize, red!70!black] {$\mathbf{1}$};
\draw[line width=2pt, red!70!black] (c2) -- ++(0,0.5) node[above, font=\bfseries\normalsize, red!70!black] {$\mathbf{2}$};

\draw[->, thick] (8.8,0) to[bend left=12] node[above] {$\dec{f}$} (10.1,0);
\node at (10.6,0) {$\G$};

\end{tikzpicture}
\end{center}

\begin{proposition}\label{prop: polynomiality}
    Let $\nabla$ be a compatible connection on $\Gfat$ and let $e=\{v,w\}\in E(\G)$.
    Let $\Gdec\in \Gamma_{0,m}(X;[C_e])$ be a decorated tree consisting of a single edge $\{\tilde{v},\tilde{w}\}$ of degree $\dec{d}(\{\tilde{v},\tilde{w}\})=1$, as well as $\dec{f}(\tilde{v})=v$ and $\dec{f}(\tilde{w})=w$.
    For the markings $\dec{s}\colon \{1,\dots,m\}\rightarrow V(\Gamma)$ and classes $x_1,\dots,x_m\in H_T^*(X;\QQ)$, consider the following three cases.
    \begin{enumerate}
        \item[(1.)] $m=0$ (no markings, no input classes),
        \item[(2.)] $m=1$ and $\dec{s}(1)=\tilde{v}$ with $x_1=\Th_v$,
        \item[(3.)] $m=2$ and $\dec{s}(1)=\tilde{v}, \dec{s}(2)=\tilde{w}$ with $x_1=\Th_v, x_2=\Th_w$.
    \end{enumerate}
    If $GW_\Gdec(x_1,\dots,x_m)\in\QQ[t_1,\dots,t_r]$ then the following conclusions hold in the three cases, respectively, where $k$ is the complexity of $(X,\omega,T)$ (see Definition~\ref{def : complexity}).
    \begin{enumerate}
        \item[(1.)] $\cN^\nabla(e,0)\le k$ and, if $\Gfat$ is $\GKM_3$, $\cN^\nabla(e,0)=0$.
        \item[(2.)] $\cN^\nabla(e, 1)\le k$ and, if $\Gfat$ is $\GKM_3$, $\cN^\nabla(e, 1)=0$.
        \item[(3.)] $\cN^\nabla(e, 2)\le k$ and, if $\Gfat$ is $\GKM_3$, $\cN^\nabla(e,2)=0$.
    \end{enumerate}
\end{proposition}
\begin{proof}
    Let $\epsilon = vw\in E(\G)^\pm$ and $E(\G)^\pm_v=\{\epsilon,\epsilon_1,\dots,\epsilon_{n-1}\}$.
    For each $i\in\{1,\dots,n-1\}$, let $\epsilon_i'=\nabla_\epsilon(\epsilon_i)$ and let $a_1,\dots,a_{n-1}\in\ZZ$ be determined by $\nabla$ via equation~\eqref{eqn : connection def a}.
    Specialising equation~\eqref{eqn: localisation contrib} to $\Gdec$ in the three cases gives
    \[
    GW_\Gdec(x_1,\dots,x_m) = 
    \begin{cases}
        \alpha(\epsilon)\alpha(\overline{\epsilon})h(e,1) & \text{(case 1.)}\\
        \alpha(v)\alpha(\overline{\epsilon})h(e,1)& \text{(case 2.)}\\
        \alpha(v)\alpha(w)h(e,1) & \text{(case 3.).}
    \end{cases}
    \]
    In each case, the factor $-\alpha(\epsilon)^{-2}$ in equation~\eqref{eqn: h factor} cancels and we obtain
    \[
    GW_\Gdec(x_1,\dots,x_m) = \prod_{i=1}^{n-1}b(\alpha(\epsilon), \alpha(\epsilon_i), a_i)\cdot \gamma_i
    \]
    where
    \begin{equation}\label{eqn: poly cases}
    \gamma_i = \begin{cases}
        1 & \text{(case 1.)}\\
        \alpha(\epsilon_i) & \text{(case 2.)}\\
        \alpha(\epsilon_i)\alpha(\epsilon_i') & \text{(case 3.)}.
    \end{cases}
    \end{equation}
    For each $i\in\{1,\dots,n-1\}$, define
    \[
    J_i := \left\{j\in \ZZ : \begin{cases}
        0\le j\le a_i&\text{(case 1.)}\\
        1\le j \le a_i&\text{(case 2.)}\\
        1\le j \le a_i-1&\text{(case 3.)}
    \end{cases}\right\}
    \]
    and let $I:=\{i\in\{1,\dots,n-1\} : J_i\neq \emptyset\}$ and $I^c:=\{1,\dots,n-1\}\setminus I$.

    Let us analyse the denominator of $GW_\Gdec(x_1,\dots,x_m)$.
    From \eqref{eqn : connection def a}, \eqref{eqn: b factor}, and \eqref{eqn: poly cases}, it is slightly tedious but elementary to see that for each $i\in I^c$ there exists $S_i\subset\ZZ$ such that
    \begin{equation}\label{eqn : denom and num factors}
        GW_\Gdec(x_1,\dots,x_m)
        =
        \frac
        {\prod_{i\in I^c}\prod_{j\in S_i}(\alpha(\epsilon_{ i}) +  j\alpha(\epsilon))}
        {\prod_{i\in I}\prod_{j\in J_i} (\alpha(\epsilon_i) - j\alpha(\epsilon))}.
    \end{equation}
    
    Since $GW_\Gdec(x_1,\dots,x_m)$ is nonzero by Lemma~\ref{lem: path Gamma nonzero} and polynomial in the equivariant parameters by assumption, each of the factors in the denominator of \eqref{eqn : denom and num factors} cancels with some factor in the numerator.
    In particular, for each $i\in I$, we have
    \[
    \alpha(\epsilon_i)\in\mathrm{span}_\RR\left\langle\alpha(\epsilon), \alpha(\epsilon_{\hat{i}}) : \hat{i}\in I^c\right\rangle.
    \]
    Since $I\cap I^c=\emptyset$, we obtain
    \[
    \dim \mathrm{span}\left\langle \alpha(\epsilon),\alpha(\epsilon_1),\dots,\alpha(\epsilon_{n-1}) \right\rangle
    \le n-|I|.
    \]
    By Remark~\ref{rem : vertex span is constant}, the left hand side is the rank $r$ of $T$.
    The above inequality is therefore
    \begin{equation}\label{eqn : |I| upper bound}
    |I| \le n-r = k.
    \end{equation}
    Moreover, by definition of $I$ we have
    \[
    |I|=\begin{cases}
        \cN^\nabla(e,0)&\text{(case 1.)}\\
        \cN^\nabla(e,1)&\text{(case 2.)}\\
        \cN^\nabla(e,2)&\text{(case 3.)}
    \end{cases}
    \]
    so equation~\eqref{eqn : |I| upper bound} gives the desired upper bound.
    When $(X,\omega,T)$ is $\GKM_3$, any cancellation between two factors in \eqref{eqn : denom and num factors} would contradict 3-independence, so $|I|=0$ as desired.
\end{proof}

\section{Curve class combinatorics}\label{sec : curve class combinatorics}

In this section, we provide the curve class combinatorics used in \S\ref{sec : proofs of technical results}-\ref{sec : quadrics}.
Throughout, let $(X,\omega,T)$ be a Hamiltonian GKM space with GKM graph $\Gfat=(\G,\alpha)$ and
fix distinct vertices $v,w\in V(\G)$ and a \emph{positive} cohomology class $a\in H^2(X;\RR)$, meaning that $\int_{C_e} a > 0$ for all $e\in E(\G)$.
For example, the class $a=[\omega]$ is always positive while $a=c_1(TX)$ is positive whenever $(X,\omega,T)$ is positive (see Definition~\ref{def : pos and almost positive}).
For a path $P=v_0\dots v_\len$, we will write
\[
a(P) := \sum_{i=1}^\len \int_{C_{v_{i-1}v_i}}a.
\]

\begin{definition}\label{def: P and Pmin}
    Let $\cP(v,w)$ be the set of paths from $v$ to $w$ (see Definition~\ref{def : paths and subpaths}), and let $\cPmin(v,w)\subseteq \cP(v,w)$ be the subset of paths $P$ minimising $a(P)$.
\end{definition}

Note that $\cP(v,w)\neq \emptyset$ as $\G$ is connected (see \S\ref{sec : gkm spaces and their graphs}), and $\cPmin(v,w)\neq\emptyset$ as $a$ is positive.

\begin{remark}
    If $P\in \cPmin(v,w)$ then $[P]$ is minimal among paths from $v$ to $w$ (cf.\  Definition~\ref{def : partial order on H2})
    as $[Q]<[P]$ implies $a(Q)<a(P)$ by positivity of $a$.
    However, the converse is not true in general.
\end{remark}

\begin{lemma}\label{lem: min path has no repeated vert}
    Each path $P\in\cPmin(v,w)$ visits each vertex at most once.
\end{lemma}
\begin{proof}
    Otherwise, if $u\in V(\G)$ is visited multiple times by $P$, we may omit the subpath of $P$ from the first to the last occurrence of $u$ to obtain a path $Q$ from $v$ to $w$ with $a(Q)<a(P)$, contradicting minimality of $a(P)$.
\end{proof}

\begin{lemma}\label{lem: vertex set}
    Let $P,Q\in \cPmin(v,w)$ visit the same set of vertices.
    Then $P=Q$.
\end{lemma}

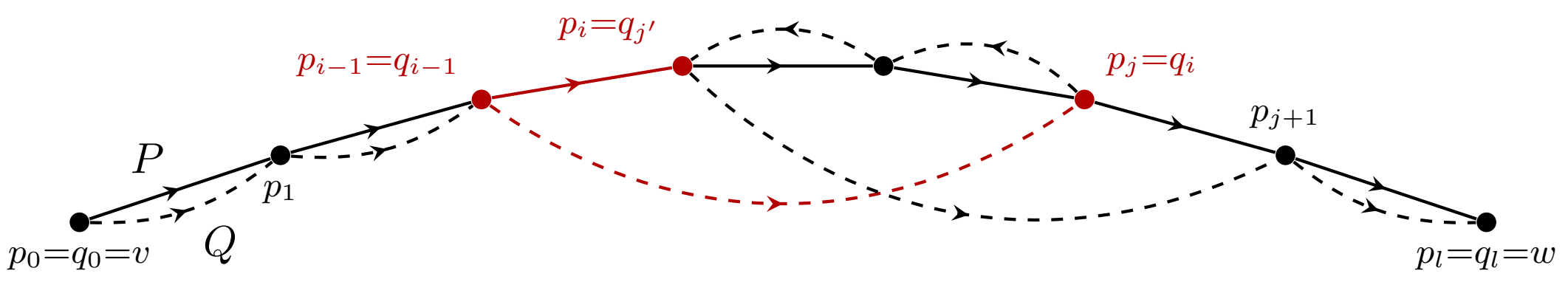
\begin{figure}[htb]
    \centering
\begin{tikzpicture}[
    vertex/.style={circle, fill=black, minimum size=5pt, inner sep=0pt},
    redvertex/.style={circle, fill=red!70!black, minimum size=5pt, inner sep=0pt},
    every label/.style={font=\footnotesize},
    >={stealth},
    midarrow/.style={postaction={decorate, decoration={markings, mark=at position 0.5 with {\arrow{>}}}}},
    P/.style={thick},
    Q/.style={thick, dashed},
    Pred/.style={thick, red!70!black},
    Qred/.style={thick, dashed, red!70!black}
]
    \node[vertex, label=below:{$p_0{=}q_0{=}v$}] (p0) at (0,0) {};
    \node[vertex, label=below:{$p_1$}] (p1) at (1.8,0.6) {};
    \node[redvertex, label={[text=red!70!black]above left:{$p_{i-1}{=}q_{i-1}$}}] (pi1) at (3.6,1.1) {};
    \node[redvertex, label={[text=red!70!black]above left:{$p_i{=}q_{j'}$}}] (pi) at (5.4,1.4) {};
    \node[vertex] (pi2) at (7.2,1.4) {};
    \node[redvertex, label={[text=red!70!black]above right:{$p_j{=}q_i$}}] (pj) at (9,1.1) {};
    \node[vertex, label=above:{$p_{j+1}$}] (pj1) at (10.8,0.6) {};
    \node[vertex, label=below:{$p_l{=}q_l{=}w$}] (pl) at (12.6,0) {};

    \draw[P, midarrow] (p0) -- node[above left] {$P$} (p1);
    \draw[P, midarrow] (p1) -- (pi1);
    \draw[Pred, midarrow] (pi1) -- (pi);
    \draw[P, midarrow] (pi) -- (pi2);
    \draw[P, midarrow] (pi2) -- (pj);
    \draw[P, midarrow] (pj) -- (pj1);
    \draw[P, midarrow] (pj1) -- (pl);

    \draw[Q, midarrow] (p0) to[bend right=20] node[below right] {$Q$} (p1);
    \draw[Q, midarrow] (p1) to[bend right=20] (pi1);
    \draw[Qred, midarrow] (pi1) to[bend right=35] (pj);
    \draw[Q, midarrow] (pj) to[bend right=35] (pi2);
    \draw[Q, midarrow] (pi2) to[bend right=35] (pi);
    \draw[Q, midarrow] (pi) to[bend right=35] (pj1);
    \draw[Q, midarrow] (pj1) to[bend right=20] (pl);
\end{tikzpicture}
    \caption{Illustration of the proof of Lemma~\ref{lem: vertex set}.}
    \label{fig: prop: vertex set}
\end{figure}

\begin{proof}
    Assume for contradiction that $P, Q\in \cPmin(v,w)$ are distinct but visit the same set of vertices.
    By Lemma~\ref{lem: min path has no repeated vert}, each vertex is visited at most once by each path, so $\len(P)=\len(Q)$ and we may write
    \[
    P=p_0p_1\dots p_{l-1}p_l
    \hspace{3mm}\text{ and }\hspace{3mm}
    Q = q_0 q_1\dots q_{l-1}q_l
    \hspace{3mm}\text{ with }\hspace{3mm}p_0=q_0=v,\,\,p_l=q_l=w,
    \]
    where $\{p_0,\dots,p_{l}\}=\{q_0,\dots,q_{l}\}$.
    Let $i:=\min\{i:p_i\neq q_i\}\in\{1,\dots,l-1\}$ and let $j\in\{1,\dots,l\}$ be such that $q_i=p_{j}$.
    Note that $j>i$.
    Similarly, let $j'\in\{1,\dots,l\}$ be such that $p_i=q_{j'}$ and note that $j'>i$.
    Both $j$ and $j'$ are unique by Lemma~\ref{lem: min path has no repeated vert}.

    We have the following inequality.
    \begin{equation}\label{eqn: lem_ineq_1}
    a(p_{i-1}p_i) < a(p_{i-1}p_i\dots p_j)
    \le a(q_{i-1}q_i)
    \end{equation}
    The left inequality follows from positivity of $a$, while the right inequality holds by minimality of $a(P)$. Indeed, if the right inequality failed, then the path 
    $$R:=p_0\dots p_{i-1}p_j\dots p_l$$
    from $v$ to $w$ would have $a(R)< a(P)$, contradicting $P\in \cPmin(v,w)$.

    By symmetry (exchanging the roles of $P$ and $Q$), the same argument yields
    \[
    a(q_{i-1}q_i) < a(p_{i-1}p_i)
    \]
    This contradicts \eqref{eqn: lem_ineq_1}, finishing our proof by contradiction.
\end{proof}

Let us also recall a complete set of relations among $\mathsf{C}:=\{[C_e] : e\in E(\G)\}\subset H_2(X;\ZZ)$.
Recall from Lemma~\ref{lem : find and psind from G} that $H_2(X;\ZZ)$ is a free $\ZZ$-module spanned by $\mathsf{C}$. Since $H_*(X;\ZZ)$ is torsion-free\footnote{as explained in the footnote on p.\pageref{footnote : torsion}.}, the Universal Coefficients Theorem gives $H_2(X;\ZZ)\otimes\RR \cong H_2(X;\RR)$.
Thus, it suffices to express the relations among $\mathsf{C}$ in $H_2(X;\RR)$.

\begin{definition}\label{def : oriented loop}
    An \emph{oriented loop} on $\Gfat$ is a path $P=p_0\dots p_{\len(P)}$ on $\Gfat$ (see Definition~\ref{def : paths and subpaths}) with $p_0=p_{\len(P)}$.
    We write $\epsilon_i:=p_{i-1}p_i\in E(\G)^\pm$ and $e_i := \{p_{i-1}, p_i\}\in E(\G)$ for each $i\in \{1,\dots,\len(P)\}$.
    If $\len(P)=3$ then we say that $e_1,e_2,e_3$ is a \emph{triangle} on $\Gfat$.
\end{definition}

With this notation, we recall the following from \cite[Proposition~2.34]{Holmes_Muratore_2025b}.

\begin{proposition}\label{prop : curve class relations}
    For any $\xi\in \Lie(T)$ and oriented loop $P$ on $\Gfat$ with $\epsilon_1,\dots,\epsilon_{\len(P)}\in E(\G)^\pm$ as above, we have the relation
    \begin{equation}\label{eqn : C1 C2 in span C3...Cn}
    \sum_{i=1}^{\len(P)} \langle \alpha(\epsilon_i), \xi \rangle [C_{e_i}] = 0
    \end{equation}
    in $H_2(X;\RR)$.
    Moreover, this is a complete set of relations in the sense that any other relation is obtained as an $\RR$-linear combination of them.
\end{proposition}

\begin{corollary}\label{cor : curve class N-gon}
    Let $P$ be an oriented loop on $\Gfat$ of length $N$ with $e_i$ and $\epsilon_i$ as in Definition~\ref{def : oriented loop}, and assume that $e_1\neq e_2$.
    Then
    \[
    [C_{e_1}],[C_{e_2}]\in\text{span}_\RR\{[C_{e_3}],\dots,[C_{e_N}]\}.
    \]
    In particular, the span of $[C_{e_1}],\dots,[C_{e_N}]$ in $H_2(X;\RR)$ has rank at most $N-2$.
\end{corollary}
\begin{proof}
    By 2-independence of $\Gfat$ and $e_1\neq e_2$, there exist $\xi_1,\xi_2\in \Lie(T)$ such that
    \[
    \langle\alpha(\epsilon_i),\xi_j\rangle = \begin{cases}
        1 & i=j\\
        0 & i \neq j
    \end{cases}
    \]
    for $i,j\in\{1,2\}$.
    Applying Proposition~\ref{prop : curve class relations} with $\xi_1,\xi_2$ thus gives \eqref{eqn : C1 C2 in span C3...Cn}.
\end{proof}

\begin{corollary}\label{cor : curve class triangle}
    Let $e_1,e_2,e_3\in E(\G)$ be a triangle on $\Gfat$ (see Definition~\ref{def : oriented loop}).
    Then $[C_{e_1}],[C_{e_2}],[C_{e_3}]$ lie on the same one-dimensional ray in $H_2(X;\RR)$.
    Moreover, if we have $[C_{e_1}]=[C_{e_2}]=[C_{e_3}]$ then writing $e_1=\{u,v\},e_2=\{v,w\}, e_3=\{w,u\}$ gives
    \[
    \alpha(uw) = \alpha(uv) + \alpha(vw).
    \]
\end{corollary}
\begin{proof}
    First, observe that the vertices of a triangle are automatically distinct, so $e_1,e_2,e_3$ are distinct.
    So Corollary~\ref{cor : curve class N-gon} implies that $[C_{e_2}]=\lambda[C_{e_3}]$ for some $\lambda\in \RR$.
    Since $\omega([C_e])>0$ for all $e\in E(\G)$, $\lambda$ is positive.
    The same argument applies to $(e_1,e_2)$ and $(e_1,e_3)$.

    The \emph{moreover} part follows immediately from $\alpha(\epsilon)=-\alpha(\overline{\epsilon})$ for all $\epsilon\in E(\G)^\pm$ and Proposition~\ref{prop : curve class relations} applied to the closed oriented loop $uvwu$.
\end{proof}

\section{Proofs of magnitude bounds}
\label{sec : proofs of technical results}

In this section, we prove the technical results announced in \S\ref{sec : intro technical} (see Tables~\ref{tab: magnitude bounds} and \ref{tab: line chern bounds}).
All minimality statements about curve classes refer to the strict partial order $<$ on $H_2(X;\ZZ)$ introduced in Definition~\ref{def : partial order on H2}.
Moreover, if $\Gdec$ is a decorated tree as in \S\ref{sec: localisation formula} and $P=p_0\dots p_l$ is a path on $\Gamma$, we write $\dec{f}(P)$ for the path $\dec{f}(p_0)\dots\dec{f}(p_l)$ on $\G$.
Finally, we write $[\dec{f}(\Gdec)]$ for $\sum_{e\in E(\Gamma)}\dec{d}(e)\left[C_{\dec{f}(e)}\right]$ and recall that $k$ is the complexity of $(X,\omega,T)$ (see Definition~\ref{def : complexity}).

\subsection{Unique paths of minimal class}

The following theorem applies in particular to any irreducible edge (considered as a path of length one) and is applied repeatedly in \S\ref{sec: mukai}-\ref{sec : quadrics} to prove our main theorems and the Kobayashi--Ochiai-type Corollary~\ref{cor : index}.

\begin{theorem}\label{thm: unique min path restated}
    Let $(X,\omega,T)$ be a Hamiltonian GKM space of dimension $2n$ with GKM graph $\Gfat$ and distinct $v,w\in V(\G)$.
    Let $P$ be a path from $v$ to $w$ such that $[P]$ is minimal among paths from $v$ to $w$, and assume that there is no other path from $v$ to $w$ in class $[P]$.
    Then
    \[
    \m(P)\le n+1
    \,\,\,\,\,\text{ and }\,\,\,\,\,
    \cN^\nabla(e,2)\le k
    \]
    for each edge $e$ of $P$ and every compatible connection $\nabla$ on $\Gfat$.
    If $\Gfat$ is $\GKM_3$, $\cN^\nabla(e,2)\le k$ may be replaced by $\cN^\nabla(e,2)=0$.
\end{theorem}

\begin{proof}
    We will show that the invariant
    \[
        GW_{0,2}^{X;[P]}(\Th_v, \Th_w)
    \]
    is nonzero.
    Then Lemma~\ref{lem:  GW key idea} gives the desired inequality
    \[
        \m(P) \le 2n - (n-3) - 2 = n+1.
    \]
    To show the desired nonvanishing, recall the localisation formula from Theorem~\ref{thm : localisation formula}
    and let $\Gdec\in \Gamma_{0,2}(X;[P])$.
    By Lemma~\ref{lem: PD at marked point}, $GW_\Gdec(\Th_v, \Th_w)$ vanishes unless
    \[
        \dec{f}(\dec{s}(1)) = v\,\,\,\,\text{and}\,\,\,\,
        \dec{f}(\dec{s}(2)) = w.
    \]
    Hence, we may assume that the above two equations hold.
    Let $\tP$ be the unique minimal path from $\dec{s}(1)$ to $\dec{s}(2)$ on the tree $\Gamma$.
    Then $\dec{f}(\tP)$ is a path from $v$ to $w$ on $\G$ satisfying
    \[
    [\dec{f}(\tP)]\le [P].
    \]
    By minimality of $[P]=[\dec{f}(\Gdec)]$, we must have
    \[
        [\dec{f}(\tP)]= [P] = [\dec{f}(\Gdec)].
    \]
    Irreflexivity of $<$ therefore implies
    \[
    \tP = \Gamma\hspace{3mm}\text{and}\hspace{3mm}
    \dec{d}(e)=1\,\,\,\forall\,e\in E(\Gamma)
    \]
    since otherwise we would have $[\dec{f}(\tP)]< [\dec{f}(\Gdec)]$.
    Since $P$ is the only path in class $[P]$ from $v$ to $w$, we have
    \[
    \dec{f}(\Gamma) = P.
    \]
    Hence, $\Gdec$ is fully determined, as illustrated below, where $P=p_0 p_1\dots p_{l-1} p_l$ with $p_0=v$ and $p_l=w$.

\begin{center}
\begin{tikzpicture}[
    vertex/.style={circle, draw, fill=black, inner sep=0pt, minimum size=5pt},
    labeledvertex/.style={circle, draw, thick, fill=white, inner sep=2pt, font=\footnotesize},
    smalllabeledvertex/.style={circle, draw, thick, fill=white, inner sep=1pt, font=\scriptsize},
    every node/.style={font=\normalsize}
]

\node at (-1.5,0) {$\Gdec$};

\node[labeledvertex] (g0) at (0,0) {$p_0$};
\node[labeledvertex] (g1) at (1.5,0) {$p_1$};
\node[labeledvertex] (g2) at (3,0) {$p_2$};
\node at (4.2,0) {$\cdots$};
\node[labeledvertex] (gl) at (5.4,0) {$p_l$};

\draw[thick] (g0) -- node[midway, fill=white, inner sep=2pt, font=\footnotesize, text=blue!60!black] {$1$} (g1);
\draw[thick] (g1) -- node[midway, fill=white, inner sep=2pt, font=\footnotesize, text=blue!60!black] {$1$} (g2);
\draw[thick] (g2) -- (3.7,0);
\draw[thick] (4.7,0) -- (gl);

\draw[line width=2pt, red!70!black] (g0) -- ++(0,0.5) node[above, font=\bfseries\normalsize, red!70!black] {$\mathbf{1}$};
\draw[line width=2pt, red!70!black] (gl) -- ++(0,0.5) node[above, font=\bfseries\normalsize, red!70!black] {$\mathbf{2}$};

\draw[->, thick] (6.2,0) to[bend left=12] node[above] {$\dec{f}$} (7.5,0);
\node at (8,0) {$\G$};

\end{tikzpicture}
\end{center}

    By Lemma~\ref{lem: path Gamma nonzero}, $GW_\Gdec(\Th_v, \Th_w)\neq 0$ so
    \[
        GW_{0,2}^{X;[P]}(\Th_v, \Th_w) =
        GW_\Gdec(\Th_v, \Th_w)
        \neq 0
    \]
    as required, which yields the desired bound on $\m(P)$.

    To obtain the result on $\cN^\nabla(e,2)$, note that each subpath of $P$ is again of minimal class among paths from its source to its target, and is unique in its class.
    In particular, applying the above argument to a single edge $e=p_{i-1}p_i$ of $P$ shows that
    \[
    GW_{0,2}^{X;[C_e]}(\Th_{p_{i-1}}, \Th_{p_{i}}) = GW_{\Gdec}(\Th_{p_{i-1}}, \Th_{p_{i}})\in H_T^*\cong\QQ[t_1,\dots,t_r]
    \]
    where $\Gdec$ is as in case 3. of Proposition~\ref{prop: polynomiality}.
    By polynomiality in the equivariant parameters, the result on $\cN^\nabla(e,2)$ follows from case 3. of Proposition~\ref{prop: polynomiality}.
\end{proof}

\begin{remark}
    The inequality $\m(P)\le n+1$ in Theorem~\ref{thm: unique min path restated} is reminiscent of Mori's cone theorem \cite{MR662120}, which states that for a Fano manifold $X$, each extremal ray of its cone of curves is generated by the class of a rational curve $C\subset X$ with $-K_X\cdot C\le n+1$.
\end{remark}

\subsection{Connectedness by low magnitudes}
An immediate consequence of the next result is the connectedness of $\Gfat$ by edges of magnitude at most $n+1$ (Corollary~\ref{cor : connectedness}).

\begin{theorem}\label{thm: connectivity restated}
    Let $(X,\omega,T)$ be a Hamiltonian GKM space of dimension $2n$ with GKM graph $\Gfat$ and complexity $k$.
    Then for any distinct vertices $v,w$ of $\G$, there exists a path $P$ from $v$ to $w$ such that each subpath $Q$ of $P$ satisfies
    \begin{equation}
        \m(Q) \le \len(Q)(n-1) + 2.
    \end{equation}
    In particular, each edge $e$ of $P$ has magnitude $\m(e)\le n+1$. 
    
    Moreover, for any compatible connection $\nabla$ on $\Gfat$, each edge $e$ of $P$ satisfies $\cN^\nabla(e,2)\le k$ and, if $\Gfat$ is $\GKM_3$, $\cN^\nabla(e,2) = 0$.
\end{theorem}

\begin{proof}
    The \emph{in particular} part of the statement follows from the first statement by taking $Q$ to be a single edge of $P$, so it remains to prove the first statement and then the \emph{moreover} part.

    Let $a\in H^2(X;\RR)$ be positive (see \S\ref{sec : curve class combinatorics}) and recall the definitions of $\cP(v,w)$ and $\cPmin(v,w)$ from Definition~\ref{def: P and Pmin}.
    Furthermore, let
    \[
    \cPminmax(v,w)\subset \cPmin(v,w)
    \]
    consist of the paths $P\in \cPmin(v,w)$ of maximal length.
    Such a maximal length exists by positivity of $a$.

    Note that any subpath $Q=q_0\dots q_{\len(Q)}$ of $P$ automatically lies in $\cPminmax(q_0,q_{\len(Q)})$, so it suffices to prove the statement for any $P\in\cPminmax(v,w)$.
    Fix such a $P$, let $l=\len(P)$, and write
    \[
    P = v_0\dots v_l\,\,\,\,\text{with}\,\,\,\,v_0=v\,\,\,\text{and}\,\,\, v_l=w.
    \]
    It suffices to show that
    \[
    GW_{0,l+1}^{X;[P]}(\Th_{v_0},\dots,\Th_{v_l})\neq 0
    \]
    since then Lemma~\ref{lem:  GW key idea} implies
    \[
    \m(P) \le (l+1)n - (n-3) - (l+1) = l(n-1)+ 2.
    \]
    As before, recall Theorem~\ref{thm : localisation formula} and let $\Gdec\in \Gamma_{0,l+1}(X;[P])$ with
    \[
        GW_{\Gdec}(\Th_{v_0},\dots,\Th_{v_l})\neq 0.
    \]
    By Lemma~\ref{lem: PD at marked point}, this implies
    \begin{equation}\label{eqn: proof2eqn}
    \dec{f}(\dec{s}(i))= v_{i-1}
    \end{equation}
    for all $i\in\{1,\dots,l+1\}$.
    Let $\tP$ be the unique minimal path on $\Gamma$ from $\dec{s}(1)$ to $\dec{s}(l+1)$.
    Since $\dec{f}(\tP)\in \cP(v,w)$, minimality of $a(P)$ together with positivity of $a$ and
    \[
        a(\dec{f}(\tP)) \le a([\dec{f}(\Gdec)]) = a(P)
    \]
    imply that $\tP = \Gamma$ and $\dec{d}(e)=1$ for all $e\in E(\Gamma)$.
    In particular, $[P]=[\dec{f}(\tP)]$.
    By Lemma~\ref{lem: min path has no repeated vert}, the $v_i$ are pairwise distinct, so \eqref{eqn: proof2eqn} implies $\len(\tP)= |V(\Gamma)|-1\ge l=\len(P)$.
    Since $\dec{f}(\tP)\in \cPmin(v,w)$ and $P\in \cPminmax(v,w)$, we have
    \[
    \len(\tP) = \len(P).
    \]
    Thus, by equation \eqref{eqn: proof2eqn},
    \[
    (\dec{s}(1),\dots,\dec{s}(l+1))
    \]
    is a permutation of $(1,\dots,l+1)$ and $\dec{f}(\tP)$ passes precisely through the vertices $v_0,\dots,v_l$.
    By Lemma~\ref{lem: vertex set}, $\dec{f}(\tP)=P$ and hence $\Gdec$ is the following decorated tree.
    
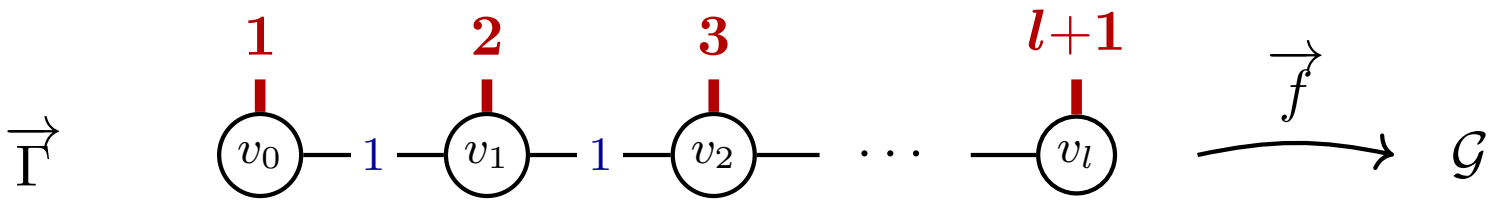
\begin{center}
\begin{tikzpicture}[
    vertex/.style={circle, draw, fill=black, inner sep=0pt, minimum size=5pt},
    labeledvertex/.style={circle, draw, thick, fill=white, inner sep=2pt, font=\footnotesize},
    smalllabeledvertex/.style={circle, draw, thick, fill=white, inner sep=1pt, font=\scriptsize},
    every node/.style={font=\normalsize}
]

\node at (-1.5,0) {$\Gdec$};

\node[labeledvertex] (g0) at (0,0) {$v_0$};
\node[labeledvertex] (g1) at (1.5,0) {$v_1$};
\node[labeledvertex] (g2) at (3,0) {$v_2$};
\node at (4.2,0) {$\cdots$};
\node[labeledvertex] (gl) at (5.4,0) {$v_l$};

\draw[thick] (g0) -- node[midway, fill=white, inner sep=2pt, font=\footnotesize, text=blue!60!black] {$1$} (g1);
\draw[thick] (g1) -- node[midway, fill=white, inner sep=2pt, font=\footnotesize, text=blue!60!black] {$1$} (g2);
\draw[thick] (g2) -- (3.7,0);
\draw[thick] (4.7,0) -- (gl);

\draw[line width=2pt, red!70!black] (g0) -- ++(0,0.5) node[above, font=\bfseries\normalsize, red!70!black] {$\mathbf{1}$};
\draw[line width=2pt, red!70!black] (g1) -- ++(0,0.5) node[above, font=\bfseries\normalsize, red!70!black] {$\mathbf{2}$};
\draw[line width=2pt, red!70!black] (g2) -- ++(0,0.5) node[above, font=\bfseries\normalsize, red!70!black] {$\mathbf{3}$};
\draw[line width=2pt, red!70!black] (gl) -- ++(0,0.5) node[above, font=\bfseries\normalsize, red!70!black] {$\boldsymbol{l}\mathbf{{+}1}$};

\draw[->, thick] (6.2,0) to[bend left=12] node[above] {$\dec{f}$} (7.5,0);
\node at (8,0) {$\G$};

\end{tikzpicture}
\end{center}

    By Lemma~\ref{lem: path Gamma nonzero}, we indeed have
    \[
    GW_{\Gdec}(\Th_{v_0},\dots,\Th_{v_l})\neq 0
    \]
    for this particular $\Gdec$.
    Hence, we have shown that
    \[
        GW_{0,l+1}^{X;[P]}(\Th_{v_0},\dots,\Th_{v_l})
        = GW_{\Gdec}(\Th_{v_0},\dots,\Th_{v_l})\neq 0
    \]
    which proves the desired bound on $\m(P)$.

    The result on $\cN^\nabla(e,2)$ follows as in the proof of Theorem~\ref{thm: unique min path restated} by invoking case 3.~of Proposition~\ref{prop: polynomiality} applied to any edge $e$ of $P$.
    By the above argument, the relevant Gromov--Witten invariant has a unique $\Gdec$ with nonzero contribution, which is a polynomial in the equivariant parameters.
\end{proof}

\subsection{Edges at a fixed vertex}
We include the next result as a precursor of the \emph{local Mukai-type} Theorem~\ref{thm: local Mukai}, and for completeness because of its statement on $\cN^\nabla(e,1)$.

\begin{theorem}\label{thm: magnitude at vertex restated}
    Let $(X,\omega,T)$ be a Hamiltonian GKM space of dimension $2n$ with GKM graph $\Gfat$.
    Let $v\in V(\G)$ and let $e\in E(\G)_v$ have class $[C_e]$ minimal among curve classes of edges at $v$.
    Then
    \[
        \m(e) \le 2
        \,\,\,\,\,\text{ and }\,\,\,\,\,
        \cN^\nabla(e,1)\le k
    \]
    for any compatible connection $\nabla$ on $\Gfat$, or there exists $e_2\in E(\G)_v\setminus\{e\}$ with $[C_{e_2}]=[C_e]$.

    Moreover, if $\Gfat$ is $\GKM_3$, then $\cN^\nabla(e,1)\le k$ may be replaced by $\cN^\nabla(e,1)=0$ in the above statement.
\end{theorem}
\begin{proof}
    Let us assume that there exists no other edge $e_2$ at $v$ with $[C_{e_2}]=[C_e]$.
    For the bound on $\m(e)$, it suffices to show that
    \[
    GW_{0,1}^{X;[C_e]}(\Th_v)\neq 0
    \]
    since Lemma~\ref{lem:  GW key idea} then gives
    \[
    \m(e)\le n - (n-3) - 1 = 2.
    \]
    To show the desired non-vanishing, recall Theorem~\ref{thm : localisation formula} and let $\Gdec\in \Gamma_{0,1}(X;[C_e])$ such that
    \[
    GW_\Gdec(\Th_v)\neq 0.
    \]
    Lemma~\ref{lem: PD at marked point} implies $\dec{f}(\dec{s}(1))=v$.
    Moreover, for each $e'\in E(\Gamma)_{\dec{s}(1)}$ the relation
    \[
    \left[C_{\dec{f}(e')}\right]\le\sum_{e\in E(\Gamma)}\dec{d}(e)\left[C_{\dec{f}(e)}\right]
    = [C_e]
    \]
    implies $\left[C_{\dec{f}(e')}\right]=[C_e]$ by minimality of $[C_e]$ among curve classes of edges at $v$, and hence $e=\dec{f}(e')$.
    Substituting this into the right equality above shows that $\Gamma$ consists of only one edge, and this edge is sent to $e$ and has multiplicity 1.
    Thus, $\Gdec$ is as depicted below, where $w$ is the other vertex of $e$.

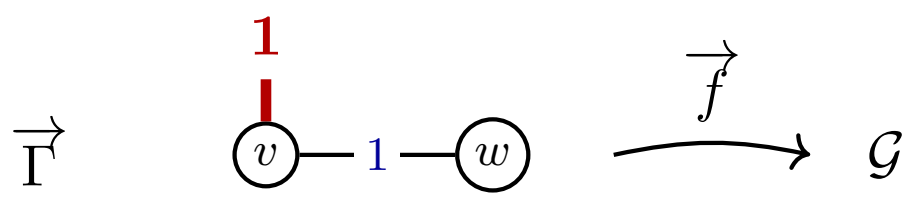
\begin{center}
\begin{tikzpicture}[
    vertex/.style={circle, draw, fill=black, inner sep=0pt, minimum size=5pt},
    labeledvertex/.style={circle, draw, thick, fill=white, inner sep=2pt, font=\footnotesize},
    smalllabeledvertex/.style={circle, draw, thick, fill=white, inner sep=1pt, font=\scriptsize},
    every node/.style={font=\normalsize}
]

\node at (-1.5,0) {$\Gdec$};

\node[labeledvertex] (g0) at (0,0) {$v$};
\node[labeledvertex] (g1) at (1.5,0) {$w$};

\draw[thick] (g0) -- node[midway, fill=white, inner sep=2pt, font=\footnotesize, text=blue!60!black] {$1$} (g1);

\draw[line width=2pt, red!70!black] (g0) -- ++(0,0.5) node[above, font=\bfseries\normalsize, red!70!black] {$\mathbf{1}$};

\draw[->, thick] (2.3,0) to[bend left=12] node[above] {$\dec{f}$} (3.6,0);
\node at (4.1,0) {$\G$};

\end{tikzpicture}
\end{center}

    For this choice of $\Gdec$, Lemma~\ref{lem: path Gamma nonzero} gives
    \[
    GW_\Gdec(\Th_v)\neq 0,
    \]
    so
    \[
    GW_{0,1}^{X;[C_e]}(\Th_v) = GW_\Gdec(\Th_v)\neq 0,
    \]
    proving the desired bound on $\m(e)$.
    The result on $\cN^\nabla(e,1)$ follows from case 2.~of Proposition~\ref{prop: polynomiality} applied to $e$, since $GW_{0,1}^{X;[C_e]}(\Th_v)\in \QQ[t_1,\dots,t_r]$.
\end{proof}

\subsection{A local Mukai-type theorem}
\label{sec : technical : local Mukai}

The next result is stronger than Theorem~\ref{thm: magnitude at vertex restated} in the sense that it produces more edges of the same curve class at a vertex $v$, but it does not make any statement about $\cN^\nabla(e,a)$.
The name \emph{local Mukai-type theorem} is motivated by its applications in \S\ref{sec: mukai}.

\begin{theorem}
\label{thm: local Mukai}
    Let $(X,\omega,T)$ be a Hamiltonian GKM space with GKM graph $\Gfat$. Let $v\in V(\G)$ and let $e=\{v,w\}\in E(\G)_v$ have class $[C_e]$ minimal among curve classes of edges at $v$.
    Then there are at least $\m(e)-2$ vertices $u\in V(\G)$ such that
    \begin{enumerate}
        \item $\{u,v\}\in E(\G)$,
        \item $\{u,w\}\in E(\G)$,
        \item $[C_{\{u,v\}}]=[C_e]$, and
        \item $[C_{\{u,w\}}]\in\RR_{>0}\cdot[C_e]$.
    \end{enumerate}
    In particular, there are at least $\m(e)-1$ edges at $v$ in class $[C_e]$.
\end{theorem}

\begin{figure}[htb]
    \centering
        \begin{tikzpicture}
\node[draw, circle, fill=black, inner sep=1.5pt] (v1) at (0, 0) {};
\node[draw, circle, fill=black, inner sep=1.5pt] (v2) at (1, 1) {};
\node[draw, circle, fill=black, inner sep=1.5pt] (v3) at (0, 2) {};
\node[draw, circle, fill=black, inner sep=1.5pt] (v4) at (-1, 1) {};
\draw[red!70!black, thick] (v2) -- node[pos=0.5, fill=white, inner sep=1pt] {\small 3} (v1);
\draw (v3) -- node[pos=0.35, fill=white, inner sep=1pt] {\small 6} (v1);
\draw[red!70!black, thick] (v3) -- node[pos=0.5, fill=white, inner sep=1pt] {\small 3} (v2);
\draw[red!70!black, thick] (v4) -- node[pos=0.5, fill=white, inner sep=1pt] {\small 3} (v1);
\draw (v4) -- node[pos=0.65, fill=white, inner sep=1pt] {\small 6} (v2);
\draw[red!70!black, thick] (v4) -- node[pos=0.5, fill=white, inner sep=1pt] {\small 3} (v3);
\end{tikzpicture}
\hspace{1cm}
\begin{tikzpicture}
\node[draw, circle, fill=black, inner sep=1.5pt, label=above:$v$] (v) at (0, 0) {};
\node[draw, circle, fill=black, inner sep=1.5pt, label=above:$w$] (w) at (2.5, 0) {};
\node[draw, circle, fill=black, inner sep=1.5pt, label=below:$u_1$] (u1) at (1.2, -1.2) {};
\node[draw, circle, fill=black, inner sep=1.5pt, label=below:$u_i$] (ui) at (2.5, -1.2) {};
\node[draw, circle, fill=black, inner sep=1.5pt, label=below:$u_{n-1}$] (un) at (3.8, -1.2) {};
\draw[thick] (v) -- node[above] {$e$} (w);
\draw (w) -- (u1);
\draw (w) -- (ui);
\draw (w) -- (un);
\node at (1.85, -1.2) {$\cdots$};
\node at (3.15, -1.2) {$\cdots$};
\draw[red!70!black, thick] (v) -- node[pos=0.35, below left] {$P$} (ui);
\end{tikzpicture}

    \caption{Left: the GKM graph of the complex quadric $X=\cQ_3$ (see \S\ref{sec : quadrics}) with $b_2(X)=1$, showing that $\RR_{>0}$ cannot be replaced by $1$ in Theorem~\ref{thm: local Mukai}.
    Each edge $e$ is labelled by $\m(e)$ and irreducible edges are drawn in red.
    Right: illustration of the proof of Theorem~\ref{thm: local Mukai}.}
    \label{fig: local Mukai}
\end{figure}
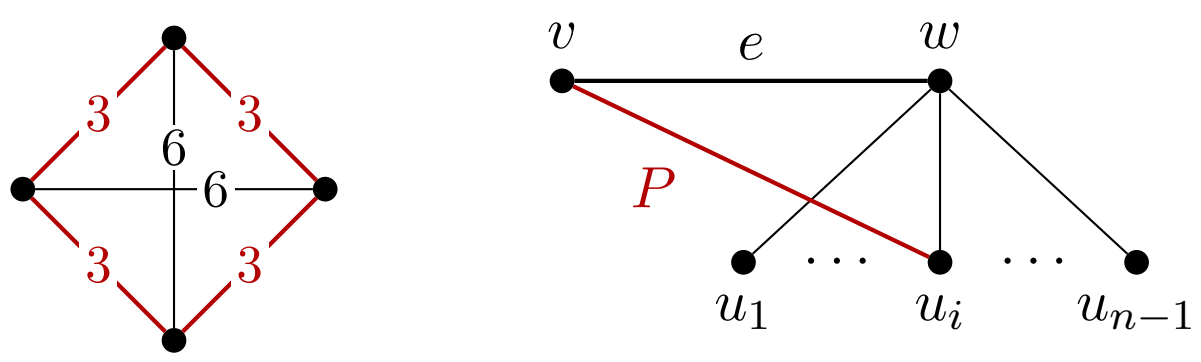

\begin{proof}
    First, note that $u=v$ and $u=w$ do not satisfy the conditions as $\G$ is a simple graph (cf.\ \S\ref{sec : gkm spaces and their graphs}).
    On the other hand, if $u\notin\{v,w\}$ satisfies (1), (2), and (3), then (4) follows from Corollary~\ref{cor : curve class triangle}.

    Let $u_1,\dots,u_{n-1},v$ be the neighbours of $w$ on $\G$.
    They automatically satisfy (2).
    We prove the theorem by producing an upper bound on the number of $i\in\{1,\dots,n-1\}$ such that $u_i$ does not satisfy properties (1) and (3).
    
    So assume that $u_i$ fails either (1) or (3) or both.
    We claim that this implies
    \begin{equation}\label{eqn : local mukai}
        GW_{0,2}^{X;[C_e]}(\Th_v,\Th_{u_i})
    =0.
    \end{equation}
    Indeed, any $\Gdec\in\Gamma_{0,2}(X;[C_e])$ with $GW_\Gdec(\Th_v,\Th_{u_i})\neq 0$ would imply the existence of a path $P$ from $v$ to $u_i$ with $[P]\le [C_e]$.
    Minimality of $[C_e]$ among $E(\G)_v$ then implies $P=vu_i$ and $[P]=[C_e]$, contradicting that $u_i$ fails (1) and (3).
    Thus, $GW_\Gdec(\Th_v,\Th_{u_i})= 0$ for all $\Gdec\in \Gamma_{0,2}(X;[C_e])$ and \eqref{eqn : local mukai} follows.

    Next, we claim that
    \begin{equation}\label{eqn : local mukai 2}
        \alpha(wu_i)\mid GW_{0,2}^{X;[C_e]}(\Th_v,\Th_w).
    \end{equation}
    To see this, note that
    \[
    GW_{0,2}^{X;[C_e]}(\Th_v,\Th_w) - \underbrace{GW_{0,2}^{X;[C_e]}(\Th_v,\Th_{u_i})}_{= 0\text{ by }\eqref{eqn : local mukai}}
    = 
    \alpha(wu_i)GW_{0,2}^{X;[C_e]}(\Th_v,\Th_{\{w,u_i\}})
    \]
    by $H_T^*$-multilinearity of $GW_{0,2}^{X;[C_e]}$ and \eqref{eqn : thom edge and vertices}, so \eqref{eqn : local mukai 2} follows.

    Finally, the path consisting of the single edge $e$ satisfies the assumption of Theorem~\ref{thm: unique min path restated}, and its proof shows that $GW_{0,2}^{X;[C_e]}(\Th_v,\Th_w)\neq 0$.
    By 2-independence of $\Gfat$, the factors $\alpha(wu_i)$ for $i\in\{1,\dots,n-1\}$ are pairwise coprime in $H_T^*\cong\QQ[t_1,\dots,t_r]$.
    Since
    \[
    \deg_\CC GW_{0,2}^{X;[C_e]}(\Th_v,\Th_w)=n+1 - \m(e),
    \]
    and $\deg_\CC \alpha(wu_i)=1$ (where $\deg_\CC$ denotes the degree in the equivariant parameters $t_1,\dots,t_r$), there are at most $n+1-\m(e)$ indices $i\in\{1,\dots,n-1\}$ for which $u_i$ fails (1) or (3).
    Consequently, there are at least $\m(e)-2$ indices $i\in\{1,\dots,n-1\}$ which satisfy (1) and (3).
    They automatically satisfy (2) and hence (4), as required.

    The \emph{in particular} part follows by counting in the original edge $e=\{v,w\}$.
\end{proof}

\subsection{Unique edges of minimal class}

Finally, we include the following result because it implies a dimension bound for the \emph{almost positive} Hamiltonian GKM spaces introduced in \cite[Definition~5.8]{Holmes_Muratore_2025b}.
A complete classification of their GKM graphs will appear in \cite{APC_upcoming}.

\begin{theorem}\label{thm: unique minimal edge restated}
    Let $(X,\omega,T)$ be a Hamiltonian GKM space of dimension $2n$ with GKM graph $\Gfat$.
    Let $e\in E(\G)$ be such that $[C_e]$ is minimal among $\{[C_e]:e\in E(\G)\}$ and no other edge has class $[C_e]$.
    Then
    \[
        \m(e) \le 3-n
        \,\,\,\,\,\text{ and }\,\,\,\,\,
        \cN^\nabla(e,0)\le k
    \]
    for any compatible connection $\nabla$ on $\Gfat$.
    If $\Gfat$ is $\GKM_3$, then $\cN^\nabla(e,0)\le k$ may be replaced by $\cN^\nabla(e,0)=0$ in the above statement.
\end{theorem}

\begin{proof}
    By assumption on $e=vw$, the only possible $\Gdec\in \Gamma_{0,0}(X;[C_e])$
    is the following.
    
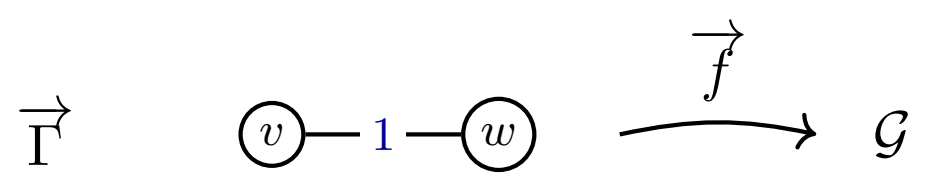
\begin{center}
\begin{tikzpicture}[
    vertex/.style={circle, draw, fill=black, inner sep=0pt, minimum size=5pt},
    labeledvertex/.style={circle, draw, thick, fill=white, inner sep=2pt, font=\footnotesize},
    smalllabeledvertex/.style={circle, draw, thick, fill=white, inner sep=1pt, font=\scriptsize},
    every node/.style={font=\normalsize}
]

\node at (-1.5,0) {$\Gdec$};

\node[labeledvertex] (g0) at (0,0) {$v$};
\node[labeledvertex] (g1) at (1.5,0) {$w$};

\draw[thick] (g0) -- node[midway, fill=white, inner sep=2pt, font=\footnotesize, text=blue!60!black] {$1$} (g1);

\draw[->, thick] (2.3,0) to[bend left=12] node[above] {$\dec{f}$} (3.6,0);
\node at (4.1,0) {$\G$};

\end{tikzpicture}
\end{center}
By Lemma~\ref{lem: path Gamma nonzero}, we have $GW_{0,0}^{X;[C_e]}=GW_\Gdec\neq 0$.
Thus, the desired bound on $\m(e)$ follows from Lemma~\ref{lem:  GW key idea}.

The result on $\cN^\nabla(e,0)$ follows directly from case 1.~of Proposition~\ref{prop: polynomiality} applied to $e$, since $GW_{0,0}^{X;[C_e]}\in \QQ[t_1,\dots,t_r]$.
\end{proof}

\section{Pseudoindex $n+1$ and projective space}
\label{sec : pseudoindex n+1 and projective space}

In this section, we prove the first half of Corollary~\ref{cor : index}, which we restate for convenience.

\begin{theorem}\label{thm: index n+1 with proof}
    Let $(X,\omega,T)$ be a Hamiltonian GKM space of dimension $2n>0$.
    Then
    \[
        \psind(X,\omega)\le n+1.
    \]
    If $\psind(X,\omega)=n+1$ then the GKM graph $\Gfat$ of $(X,\omega,T)$ is isomorphic to a projection of the\footnote{Here, we mean the GKM graph of the standard $(\CC\PP^n, \omega_{FS}, T)$ as in Example~\ref{ex: GKM graph of CPn}, see also Example~\ref{ex: CPn proj criterion}.} GKM graph of $\CC\PP^n$ and $X$ is homotopy equivalent to $\CC\PP^n$.
\end{theorem}

\begin{proof}
    Since $n>0$, $\G$ has at least two vertices.
    By Corollary~\ref{cor : connectedness} they remain connected after discarding all edges $e$ with $\m(e)>n+1$.
    Thus, there exists $e\in E(\G)$ with $\m(e)\le n+1$ and hence $\psind(X,\omega)\le n+1$ by Lemma~\ref{lem : find and psind from G}.

    Assume now that $\psind(X,\omega)=n+1$ and let $e=\{v,w\}\in E(\G)$ have $\m(e)\le n+1$ as above.
    By Lemma~\ref{lem : find and psind from G}, $\m(e)=n+1$.
    Since no edge has lower magnitude, $[C_e]$ is minimal among edge curve classes.
    By Theorem~\ref{thm: local Mukai}, all edges at $v$ and $w$ have the same curve class, and every $u\in V(\G)\setminus\{v,w\}$ is a neighbour of $v$ if and only if it is a neighbour of $w$.
    Iterating this application of Theorem~\ref{thm: local Mukai}, we obtain that $\G$ is a complete graph and all edges have the same curve class.
    By Corollary~\ref{cor : curve class triangle} and Example~\ref{ex: CPn proj criterion}, $\Gfat$ is isomorphic to a projection of $\Gfat_{\CC\PP^n}$.

    Since the GKM graph determines the Betti numbers \cite[Theorem~1.3.2]{GZ01}, we have $b_{2i}(X)=b_{2i}(\CC\PP^n)=1$ for all $0\le i\le n$, which is unimodal.
    The desired homotopy equivalence now follows from \cite[Theorem~1.7]{MR4029085}, noting that every Hamiltonian GKM space has the required \emph{toric 1-skeleton} \cite[Lemma~4.11(i)]{MR3695881}.
\end{proof}

\section{The symplectic Mukai inequality}
\label{sec: mukai}

In this section, we apply the tools developed in \S\ref{sec : curve class combinatorics}--\ref{sec : proofs of technical results} to prove the \emph{generalised symplectic Mukai inequality} \cite{MR4671382}
\[
b_2(X)(\psind(X,\omega) - 1)\le n
\]
for Hamiltonian GKM spaces $(X,\omega,T)$ when either there exists a \emph{special vertex} or the pseudoindex is large (i.e., $\psind(X,\omega) \ge (n+2)/2$).
Here, $2n$ is the dimension of $X$, $b_2(X)$ is its second Betti number, and $\psind(X,\omega)$ is its pseudoindex (see \S\ref{sec : indices}).

\subsection{The large pseudoindex case}

The analogous case of the generalised Mukai conjecture in algebraic geometry was considered by Wi{\'s}niewski in \cite{MR1063222}.
Let $\Gfat=(\G,\alpha)$ be the GKM graph of $(X,\omega,T)$ and let us write $\psind$ for $\psind(X,\omega)$.
By convention, we will call curve classes \emph{linearly independent} if they are linearly independent as elements of $H_2(X;\RR)\cong H_2(X;\ZZ)\otimes\RR$ (cf.\ \S\ref{sec : curve class combinatorics}).

\begin{lemma}\label{lem : large psind primitivity}
    If $\psind\ge (n+2)/2$, the following are equivalent for any $e\in E(\G)$.
    \begin{enumerate}
        \item $\m(e)\le n+1$.
        \item $e$ is irreducible (see Definition~\ref{def : irreducible}).
    \end{enumerate}
\end{lemma}
\begin{proof}
    If $e$ is irreducible then $\m(e)\le n+1$ by Theorem~\ref{thm: unique min path restated}.
    Conversely, if $e$ is not irreducible then $[C_e]= \sum_{i=1}^{s}[C_{e_i}]$ for some $e_1,\dots,e_s\in E(\G)$ with $s\ge 2$ and hence, by Lemma~\ref{lem : find and psind from G}, $\m(e)\ge s\psind\ge 2\psind > n+1$.
\end{proof}

\begin{lemma}\label{lem : large psind primitive connectedness}
    If $\psind\ge (n+2)/2$, any two vertices of $\G$ are connected by a path consisting of irreducible edges.
\end{lemma}
\begin{proof}
    By Corollary~\ref{cor : connectedness}, any two vertices of $\G$ are connected by a path of edges with magnitude at most $n+1$.
    The claim thus follows from Lemma~\ref{lem : large psind primitivity}.
\end{proof}

\begin{lemma}\label{lem : large psind key vertex}
    If $\psind\ge (n+2)/2$ and $b_2(X)\ge 2$, then there exists $v\in V(\G)$ and irreducible $e_1,e_2\in E(\G)_v$ such that $[C_{e_1}], [C_{e_2}]$ are linearly independent.
\end{lemma}
\begin{proof}
    Since $b_2(X)\ge 2$, by Lemma~\ref{lem : find and psind from G} there exist irreducible $\widetilde{e_1},\widetilde{e_2}\in E(\G)$ with $[C_{\widetilde{e_1}}]$, $[C_{\widetilde{e_2}}]$ linearly independent.
    By Lemma~\ref{lem : large psind primitive connectedness} applied to a vertex $v_1$ of $\widetilde{e_1}$ and a vertex $v_2$ of $\widetilde{e_2}$, there exists 
    a path $P=w_0\dots w_\len$ from $v_1$ to $v_2$ consisting of irreducible edges.
    If $[C_{\{w_{i-1}w_i\}}],[C_{\widetilde{e_1}}]$ are linearly dependent for all $i\in\{1,\dots,\len\}$ then $v=v_2$ is the desired vertex.
    Otherwise, let $i\in\{1,\dots,\len\}$ be the minimal index such that $[C_{\{w_{i-1}w_i\}}],[C_{\widetilde{e_1}}]$ are linearly independent.
    Then $v=w_{i-1}$ is the desired vertex.
\end{proof}

\begin{proposition}\label{prop : large psind local product}
    Let $\psind= (n+2)/2$, $v\in V(\G)$, and let $e_1=\{v,w\},e_2=\{v,u\}$ be irreducible edges with $[C_{e_1}], [C_{e_2}]$ linearly independent and $\m(e_1),\m(e_2)=\psind$.
    Then there exists a unique $z\in V(\G)$ such that:
    \begin{enumerate}
        \item $e_1':=\{u,z\}\in E(\G)$, $\alpha(uz)=\alpha(vw)$, $[C_{e_1}]=[C_{e_1'}]$, and
        \item $e_2':=\{w,z\}\in E(\G)$, $\alpha(wz)=\alpha(vu)$, $[C_{e_2}]=[C_{e_2'}]$.
    \end{enumerate}
\end{proposition}

\begin{figure}[htb]
    \centering
    \begin{tikzpicture}
\node[draw, circle, fill=black, inner sep=1.5pt, label=below left:$v$] (v) at (0, 0) {};
\node[draw, circle, fill=black, inner sep=1.5pt, label=below right:$w$] (w) at (1.5, 0) {};
\node[draw, circle, fill=black, inner sep=1.5pt, label=above left:$u$] (u) at (0.5, 1.1) {};
\node[draw, circle, fill=black, inner sep=1.5pt, label=above right:$z$] (z) at (2.0, 1.1) {};
\node at (0, -1) {\,};
\draw[red!70!black, thick] (v) -- node[pos=0.5, fill=white, inner sep=1pt] {\small $e_1$} (w);
\draw[blue!60!black, thick] (v) -- node[pos=0.5, fill=white, inner sep=1pt] {\small $e_2$} (u);
\draw[blue!60!black, thick] (w) -- node[pos=0.5, fill=white, inner sep=1pt] {\small $e_2'$} (z);
\draw[red!70!black, thick] (u) -- node[pos=0.5, fill=white, inner sep=1pt] {\small $e_1'$} (z);
\end{tikzpicture}
\hspace{1cm}
\begin{tikzpicture}
\node[draw, circle, fill=black, inner sep=1.5pt, label=left:$v$] (v) at (0, 0) {};
\node[draw, circle, fill=black, inner sep=1.5pt, label=right:$w$] (w) at (2.5, 0) {};
\node[draw, circle, fill=black, inner sep=1.5pt, label=below:$y_1$] (y1) at (0.6, -1) {};
\node[draw, circle, fill=black, inner sep=1.5pt, label=below:$y_2$] (y2) at (1.15, -1) {};
\node[draw, circle, fill=black, inner sep=1.5pt, label=below:$y_{\frac{n}{2}-1}$] (yn2) at (2, -1) {};

\node[draw, circle, fill=black, inner sep=1.5pt, label=above left:$u_{1}$] (u1) at (-0.6, 1) {};
\node[draw, circle, fill=black, inner sep=1.5pt, label=above:$u_{2}$] (u2) at (-0.2, 1) {};
\node[draw, circle, fill=black, inner sep=1.5pt, label=above:$u_{n/2}$] (un2) at (0.6, 1) {};

\node[draw, circle, fill=black, inner sep=1.5pt, label=above left:$z_{1}$] (z1) at (1.9, 1) {};
\node[draw, circle, fill=black, inner sep=1.5pt, label=above:$z_{2}$] (z2) at (2.3, 1) {};
\node[draw, circle, fill=black, inner sep=1.5pt, label=above:$z_{n/2}$] (zn2) at (3.1, 1) {};

\draw[red!70!black, thick] (v) -- node[pos=0.5, fill=white, inner sep=1pt] {\small $e_1$} (w);
\draw[red!70!black] (w) -- (y1) -- (v);
\draw[red!70!black] (w) -- (y2) -- (v);
\draw[red!70!black] (w) -- (yn2) -- (v);

\draw[blue!70!black] (v) -- (u1);
\draw[blue!60!black] (v) -- (u2);
\draw[blue!60!black] (v) -- (un2);

\draw[blue!70!black] (w) -- (z1);
\draw[blue!60!black] (w) -- (z2);
\draw[blue!60!black] (w) -- (zn2);

\node at (1.6, -1) {$\cdots$};
\node at (0.25, 1) {$\cdots$};
\node at (2.75, 1) {$\cdots$};
\end{tikzpicture}
\hspace{1cm}
\begin{tikzpicture}
\node[draw, circle, fill=black, inner sep=1.5pt, label=below left:$v$] (v) at (0, 0) {};
\node[draw, circle, fill=black, inner sep=1.5pt, label=below right:$w$] (w) at (1.5, 0) {};
\node[draw, circle, fill=black, inner sep=1.5pt, label=below right:$u$] (u) at (0.5, 1.1) {};
\node[draw, circle, fill=black, inner sep=1.5pt, label=above left:$\widetilde{u}$] (up) at (-0.2, 2) {};
\node[draw, circle, fill=black, inner sep=1.5pt, label=above right:$z$] (z) at (2.0, 1.1) {};
\node at (0, -0.5) {\,};
\draw[red!70!black, thick] (v) -- node[pos=0.5, fill=white, inner sep=1pt] {\small $e_1$} (w);
\draw[blue!60!black, thick] (v) -- node[pos=0.5, fill=white, inner sep=1pt] {\small $e_2$} (u);
\draw[blue!60!black, thick] (v) -- node[pos=0.65, fill=white, inner sep=1pt] {\small $\widetilde{e_2}$} (up);
\draw[blue!60!black, thick] (w) -- node[pos=0.5, fill=white, inner sep=1pt] {\small $e_2'$} (z);
\draw[red!70!black, thick] (u) -- node[pos=0.5, fill=white, inner sep=1pt] {\small $e_1'$} (z);
\draw[red!70!black, thick] (up) -- node[pos=0.35, fill=white, inner sep=1pt] {\small $\widetilde{e_1'}$} (z);
\end{tikzpicture}

    \caption{Illustration of the proof of Proposition~\ref{prop : large psind local product}.}
    \label{fig: large psind local product}
\end{figure}
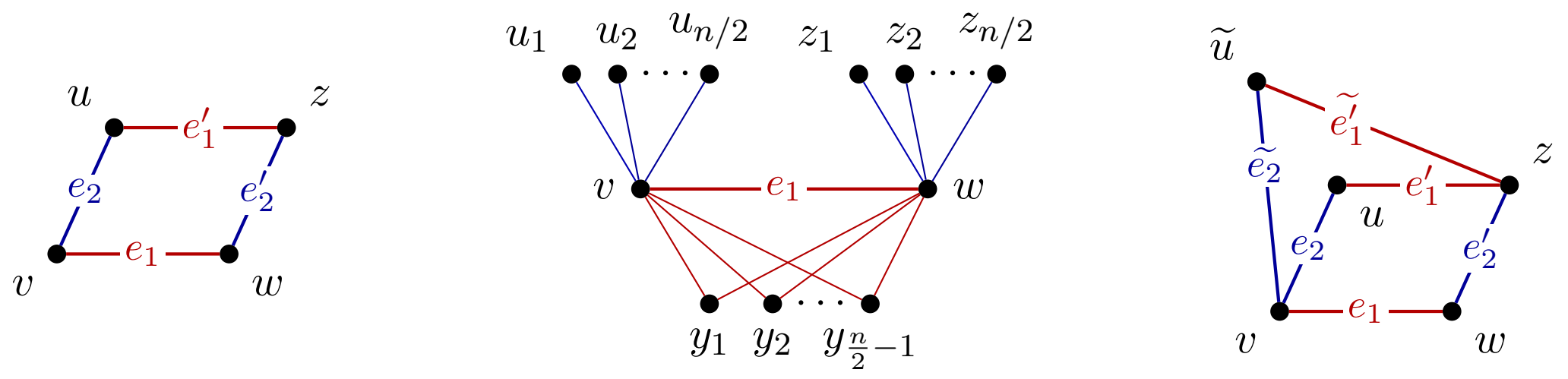

\begin{proof}
    The uniqueness of $z$ is clear by 2-independence of GKM graphs, so let us show existence.
    Consider the path $P=uvw$.
    By Theorem~\ref{thm: unique min path restated} and
    \[
    \m(P)=\m(e_1)+\m(e_2)= 2\psind>n+1,
    \]
    there is another path $Q$ from $u$ to $w$ with $[Q]\le [P]$.
    
    If $Q=uw$ then $e:=\{u,w\}\in E(\G)$ and Corollary~\ref{cor : curve class triangle} applied to $e,e_1,e_2$ shows that $[C_{e_1}]=\lambda [C_{e_2}]$ for some $\lambda\in \RR_{>0}$, contradicting the linear independence assumption.
    Thus, $\len(Q)\ge 2$.
    On the other hand,
    \[
    \len(Q)\psind \le \m(Q)\le \m(P)=2\psind
    \]
    gives $\len(Q)=2$ and $\m(Q)=\m(P)$.
    With $[Q]\le [P]$, this implies $[P]=[Q]$.
    That is, there exists $z\in V(\G)\setminus\{v\}$ such that $Q=uzw$.
    Letting $e_1'$ and $e_2'$ be as in the statement, we get
    \begin{equation}\label{eqn : m e1 e2 e1' e2'}
        \m(e_1)=\m(e_2)=\m(e_1')=\m(e_2')=\psind
    \end{equation}
    and
    \begin{equation}\label{eqn : C e1 e2 e1' e2' sum}
        [Q] = [C_{e_1'}] + [C_{e_2'}]
        = [C_{e_1}] + [C_{e_2}] = [P].
    \end{equation}
    Note that $e_1',e_2'$ are irreducible by Lemma~\ref{lem : large psind primitivity} and that $[C_{e_1}],[C_{e_2'}]$ are linearly independent (otherwise, Corollary~\ref{cor : curve class N-gon} applied to $vwzuv$ contradicts linear independence of $[C_{e_1}],[C_{e_2}]$).
    Similarly, $[C_{e_1'}],[C_{e_2}]$ are linearly independent.
    Theorem~\ref{thm: local Mukai} applied to each of $e_1,e_2,e_2'$ now completely describes the edges of $\G$ at $v$ and $w$ (see middle of Figure~\ref{fig: large psind local product}).
    Namely, there are:
    \begin{itemize}
        \item precisely $n/2$ edges of class $[C_{e_1}]$ at $v$ (call them $e_1, vy_1,\dots,vy_{\frac{n}{2}-1}$),
        \item precisely $n/2$ edges of class $[C_{e_1}]$ at $w$ (call them $e_1, wy_1',\dots,wy_{\frac{n}{2}-1}'$),
        \item precisely $n/2$ edges of class $[C_{e_2}]$ at $v$ (call them $vu_1,\dots,vu_{n/2}$), and
        \item precisely $n/2$ edges of class $[C_{e_2'}]$ at $w$ (call them $wz_1,\dots,wz_{n/2}$).
    \end{itemize}
    In fact, without loss of generality, we have $y_i=y_i'$ for all $i$ (by Theorem~\ref{thm: local Mukai}, $\{v,y_i'\}\in E(\G)$ and $[C_{\{v,y_i'\}}]\in\QQ\cdot[C_{e_1}]$, so $y_i'\in \{y_1,\dots,y_{\frac{n}{2}-1}\}$ by the above classification of edges at $v$).
    
    Since the above construction can be performed for any $u_i$ instead of $u$ (with the same $v$ and $w$), the assignment $u\mapsto z$ gives rise to a map
    \[
    \varphi\colon \{u_1,\dots,u_{n/2}\}\longrightarrow\{z_1,\dots,z_{n/2}\}.
    \]

    \textit{Claim 1: $\varphi$ is injective.}

    \textit{Proof of Claim 1:}
    Distinct $u,\widetilde{u}\in\varphi^{-1}(\{z\})$ would imply via Corollary~\ref{cor : curve class N-gon} applied to $uz\widetilde{u}vu$ that $[C_{e_1'}]\in\mathrm{span}_\RR\langle[C_{e_2}], [C_{\widetilde{e_2}}]\rangle = \RR\cdot[C_{e_2}]$ where $\widetilde{e_2}=\{v,\widetilde{u}\}$ (see right of Figure~\ref{fig: large psind local product}).
    This would contradict linear independence of $[C_{e_1'}],[C_{e_2}]$, proving Claim 1.
    
    Thus, $\varphi$ is a bijection, so let us write without loss of generality $\varphi(u_i)=z_i$.
    In particular, the edge $\{u_i,z_i\}$ exists for each $i$ and has the same properties as $e_1'=\{u,z\}$ above, with $u$ replaced by $u_i$ and $z$ replaced by $z_i$.
    In particular, by \eqref{eqn : C e1 e2 e1' e2' sum},
    \[
    [C_{\{u_i,z_i\}}] = [C_{e_1}] + [C_{\{v,u_i\}}] - [C_{\{w,z_i\}}]
    = [C_{e_1}] + [C_{e_2}] - [C_{e_2'}] = [C_{e_1'}].
    \]

    \textit{Claim 2: We have
    \[
    \sum_{i=1}^{n/2}\left(\alpha(wz_i) - \alpha(vu_i)\right) = 0.
    \]
    Proof of Claim 2:}
    By equation~\eqref{eqn : magnitude from G} for $\m(e_1)$, we have
    \begin{align*}
    \m(e_1)\alpha(vw) = 
    \left(\alpha(vw) - \alpha(wv)\right) +\sum_{i=1}^{n/2}\left(\alpha(vu_i) - \alpha(wz_i)\right)
    + \sum_{i=1}^{\frac{n}{2}-1}\left(\alpha(vy_i)-\alpha(wy_i)\right).
    \end{align*}
    Using $\alpha(vy_i)-\alpha(wy_i)=\alpha(vw)$ (by Corollary~\ref{cor : curve class triangle}) and $\alpha(vw)-\alpha(wv)=2\alpha(vw)$, we obtain
    \[
    \left(\m(e_1) - \left(\frac{n}{2}+1\right)\right)\alpha(vw) = \sum_{i=1}^{n/2}\left(\alpha(vu_i) - \alpha(wz_i)\right)
    \]
    so Claim 2 follows from $\m(e_1)=\psind=(n+2)/2$.

    \textit{Claim 3: There exists $b\in\RR$ such that for all $i\in \{1,\dots,\frac{n}{2}\}$, we have
    \[
    \alpha(wz_i)-\alpha(vu_i)=(b-1)(\alpha(vu_i)+\alpha(vw)).
    \]
    Proof of Claim 3:}
    By Corollary~\ref{cor : curve class N-gon} applied to $vwz_iu_iv$ and using $[C_{\{u_i,z_i\}}]=[C_{e_1'}]$, we may express $[C_{e_1}],[C_{e_2}]$ as linear combination of $[C_{e_1'}],[C_{e_2'}]$.
    Combining this with \eqref{eqn : m e1 e2 e1' e2'} and \eqref{eqn : C e1 e2 e1' e2' sum}, there exists $b\in\RR$ such that
    \begin{equation}\label{eqn : b matrix}
    \begin{pmatrix}
        [C_{e_1}] \\ [C_{e_2}]
    \end{pmatrix}
    =
    \begin{pmatrix}
        b & 1-b \\
        1-b & b
    \end{pmatrix}
    \begin{pmatrix}
        [C_{e_1'}] \\ [C_{e_2'}]
    \end{pmatrix}.
    \end{equation}
    Since $[C_{e_1}],[C_{e_2}]$ are linearly independent, so are $[C_{e_1'}],[C_{e_2'}]$.
    Applying Proposition~\ref{prop : curve class relations} to $vwz_iu_iv$ and taking the $[C_{e_2'}]$-component gives
    \begin{equation}\label{eqn : b and alpha for ui wi}
    (1-b)\alpha(vw)+\alpha(wz_i) + b\alpha(u_iv) = 0.
    \end{equation}
    Rearranging and using $\alpha(\epsilon) = -\alpha(\overline{\epsilon})$ gives the desired equation
    \[
    \alpha(wz_i) - \alpha(vu_i)
    =
    (b-1)\left(\alpha(vw) + \alpha(vu_i)\right).
    \]
    As $b$ only depends on $[C_{e_1}],[C_{e_2}],[C_{e_1'}],[C_{e_2'}]$, it is independent of $i$.
    So Claim 3 is proven.

    \textit{Conclusion:}
    By 2-independence, we have $\alpha(vu_i)+\alpha(vw)\neq 0$, so Claim 3 gives
    \begin{equation}\label{eqn : b not 1 equivalences}
        \exists\,i: \alpha(wz_i)\neq \alpha(vu_i)
        \,\,\,\,\,\Longleftrightarrow\,\,\,\,\,
        b\neq 1
        \,\,\,\,\,\Longleftrightarrow\,\,\,\,\,
        \forall\,i: \alpha(wz_i)\neq \alpha(vu_i).
    \end{equation}
    If $b\neq 1$ then summing Claim 3 over $i$, combining with Claim 2, and dividing by $b-1$ gives
    \begin{equation}\label{eqn : nearly there 1}
    \sum_{i=1}^{n/2}\alpha(vu_i) = -\frac{n}{2}\alpha(vw)
    \end{equation}
    By symmetry, repeating Claims 2 and 3 after interchanging the roles of $v$ and $w$, and using \eqref{eqn : b not 1 equivalences} to see that the corresponding coefficient is not $1$, we similarly obtain
    \begin{equation}\label{eqn : nearly there 2}
    \sum_{i=1}^{n/2}\alpha(wz_i) = -\frac{n}{2}\alpha(wv) = \frac{n}{2}\alpha(vw)
    \end{equation}
    Combining Claim 2, \eqref{eqn : nearly there 1}, and \eqref{eqn : nearly there 2} finally gives the contradiction
    \[
    0 = \sum_{i=1}^{n/2}\left(\alpha(vu_i) - \alpha(wz_i)\right)
    =
    -n\alpha(vw)\neq 0.
    \]
    Thus, we must have $b=1$ and hence \eqref{eqn : b not 1 equivalences} implies $\alpha(wz_i)=\alpha(vu_i)$ for all $i$, and \eqref{eqn : b matrix} implies $[C_{e_1}]=[C_{e_1'}]$ and $[C_{e_2}]=[C_{e_2'}]$.
    Using Proposition~\ref{prop : curve class relations} on $vwzuv$ and extracting the $[C_{e_1}]$-component and $[C_{e_2}]$-component gives the desired $\alpha(vw)=\alpha(uz)$ and $\alpha(vu)=\alpha(wz)$, respectively.
    This finishes the proof of Proposition~\ref{prop : large psind local product}.
\end{proof}

\begin{remark}
    The proof can be substantially shortened if $\Gfat$ is assumed to be $\GKM_3$ (see Definition~\ref{def: k-indep}).
    In that case, we do not need $\varphi$ and can conclude directly from equation \eqref{eqn : b and alpha for ui wi} as follows.
    By 3-independence, there exists at most one $\epsilon\in E(\G)_w$ with $\alpha(\epsilon)\in\mathrm{span}_\RR(\alpha(vw), \alpha(vu_i))$.
    By existence of a compatible connection for $\Gfat$, there exists some $\epsilon'\in E(\G)_w$ with $\alpha(\epsilon')\equiv\alpha(vu_i)$ mod $\alpha(vw)$.
    Combining with \eqref{eqn : b and alpha for ui wi} therefore gives $\epsilon=\epsilon'=wz_i$ and $b=1$.
\end{remark}

\begin{theorem}\label{thm : Mukai large index}
    Let $(X,\omega,T)$ be a Hamiltonian GKM space of dimension $2n>0$ with GKM graph $\Gfat$.
    \begin{enumerate}
        \item If $\psind > (n+2)/2$ then $b_2(X)=1$.
        \item If $\psind = (n+2)/2$ then either
        $b_2(X)=1$; or $b_2(X)=2$ and $\Gfat$ is isomorphic to a projection of the\footnote{Here, we mean the standard GKM action of Example~\ref{ex: GKM graph of CPn} on $\CC\PP^{n/2}$, see also Example~\ref{ex: CPn proj criterion}.} GKM graph of $(\CC\PP^{n/2})^2$.
    \end{enumerate}
\end{theorem}

\begin{proof}
    Note that $b_2(X)\ge 1$ by existence of a symplectic form on $X$.
    So assume $b_2(X)\ge 2$ and $\psind\ge (n+2)/2$.
    
    By Lemma~\ref{lem : large psind key vertex} there exists $v\in V(\G)$ and irreducible $e_1,e_2\in E(\G)_v$ with $[C_{e_1}], [C_{e_2}]$ linearly independent.
    By Theorem~\ref{thm: local Mukai}, there exist at least $\m(e_1)-1$ and $\m(e_2)-1$ edges at $v$ in class $[C_{e_1}]$ and $[C_{e_2}]$, respectively.
    Thus, by Lemma~\ref{lem : find and psind from G} we have
    \[
    2(\psind - 1)\le (\m(e_1) -1) + (\m(e_2) - 1) \le |E(\G)_v|=n
    \]
    and hence
    \[
    \psind = \frac{n+2}{2}
    \]
    and $n$ is even. Thus, $\m(e_1)=\m(e_2)=\psind$ and there are precisely $\psind-1=n/2$ edges in class $[C_{e_1}]$ and $[C_{e_2}]$ at $v$, respectively.

    \textit{Claim 1: Each vertex of $\G$ is adjacent to precisely $n/2$ edges of class $[C_{e_1}]$ and $[C_{e_2}]$ each.}
    
    \textit{Proof of Claim 1:}
    By Proposition~\ref{prop : large psind local product} and Theorem~\ref{thm: local Mukai}, all neighbours of $v$ also have exactly $n/2$ adjacent edges in class $[C_{e_1}]$ and $[C_{e_2}]$, respectively.
    By induction and since $\G$ is connected, the same applies to all vertices of $\G$, proving Claim 1.
    
    By Claim 1 and Lemma~\ref{lem : find and psind from G}, $H_2(X;\ZZ)$ is spanned by $[C_{e_1}],[C_{e_2}]$ and hence $b_2(X)=2$.
    Moreover, each edge is irreducible and has class $[C_{e_1}]$ or $[C_{e_2}]$.
    
    It remains to deduce the desired product structure of $\G$.
    Let $\Gfat_1=(\G_1,\alpha_1)$ where $\G_1$ is the maximal connected subgraph of $\G$ containing $v$ and only edges of class $[C_{e_1}]$, and $\alpha_1=\alpha|_{E(\G_1)^\pm}$.
    Then by Theorem~\ref{thm: local Mukai} and Claim 1, any two vertices on $\G_1$ of distance at most 2 have distance 1.
    That is, $\G_1$ is a complete graph on $n/2+1$ vertices.
    The same holds for $\Gfat_2=(\G_2,\alpha_2)$, where $\G_2$ is the maximal connected subgraph of $\G$ containing $v$ and only edges of class $[C_{e_2}]$, and $\alpha_2=\alpha|_{E(\G_2)^\pm}$.

    \textit{Claim 2: $\Gfat$ is isomorphic to  $\Gfat_1\times \Gfat_2$, in the sense that $V(\G)=V(\G_1)\times V(\G_2)$ and $E(\G)=V(\G_1)\times E(\G_2)\sqcup E(\G_1)\times V(\G_2)$ and $\alpha$ is induced from $\alpha_1,\alpha_2$.}

    \textit{Proof of Claim 2:}
    Consider the map
    \begin{equation}\label{eqn : product map}
    \psi_V\colon V(\G_1)\times V(\G_2)\longrightarrow V(\G),
    \hspace{7mm}
    (w,u)\longmapsto
    \begin{cases}
        v & u=w=v,\\
        u & u\neq v, w=v,\\
        w & u=v, w\neq v,\\
        z & u,w\neq v,
    \end{cases}
    \end{equation}
    where, in the last case, $z$ is obtained from Proposition~\ref{prop : large psind local product} applied to $(u,v,w)$.
    This is well-defined since $\G_1$ and $\G_2$ are both complete graphs.

    Let us prove that $\psi_V$ is injective.
    Define $d\colon V(\G_1)\times V(\G_2)\rightarrow \ZZ_{\ge 0}$ such that $d(w,u)$ is the distance between $v$ and $\psi_V(w,u)$ on $\G$.
    Observe that $d(w,u)$ is precisely 0, 1, 1, and 2 in the four cases of equation~\eqref{eqn : product map}, respectively
    (the only non-trivial step in verifying this is to see that $\{v,z\}\notin E(\G)$, which follows from linear independence of $[C_{e_1}],[C_{e_2}]$ via Corollary~\ref{cor : curve class triangle} applied to $vzuv$).
    Therefore, it suffices to prove injectivity of $\psi_V|_{d^{-1}(i)}$ for each $i\in \{0,1,2\}$.
    \begin{itemize}
        \item For $i\in\{0,1\}$, this is easy to see from $V(\G_1)\cap V(\G_2)=\{v\}$.
        \item For $i=2$, this follows from the right of Figure~\ref{fig: large psind local product}:
        if $\psi_V(\widetilde{w},\widetilde{u})=\psi_V(w,u)= z$ with $u,w,\widetilde{u},\widetilde{w}\neq v$ and, without loss of generality, $u\neq\widetilde{u}$, the closed loop $uz\widetilde{u}vu$ would contradict linear independence of $[C_{e_1}],[C_{e_2}]$ via Corollary~\ref{cor : curve class N-gon}.
    \end{itemize}
    Thus, $\psi_V$ is injective.

    Next, we prove that $\psi_V$ extends to a graph homomorphism $\psi\colon\G_1\times\G_2\rightarrow \G$ compatible with the axial function values.
    By symmetry, it suffices to show for any $\{u_1,u_2\}\in E(\G_2)$ and $w\in V(\G_1)$ that $z_1:=\psi_V(w,u_1)$ and $z_2:=\psi_V(w,u_2)$ satisfy $\{z_1,z_2\}\in E(\G)$ and $\alpha(z_1z_2) = \alpha_2(u_1u_2)$.
    \begin{itemize}
        \item If $u_1=v$ or $u_2=v$, this is immediate from equation~\eqref{eqn : product map} and Proposition~\ref{prop : large psind local product}.
        \item If $u_1,u_2\neq v$ then Theorem~\ref{thm: local Mukai} applied at $w$ gives $\{z_1,z_2\}\in E(\G)$.
        By Corollary~\ref{cor : curve class triangle}, all edges in the triangles $v,u_1,u_2$ and $w,z_1,z_2$ on $\G$ have curve class $[C_{e_2}]$.
        The \emph{moreover} part of the same corollary now gives
        \[
        \alpha(z_1z_2) = \alpha(wz_2)-\alpha(wz_1)
        = \alpha(vu_2) - \alpha(vu_1)
        = \alpha(u_1u_2)
        = \alpha_2(u_1u_2)
        \]
        where the last equality is by definition of $\alpha_2$.
    \end{itemize}
    Thus, $\psi_V$ extends to $\psi\colon\G_1\times\G_2\rightarrow \G$ as desired.
    Finally, $\psi$ is injective since $\psi_V$ is, and $\psi$ is surjective since it is injective and $\G,\G_1\times \G_2$ are both connected, finite, and regular of the same degree $n$.
    Thus, $\psi$ is the desired isomorphism, proving Claim 2.

    \textit{Claim 3: $\Gfat_1\times\Gfat_2$ is isomorphic to a projection of the GKM graph of $\CC\PP^{n/2}\times\CC\PP^{n/2}$ (see Example~\ref{ex: GKM graph of CPn})} 

    \textit{Proof of Claim 3:}
    First, note that $\G_1,\G_2$ are both complete graphs on $n/2 + 1$ vertices.
    Generalising Example~\ref{ex: CPn proj criterion}, one easily checks that $\Gfat_1\times\Gfat_2$ is isomorphic to a projection of  $\CC\PP^{n/2}\times\CC\PP^{n/2}$ if and only if for each $i\in\{1,2\}$ and any distinct $u,w\in V(\G_i)\setminus\{v\}$, we have
    \[
    \alpha(vu)+\alpha(uw)=\alpha(vw).
    \]
    But this is precisely the \emph{moreover} part of Corollary~\ref{cor : curve class triangle}, proving Claim 3.

    Combining Claim 2 and Claim 3 finishes the proof of Theorem~\ref{thm : Mukai large index}.
\end{proof}

\subsection{The special vertex}

In this section, we prove Theorem~\ref{thm : intro Mukai} which states that the generalised symplectic Mukai inequality holds in the presence of a \emph{special vertex}.
Let $(X,\omega,T)$ be a Hamiltonian GKM space with GKM graph $\Gfat$.

\begin{definition}\label{def : Cv and Cvmin}
    For any $v\in V(\G)     $, let $\Cmin_v\subseteq \C_v\subset H_2(X;\RR)$ be the closed convex cones
    \begin{align*}
        \C_v &:= \mathrm{span}_{\RR_{\ge 0}}\left\{[C_e] : e\in E(\G)_v\right\}
        \\
        \Cmin_v &:= \mathrm{span}_{\RR_{\ge 0}}\left\{[C_e] : \begin{matrix}e\in E(\G)_v, \\ [C_e]\text{ is minimal among} \\ \text{curve classes of edges at }v\end{matrix}\right\}
    \end{align*}
    where minimality is as in Definition~\ref{def : partial order on H2}.
    The full \emph{cone of curves} is
    \[
    \C := \mathrm{span}_{\RR_{\ge 0}}\{[C_e] : e\in E(\G)\} \subset H_2(X;\RR).
    \]
\end{definition}

\begin{definition}
    \label{def : special vertex}
    A vertex $v\in V(\G)$ is \emph{special} if $\C_v = \C$.
\end{definition}

\begin{proposition}
\label{prop : vertex mukai inequality}
    Let $(X,\omega,T)$ be a Hamiltonian GKM space of dimension $2n$ with GKM graph $\Gfat$, and let $v\in V(\G)$. Then
    \[
    \dim_\RR(\Cmin_v)(\psind - 1) \le n.
    \]
\end{proposition}
\begin{proof}
    First, note that the statement holds for $\psind \le 0$, so assume $\psind > 0$.
    
    Let $E(\G)_v^\mathrm{min}\subset E(\G)_v$ be the set of edges at $v$ whose curve class is minimal among curve classes of edges at $v$.
    Let $n_v := \#\{ [C_e] : e\in E(\G)_v^\mathrm{min}\}$.
    By Theorem~\ref{thm: local Mukai}, for each $e\in E(\G)_v^\mathrm{min}$ there exist at least $\m(e)-1\ge \psind-1$ edges at $v$ with curve class $[C_e]$.
    Thus,
    \[
    n_v(\psind-1)\le n.
    \]
    Since the dimension of a cone is bounded by the number of its generators, we obtain
    \[
    \dim_\RR (\Cmin_v) (\psind-1) \le n_v(\psind -1) \le n
    \]
    as desired.
\end{proof}

\begin{corollary}\label{cor : key vertex}
    Let $(X,\omega,T)$ be a Hamiltonian GKM space of dimension $2n$ with GKM graph $\Gfat$ and special vertex $v\in V(\G)$.
    Then the generalised symplectic Mukai inequality
    \[
        b_2(X)(\psind-1)\le n
    \]
    holds.
\end{corollary}
\begin{proof}
    By Proposition~\ref{prop : vertex mukai inequality}, it suffices to prove $\Cmin_v=\C_v$ since then $\Cmin_v=\C$ and hence
    \[
    \dim_\RR(\Cmin_v)=\dim_\RR(\C)=b_2(X),
    \]
    where the last equality follows from Lemma~\ref{lem : find and psind from G}.
    
    So let $R\subseteq \C_v=\C$ be an extremal ray.
    By assumption, there is $e\in E(\G)_v$ such that $[C_e]\in R$.
    In particular, there exists $e_{\min}\in E(\G)_v$ such that $[C_{e_{\min}}]\in R$ and no edge at $v$ has curve class $\lambda[C_{e_{\min}}]$ with $\lambda < 1$.
    Since $R$ is an extremal ray of $\C$, $[C_{e_{\min}}]$ is minimal among curve classes of edges at $v$.
    Thus, $[C_{e_{\min}}]\in \Cmin_v$ and hence $R\subseteq \Cmin_v$.

    Since this holds for any extremal ray $R\subseteq\C_v$, we have $\C_v\subseteq \Cmin_v$ and hence $\C_v=\Cmin_v$, which finishes the proof.
\end{proof}

\subsubsection{Examples with special vertex}
\label{sec : examples with special vertex}

We show that homogeneous spaces and smooth Schubert varieties have a special vertex.

\begin{lemma}\label{lem : big aut gives special vertices}
    Let $(X,\omega,T)$ be a Hamiltonian GKM space with GKM graph $\Gfat$ and let $W$ be a group acting on $\Gfat$ by automorphisms (cf.\ Definition~\ref{def : iso and projection of GKM graphs}) such that:
    \begin{enumerate}
        \item $W$ acts transitively on $V(\G)$, and
        \item $W$ preserves curve classes, that is, $[C_{\varphi(e)}]=[C_e]$ for each $\varphi\in W$ and $e\in E(\G)$.
    \end{enumerate}
    Then every vertex $v\in V(\G)$ is special.
\end{lemma}
\begin{proof}
    Assumptions (1) and (2) imply $\C_v=\C_w$ for all $v,w\in V(\G)$ and hence $\C_v=\C$ for all $v\in V(\G)$.
    That is, every vertex of $\G$ is special.
\end{proof}

Typical examples where the assumptions of Lemma~\ref{lem : big aut gives special vertices} hold are homogeneous spaces.
Let $G$ be a connected complex semisimple Lie group and $T_\CC\subset B\subset P\subset G$ a maximal (complex) torus inside a Borel subgroup inside a parabolic subgroup of $G$.
Moreover, let $T \cong (S^1)^r\subset (\CC^\times)^r\cong T_\CC$ be the maximal compact torus inside $T_\CC$.
Then $T_\CC$ acts on $G/P$ by left multiplication and, after replacing $T_\CC$ by its effective quotient, $G/P$ is an algebraic GKM space with respect to $T_\CC$.
In particular, $G/P$ is a Hamiltonian GKM space with respect to the effectively acting quotient torus $T^\mathrm{eff}$ of $T$ and suitable $\omega$ (cf.\ \S\ref{sec : alg Mukai}).

Its GKM graph $\Gfat=(\G,\alpha)$ is constructed as follows \cite{Guillemin_Holm_Zara_2006}.
Let $L$ be the Levi factor of $P$ containing $T_\CC$.
Let $W_G:=N_G(T_\CC)/T_\CC$ and $W_P:=N_L(T_\CC)/T_\CC$ be the Weyl groups of $G$ and $L$, respectively, and let $\Delta_G^+$ and $\Delta_P^+$ be the sets of positive roots of $G$ and $L$, respectively.
For each $\beta\in \Delta_G^+$, let $s_\beta\in W_G$ be the associated reflection.
Then
\begin{itemize}
    \item $V(\G)=W_G/W_P$,
    \item $E(\G)=\{\{[u],[us_\beta]\} : u\in W_G\,,\,\beta\in \Delta_G^+\setminus \Delta_P^+\}$, and
    \item $\alpha(([u],[us_\beta])) = -u(\beta)$.
\end{itemize}

Now, the left multiplication action of $G$ on $G/P$ restricts to an action of $N_G(T_\CC)$, which induces an action of $W_G$ on $\Gfat$ by automorphisms which, on
$V(\G)=W_G/W_P$, is given by left multiplication.
This action is transitive, so the $W_G$-action on $\Gfat$ satisfies assumption (1) of Lemma~\ref{lem : big aut gives special vertices}.
On the other hand, $G$ is connected and $H_2(G/P;\ZZ)$ is discrete, so the former acts trivially on the latter.
Hence, the $W_G$-action on $\Gfat$ also satisfies assumption (2) of the lemma.
Thus, every vertex of $\Gfat$ is special.
By Corollary~\ref{cor : key vertex}, we have obtained a new proof of the following result, which was previously obtained by \cite{MR4671382} (see also \cite{MR2415341}) using different methods.

\begin{corollary}
    Let $G/P$ be a homogeneous space as above. Then
    $G/P$ satisfies the generalised symplectic Mukai inequality.
\end{corollary}

Next, we extend this result to smooth Schubert varieties.
We recall the following background from \cite[\S3--4]{MR2072765}.
For any $[w]\in W_G/W_P$, the associated Schubert variety is
\[
X_{[w]}:= \overline{BwP/P}\subseteq G/P.
\]
Here, we have $X_{[w_1]}\subseteq X_{[w_2]}$ if and only if $[w_1]\le [w_2]$ in the Bruhat order of $W_G/W_P$.
In particular,
\begin{equation}\label{eqn : schubert homology}
H_*(X_{[w]};\ZZ) = \bigoplus_{[w']\le [w]} \ZZ\cdot [X_{[w']}]
\end{equation}
where the homological degree of $[X_{[w']}]$ is $2\len(w'^P)$, i.e., twice the length of the minimum length representative $w'^P$ of $[w']$.
In particular, homology is supported in even degree.
Since $G/P$ is itself a Schubert variety (associated to $[w_0]$ where $w_0\in W_G$ is the longest element), the map
\[
H_*(X_{[w]};\ZZ)\longrightarrow H_*(G/P;\ZZ)
\]
is injective, since a free basis of the former is a subset of a free basis of the latter.
When $X_{[w]}$ is smooth, $(X_{[w]},\omega|_{X_{[w]}},T^\mathrm{eff})$ is a Hamiltonian GKM space\footnote{Here, $T^\mathrm{eff}$ now means the effective quotient of $T$ for the action on $X_{[w]}$.} and its GKM graph $\Gfat_{[w]}$ is the full subgraph of $\Gfat$ induced by the vertex set $\{[u]:u\in W_G\,,\, [u] \le [w]\}$.
In fact, $(X_{[w]},T_{\CC}^\text{eff})$ is an algebraic GKM space where $T_{\CC}^\text{eff}$ is the effective quotient of $T_\CC$ (cf.\ \S\ref{sec : alg Mukai}).

\begin{lemma}\label{lem : smooth schubert has special vertex}
    Let $X_{[w]}\subseteq G/P$ be the Schubert variety associated to $[w]\in W_G/W_P$ and assume $X_{[w]}$ is smooth.
    Then the vertex $[1]\in V(\G_{[w]})$ is special.
\end{lemma}
\begin{proof}
    Let $e=\{[u],[us_\alpha]\}\in E(\G_{[w]})$.
    By equation~\eqref{eqn : schubert homology}, we have integers $c_{[w']}$ such that
     \[
     [C_e] = \sum_{\substack{[w']\le [w] \\ \len(w'^P) = 1}} c_{[w']}[X_{[w']}]
     \]
     In each summand, the condition $\len(w'^P)=1$ implies that $[w'] = [s_\gamma]$ for some simple root $\gamma\in\Delta_G^+\setminus\Delta_P^+$.
     It follows that $X_{[w']} = C_{\{[1],[s_\gamma]\}}$ and $\{[1],[s_\gamma]\}\in E(\G_{[w]})$ because $[s_\gamma]=[w']\le [w]$.
     Moreover, $c_{[w']}$ is the coefficient of $\gamma^\vee$ in $\alpha^\vee$ which is non-negative as $\alpha$ is positive (see \cite[Lemma~3.4]{MR2072765}).
     Therefore, 
     \[
     [C_e] \in \C_{[1]} = \mathrm{span}_{\RR_{\ge 0}}\left\{[C_{e'}] : e'\in E(\G_{[w]})_{[1]}\right\}.
     \]
     This holds for all $e\in E(\G_{[w]})$, so the vertex $[1]\in V(\G_{[w]})$ is special, as desired.
\end{proof}

Combining Lemma~\ref{lem : smooth schubert has special vertex} and Corollary~\ref{cor : key vertex}, we obtain Corollary~\ref{cor : intro Mukai Schubert} from \S\ref{sec : intro}.

\begin{corollary}\label{cor : mukai for schubert}
    Let $X_{[w]}\subseteq G/P$ be a smooth Schubert variety.
    Then $(X_{[w]},\omega|_{X_{[w]}})$ satisfies the generalised symplectic Mukai inequality.
    By the same argument performed in the algebraic category (see \S\ref{sec : alg Mukai}) and using that $(X_{[w]},T_\CC^\text{eff})$ is an algebraic GKM space, the generalised Mukai conjecture holds for $X_{[w]}$.
\end{corollary}

\subsubsection{Example without special vertex}
\label{sec : no special vertex}
Not every Hamiltonian GKM space has a special vertex.
For example, consider the toric blowup
\[
X:=\text{Bl}_{\PP^2\sqcup\PP^2}\PP^5
\]
which is a smooth projective toric variety and therefore a Hamiltonian GKM space.
It has $b_2(X)=3$ and $\find(X)=\psind(X)=2$ and hence satisfies the generalised symplectic Mukai inequality.
However, at each vertex $v$ of its GKM graph $\Gfat$, the set
\[
E(\G)_v = \{e_1,e_2,e_3,e_4,e_5\}
\]
of edges at $v$ has the following properties.
\begin{itemize}
    \item $[C_{e_1}]=[C_{e_2}]$ with $\m(e_1)=\m(e_2)=2$ and $e_1,e_2$ irreducible (see Definition~\ref{def : irreducible}).
    \item $[C_{e_4}]=[C_{e_5}]$ with $\m(e_4)=\m(e_5)=4$ and $e_4,e_5$ are not irreducible.
    \item $\m(e_3)=2$ and $e_3$ is irreducible.
    \item $[C_{e_4}]$ is not minimal among classes at $v$ because $[C_{e_3}]<[C_{e_4}]$. In particular,
\[
    \Cmin_v = \mathrm{span}_{\RR_{\ge 0}}\{[C_{e_1}], [C_{e_3}]\}
\]
    \item $\C_v\neq \C$ even though $\C_v$ is three-dimensional.
\end{itemize}

Thus, no vertex is special and Proposition~\ref{prop : vertex mukai inequality} does not imply the generalised Mukai inequality.
The claims on $X$ above can easily be verified using \textsf{GKMtools.jl} \cite{Holmes_Muratore_2025a}.

\section{Pseudoindex $n$ and complex quadrics}
\label{sec : quadrics}

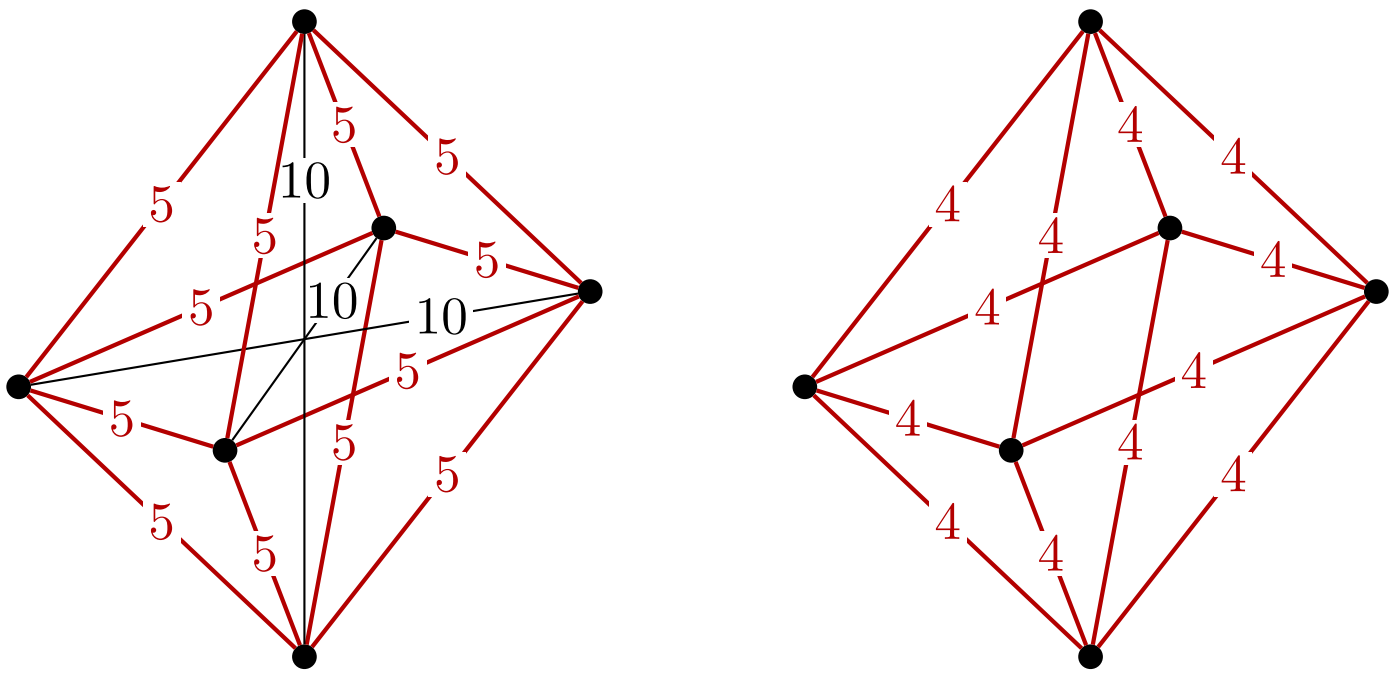
\begin{figure}[h]
    \centering
    
\begin{tikzpicture}[scale=2.1]
\node[draw, circle, fill=black, inner sep=1.5pt] (v1) at (0.0, 0.0) {};
\node[draw, circle, fill=black, inner sep=1.5pt] (v2) at (-0.9, 0.85) {};
\node[draw, circle, fill=black, inner sep=1.5pt] (v3) at (-0.25, 0.65) {};
\node[draw, circle, fill=black, inner sep=1.5pt] (v4) at (0.25, 1.35) {};
\node[draw, circle, fill=black, inner sep=1.5pt] (v5) at (0.9, 1.15) {};
\node[draw, circle, fill=black, inner sep=1.5pt] (v6) at (0.0, 2.0) {};
\draw[red!70!black, thick] (v2) -- node[pos=0.5, fill=white, inner sep=1pt] {\small 5} (v1);
\draw[red!70!black, thick] (v3) -- node[pos=0.5, fill=white, inner sep=1pt] {\small 5} (v1);
\draw[red!70!black, thick] (v3) -- node[pos=0.5, fill=white, inner sep=1pt] {\small 5} (v2);
\draw[red!70!black, thick] (v4) -- node[pos=0.5, fill=white, inner sep=1pt] {\small 5} (v1);
\draw[red!70!black, thick] (v4) -- node[pos=0.5, fill=white, inner sep=1pt] {\small 5} (v2);
\draw (v4) -- node[pos=0.31, fill=white, inner sep=1pt] {\small 10} (v3);
\draw[red!70!black, thick] (v5) -- node[pos=0.5, fill=white, inner sep=1pt] {\small 5} (v1);
\draw (v5) -- node[pos=0.25, fill=white, inner sep=1pt] {\small 10} (v2);
\draw[red!70!black, thick] (v5) -- node[pos=0.5, fill=white, inner sep=1pt] {\small 5} (v3);
\draw[red!70!black, thick] (v5) -- node[pos=0.5, fill=white, inner sep=1pt] {\small 5} (v4);
\draw (v6) -- node[pos=0.24, fill=white, inner sep=1pt] {\small 10} (v1);
\draw[red!70!black, thick] (v6) -- node[pos=0.5, fill=white, inner sep=1pt] {\small 5} (v2);
\draw[red!70!black, thick] (v6) -- node[pos=0.5, fill=white, inner sep=1pt] {\small 5} (v3);
\draw[red!70!black, thick] (v6) -- node[pos=0.5, fill=white, inner sep=1pt] {\small 5} (v4);
\draw[red!70!black, thick] (v6) -- node[pos=0.5, fill=white, inner sep=1pt] {\small 5} (v5);
\end{tikzpicture}
\hspace{1cm}
    \begin{tikzpicture}[scale=2.1]
\node[draw, circle, fill=black, inner sep=1.5pt] (v1) at (0.0, 0.0) {};
\node[draw, circle, fill=black, inner sep=1.5pt] (v2) at (0.9, 1.15) {};
\node[draw, circle, fill=black, inner sep=1.5pt] (v3) at (0.25, 1.35) {};
\node[draw, circle, fill=black, inner sep=1.5pt] (v4) at (-0.25, 0.65) {};
\node[draw, circle, fill=black, inner sep=1.5pt] (v5) at (-0.9, 0.85) {};
\node[draw, circle, fill=black, inner sep=1.5pt] (v6) at (0.0, 2.0) {};
\draw[red!70!black, thick] (v2) -- node[pos=0.5, fill=white, inner sep=1pt] {\small 4} (v1);
\draw[red!70!black, thick] (v3) -- node[pos=0.5, fill=white, inner sep=1pt] {\small 4} (v1);
\draw[red!70!black, thick] (v3) -- node[pos=0.5, fill=white, inner sep=1pt] {\small 4} (v2);
\draw[red!70!black, thick] (v4) -- node[pos=0.5, fill=white, inner sep=1pt] {\small 4} (v1);
\draw[red!70!black, thick] (v4) -- node[pos=0.5, fill=white, inner sep=1pt] {\small 4} (v2);
\draw[red!70!black, thick] (v5) -- node[pos=0.5, fill=white, inner sep=1pt] {\small 4} (v1);
\draw[red!70!black, thick] (v5) -- node[pos=0.5, fill=white, inner sep=1pt] {\small 4} (v3);
\draw[red!70!black, thick] (v5) -- node[pos=0.5, fill=white, inner sep=1pt] {\small 4} (v4);
\draw[red!70!black, thick] (v6) -- node[pos=0.5, fill=white, inner sep=1pt] {\small 4} (v2);
\draw[red!70!black, thick] (v6) -- node[pos=0.5, fill=white, inner sep=1pt] {\small 4} (v3);
\draw[red!70!black, thick] (v6) -- node[pos=0.5, fill=white, inner sep=1pt] {\small 4} (v4);
\draw[red!70!black, thick] (v6) -- node[pos=0.5, fill=white, inner sep=1pt] {\small 4} (v5);
\end{tikzpicture}

    \caption{The GKM graphs of $\cQ_5$ (left) and $\cQ_4$ (right). Each edge $e$ is labelled by its magnitude $\m(e)$ and irreducible edges are coloured red.
    The GKM graph of $\cQ_3$ is depicted in Figure~\ref{fig: local Mukai} (left).}
    \label{fig : quadric GKM graphs}
\end{figure}

In this section, we investigate Hamiltonian GKM spaces of dimension $2n$ and pseudoindex $n$, continuing the analogy with the classical Kobayashi--Ochiai classification of projective space and complex quadrics \cite{MR316745, MR1929792, MR1922103, MR3098425, MR3592468} started in Theorem~\ref{thm: index n+1 with proof}.
For $n=2N$ and $n=2N-1$, the complex quadric $\cQ_n$ is given by
\begin{align*}
        \cQ_{2N} &= \left\{[z_1:\cdots:z_{2N+2}] \in \CC\PP^{2N+1} : \sum_{i=1}^{N+1} z_iz_{2N+3-i} = 0
        \right\},
        \\
        \cQ_{2N-1} &=\left\{[z_1:\cdots:z_{2N+1}] \in \CC\PP^{2N} : \sum_{i=1}^{N+1} z_iz_{2N+1-i} + z_{2N+1}^2= 0
        \right\}.
\end{align*}
It carries a torus action by $T=(S^1)^{N+1}$ and $T=(S^1)^N$, respectively, given by
\begin{align*}
    (t_1,\dots,t_{N+1})\cdot [z_1:\cdots:z_{2N+2}] &:= [t_1z_1:\cdots t_{N+1}z_{N+1}:t_{N+1}^{-1}z_{N+2}:\cdots :t_1^{-1}z_{2N+2}],
    \\
    (t_1,\dots,t_N)\cdot [z_1:\cdots:z_{2N+1}] &:= [t_1z_1:\cdots t_{N}z_{N}:t_{N}^{-1}z_{N+1}:\cdots :t_1^{-1}z_{2N}:z_{2N+1}]
\end{align*}
After replacing $T$ by its effectively acting quotient torus, $\cQ_n$ becomes a Hamiltonian\footnote{
Indeed, the $T$-action is the restriction of an algebraic GKM action of a complex torus on $\cQ_n$, so the $T$-space $\cQ_n$ has the structure of a Hamiltonian GKM space (cf. \S\ref{sec : alg Mukai}).
} GKM space with respect to the given $T$-action \cite{Kuroki2025, KurokiPaul2026}.
In loc.\ cit., an explicit description of its GKM graph $\Gfat_{\cQ_n}$ is given (see Figure~\ref{fig : quadric GKM graphs} for two examples).

The goal of this section is to prove the following theorem, which constitutes the second half of Corollary~\ref{cor : index}.

\begin{theorem}\label{thm : quadrics}
    Let $(X,\omega,T)$ be a Hamiltonian GKM space of dimension $2n>0$ and pseudoindex $\psind(X,\omega)=n$.
    Then the GKM graph of $X$ is isomorphic to a projection of the GKM graph $\Gfat_{\cQ_n}$ of the complex quadric $\cQ_n$.
\end{theorem}

To prove Theorem~\ref{thm : quadrics}, we will first give a criterion in \S\ref{sec : GKM graph of quadrics} for a GKM graph to be a projection of $\Gfat_{\cQ_n}$.
In \S\ref{sec : GKM graph of pseudoindex n}, we prove that this criterion is indeed satisfied by any Hamiltonian GKM space as in the theorem.

\subsection{The GKM graphs of complex quadrics}\label{sec : GKM graph of quadrics}

\begin{proposition}\label{prop : quadric criterion}
    A GKM graph $\Gfat$ of valency $n$ is isomorphic to a projection of $\Gfat_{\cQ_n}$ if and only if the following condition holds.
     \begin{itemize}
        \item If $2\mid n$ then $\G$ is a graph on $n+2$ vertices such that for each $v\in V(\G)$ there exists a unique $\overline{v}\in V(\G)$ with $\{v,\overline{v}\}\notin E(\G)$.
        Moreover:
        \begin{itemize}
            \item For each $\{v,w\}\in E(\G)$, we have $\alpha(vw)=\alpha(\overline{w}\,\overline{v})$.
            \item For each $u,v,w\in V(\G)$ forming a triangle, we have $\alpha(uv)+\alpha(vw)=\alpha(uw)$.
        \end{itemize}

        \item If $2\nmid n$ then $\G$ is a complete graph on $n+1$ vertices and there is a fixed-point free involution $v\mapsto \overline{v}$ on $V(\G)$ such that:
        \begin{itemize}
            \item For each $\{v,w\}\in E(\G)$, we have $\alpha(vw)=\alpha(\overline{w}\,\overline{v})$.
            \item For each $v\in V(\G)$ and $w\in V(\G)\setminus\{v,\overline{v}\}$, we have $\alpha(wv)+\alpha(v\overline{w})=2\alpha(w\overline{w})$.
            \item For any distinct $u,v,w\in V(\G)$ with $\{\overline{u},\overline{v},\overline{w}\}\cap \{u,v,w\}=\emptyset$, we have $\alpha(uv)+\alpha(vw)=\alpha(uw)$.
        \end{itemize}
    \end{itemize}
\end{proposition}
\begin{proof}
    It follows immediately from \cite[(2.3)]{Kuroki2025} and \cite[\S2.1]{KurokiPaul2026} that $\Gfat_{\cQ_n}$ satisfies the given condition, which is invariant under isomorphisms and projections of GKM graphs.
    Therefore, any GKM graph isomorphic to a projection of $\Gfat_{\cQ_n}$ also satisfies the condition.

    Conversely, let $\Gfat=(\G,\alpha)$ satisfy the condition.
    We need to show that $\Gfat$ is isomorphic to a projection of $\Gfat_{\cQ_n}=(\G_{\cQ_n},\alpha_{\cQ_n})$.
    First, note that the condition completely determines the isomorphism type of the underlying graph, so without loss of generality $\G=\G_{\cQ_n}$.
    Let $T,T_{\cQ_n}$ be the tori acting on $X, \cQ_n$, respectively.
    It remains to show that there exists a surjection $\psi\colon\Lie(T_{\cQ_n})^\vee\rightarrow \Lie(T)^\vee$ such that, for all $vw\in E(\G)^\pm$,
    \begin{equation}\label{eqn : quadric projection}
    \alpha(vw)=\psi(\alpha_{\cQ_n}(vw)).
    \end{equation}
    Let $V := \RR^{E(\G)^\pm}$ with basis vectors $(\mathbf{1}_\epsilon)_{\epsilon\in E(\G)^\pm}$.
    Let $S\le V$ be the subspace spanned by $ \mathbf{1}_\epsilon + \mathbf{1}_{\overline{\epsilon}}$ for each $\epsilon\in E(\G)^\pm$ and the relations given by the condition.
    For example, the condition $\alpha(uv)+\alpha(vw)=\alpha(uw)$ translates to the element $\mathbf{1}_{uv}+\mathbf{1}_{vw}-\mathbf{1}_{uw}$.
    Then
    for each $(\alpha',T')\in \{(\alpha, T),(\alpha_{\cQ_n}, T_{\cQ_n})\}$, there exists a unique linear map $\widetilde{\alpha'}$ making the following diagram commute.
    \begin{equation}\label{eqn : quadric eqn diagram}
        \begin{tikzcd}
            V\ar["\alpha'", twoheadrightarrow]{r}\ar[twoheadrightarrow]{d} & \Lie(T')^\vee
            \\
            V/S\ar["\widetilde{\alpha'}"', dashed, twoheadrightarrow]{ur}
        \end{tikzcd}
    \end{equation}
    Here, surjectivity of $\widetilde{\alpha'}$ follows from surjectivity of $\alpha'$, which holds since the $T'$-action is assumed to be effective.

    \textit{Claim: $\widetilde{\alpha_{\cQ_n}}$ is an isomorphism.} 
        
    \textit{Proof of Claim:}
    By surjectivity, it suffices to prove that $\mathrm{rank}(T_{\cQ_n})\ge \dim V/S$.
    By \cite{Kuroki2025, KurokiPaul2026}, we have
    \[
    \mathrm{rank}(T_{\cQ_n})= \left\lceil \frac{n+1}{2} \right\rceil =: r_n.
    \]
    On the other hand, we may write
    \[
    V(\G) = \{1,\dots,r_n,\overline{1},\dots,\overline{r_n}\}
    \]
    and apply the given conditions to see that $V/S$ is spanned by the images of $\{\mathbf{1}_\epsilon : \epsilon\in E_0\}$, where $E_0\subset E(\G)^\pm$ is defined as
    \[
    E_0 := \begin{cases}
        \{12, 13, \dots, 1 r_n, 1\overline{1}\}
        & 2\nmid n,
        \\
        \{12,13,\dots,1r_n, 1\overline{2}\}
        & 2\mid n.
    \end{cases}
    \]
    Thus, we obtain $\mathrm{dim}V/S\le |E_0| = r_n$, as desired, which proves the Claim.

    \textit{Conclusion:}
    One easily uses diagram \eqref{eqn : quadric eqn diagram} to verify that a projection satisfying \eqref{eqn : quadric projection} is given by $\psi:= \widetilde{\alpha}\circ\widetilde{\alpha_{\cQ_n}}^{-1}$.
\end{proof}

\subsection{Hamiltonian GKM graphs of pseudoindex $n$}\label{sec : GKM graph of pseudoindex n}

Throughout this section, we let $(X,\omega,T)$ be a Hamiltonian GKM space of dimension $2n>0$ and assume that its pseudoindex is $\psind:=\psind(X,\omega)=n$.

\begin{lemma}\label{lem : quadric n ge 3, b2 is 1}
     Either $n\ge 3$ and $b_2(X)=1$, or $n=2$ and Theorem~\ref{thm : quadrics} holds.
\end{lemma}
\begin{proof}
    We have
\[
\psind = n\Longrightarrow
\begin{cases}
    \psind < (n+2)/2
    & \text{if }n=1,
    \\
    \psind = (n+2)/2
    & \text{if } n=2,
    \\
    \psind > (n+2)/2
    &\text{if }n\ge 3.
\end{cases}
\]

If $n=1$ then $\Gfat$ is isomorphic to the GKM graph of $\CC\PP^1$, which is impossible since its index is $\find=2$ which does not divide $\psind=n=1$, contradicting Lemma~\ref{lem : find and psind from G}.

If $n=2$ then by Theorem~\ref{thm : Mukai large index} either $b_2(X)=1$, or $b_2(X)=2$ and $\Gfat$ is isomorphic to a projection of the GKM graph of $\CC\PP^1\times\CC\PP^1$.
One easily checks using Example~\ref{ex: GKM graph of CPn} and Proposition~\ref{prop : quadric criterion} that the GKM graphs of $\CC\PP^1\times\CC\PP^1$ and $\cQ_2$ are isomorphic, so Theorem~\ref{thm : quadrics} holds in this case.
On the other hand, if $b_2(X)=1$ then $|V(\G)|=\chi(X)=3$. Using the existence of a compatible connection on $\Gfat$, one easily sees that $\Gfat$ is isomorphic to a projection of the GKM graph of $\CC\PP^2$.
But then its index is $\find=3$, which does not divide $\psind=n=2$, contradicting again Lemma~\ref{lem : find and psind from G}.

Finally, if $n\ge3$ then $b_2(X)=1$ by Theorem~\ref{thm : Mukai large index}.
\end{proof}

By the lemma above, we may assume from now on that $n\ge 3$ and $b_2(X)=1$.
In particular, since $c_1(TX)\neq 0$ (cf.\ Lemma~\ref{lem : exists positive magnitude edge}) the curve class $[C_e]$ of an edge $e\in E(\G)$ is uniquely determined by $\m(e)$.

\begin{lemma}\label{lem : quadric irred}
    An edge $e\in E(\G)$ is irreducible if and only if $\m(e)=n$.
    Moreover, any two vertices of $\G$ are connected by a path of irreducible edges.
\end{lemma}
\begin{proof}
    Since $n\ge 3$ gives $\psind=n\ge (n+2)/2$, Lemma~\ref{lem : large psind primitivity} implies that an edge $e$ is irreducible if and only if $\m(e)\le n+1$.
    By Lemma~\ref{lem : find and psind from G}, we also have $\m(e)\ge \psind =n$, so $e$ is irreducible if and only if $\m(e)\in\{n,n+1\}$.

    It remains to prove that there is no irreducible edge $e$ with $\m(e)=n+1$.
    Indeed,
    if $e=vw$ is an irreducible edge with $\m(e)=n+1$ then $[C_e]$ is minimal among curve classes of edges at $v$.
    Thus, by Theorem~\ref{thm: local Mukai}, all edges at $v$ have magnitude $n+1$ and are therefore irreducible.
    Iterating this argument shows that all $e\in E(\G)$ satisfy $\m(e)=n+1$.
    In particular, $\find(X,\omega)=n+1$ by Lemma~\ref{lem : find and psind from G}.
    On the other hand, we have $\find(X,\omega)\mid\psind$ by the same lemma, which contradicts $\psind=n$.
    Therefore, there is no edge $e\in E(\G)$ with $\m(e)=n+1$, which finishes the proof.

    The \emph{moreover} part follows from Lemma~\ref{lem : large psind primitive connectedness} since we have $n\ge 2$ and $\psind=n$.
\end{proof}

\begin{figure}[htb]
    \centering
    \begin{tikzpicture}
\node at (0,-2.45) {\textup{(T1)}};

\node[draw, circle, fill=black, inner sep=1.5pt, label=left:$v$] (TOnev) at (-0.75,0) {};
\node[draw, circle, fill=black, inner sep=1.5pt, label=right:$w$] (TOneW) at (0.75,0) {};
\node[draw, circle, fill=black, inner sep=1.5pt, label=above:$u_1$] (TOneuOne) at (-1.5,1.65) {};
\node[draw, circle, fill=black, inner sep=1.5pt, label=above:$u_{n-3}$] (TOneuLast) at (1.5,1.65) {};
\node at (0,1.65) {$\cdots$};
\node[draw, circle, fill=black, inner sep=1.5pt, label=below:$x$] (TOnex) at (-1.5,-1.65) {};
\node[draw, circle, fill=black, inner sep=1.5pt, label=below:$y$] (Toney) at (1.5,-1.65) {};

\draw[red!70!black, thick] (TOnev) --
    node[pos=0.50, fill=white, inner sep=1pt] {\small $n$} (TOneW);

\draw[red!70!black, thick] (TOnev) -- node[pos=0.52, fill=white, inner sep=1pt] {\small $n$} (TOnex);
\draw[red!70!black, thick] (TOnev) -- node[pos=0.50, fill=white, inner sep=1pt] {\small $n$} (TOneuOne);
\draw[red!70!black, thick] (TOnev) -- node[pos=0.50, fill=white, inner sep=1pt] {\small $n$} (TOneuLast);
\draw (TOnev) -- node[pos=0.52, fill=white, inner sep=1pt] {\small $2n$} (Toney);

\draw (TOneW) -- node[pos=0.52, fill=white, inner sep=1pt] {\small $2n$} (TOnex);
\draw[red!70!black, thick] (TOneW) -- node[pos=0.50, fill=white, inner sep=1pt] {\small $n$} (TOneuOne);
\draw[red!70!black, thick] (TOneW) -- node[pos=0.50, fill=white, inner sep=1pt] {\small $n$} (TOneuLast);
\draw[red!70!black, thick] (TOneW) -- node[pos=0.52, fill=white, inner sep=1pt] {\small $n$} (Toney);
    \end{tikzpicture}
    \hspace{0.5cm}
    \begin{tikzpicture}
\node at (0,-2.45) {\textup{(T2)}};

\node[draw, circle, fill=black, inner sep=1.5pt, label=left:$v$] (Ttwov) at (-0.75,0) {};
\node[draw, circle, fill=black, inner sep=1.5pt, label=right:$w$] (TtwoW) at (0.75,0) {};
\node[draw, circle, fill=black, inner sep=1.5pt, label=above:$u_1$] (TtwouOne) at (-1.5,1.65) {};
\node[draw, circle, fill=black, inner sep=1.5pt, label=above:$u_{n-2}$] (TtwouLast) at (1.5,1.65) {};
\node at (0,1.65) {$\cdots$};
\node[draw, circle, fill=black, inner sep=1.5pt, label=below:$x$] (Ttwox) at (-1.5,-1.65) {};
\node[draw, circle, fill=black, inner sep=1.5pt, label=below:$y$] (Ttwoy) at (1.5,-1.65) {};

\draw[red!70!black, thick] (Ttwov) --
    node[pos=0.50, fill=white, inner sep=1pt] {\small $n$} (TtwoW);

\draw[red!70!black, thick] (Ttwov) -- node[pos=0.50, fill=white, inner sep=1pt] {\small $n$} (TtwouOne);
\draw[red!70!black, thick] (TtwoW) -- node[pos=0.50, fill=white, inner sep=1pt] {\small $n$} (TtwouOne);
\draw[red!70!black, thick] (Ttwov) -- node[pos=0.50, fill=white, inner sep=1pt] {\small $n$} (TtwouLast);
\draw[red!70!black, thick] (TtwoW) -- node[pos=0.50, fill=white, inner sep=1pt] {\small $n$} (TtwouLast);

\draw (Ttwov) -- (Ttwox);
\draw (TtwoW) -- (Ttwoy);
    \end{tikzpicture}

    \caption{Local model of an irreducible edge $e=\{v,w\}$ of type (T1) and (T2) as in Lemma~\ref{lem : quadrics T1 T2}. Edge labels indicate the magnitude of that edge.}
    \label{fig: quadrics T1 T2}
\end{figure}
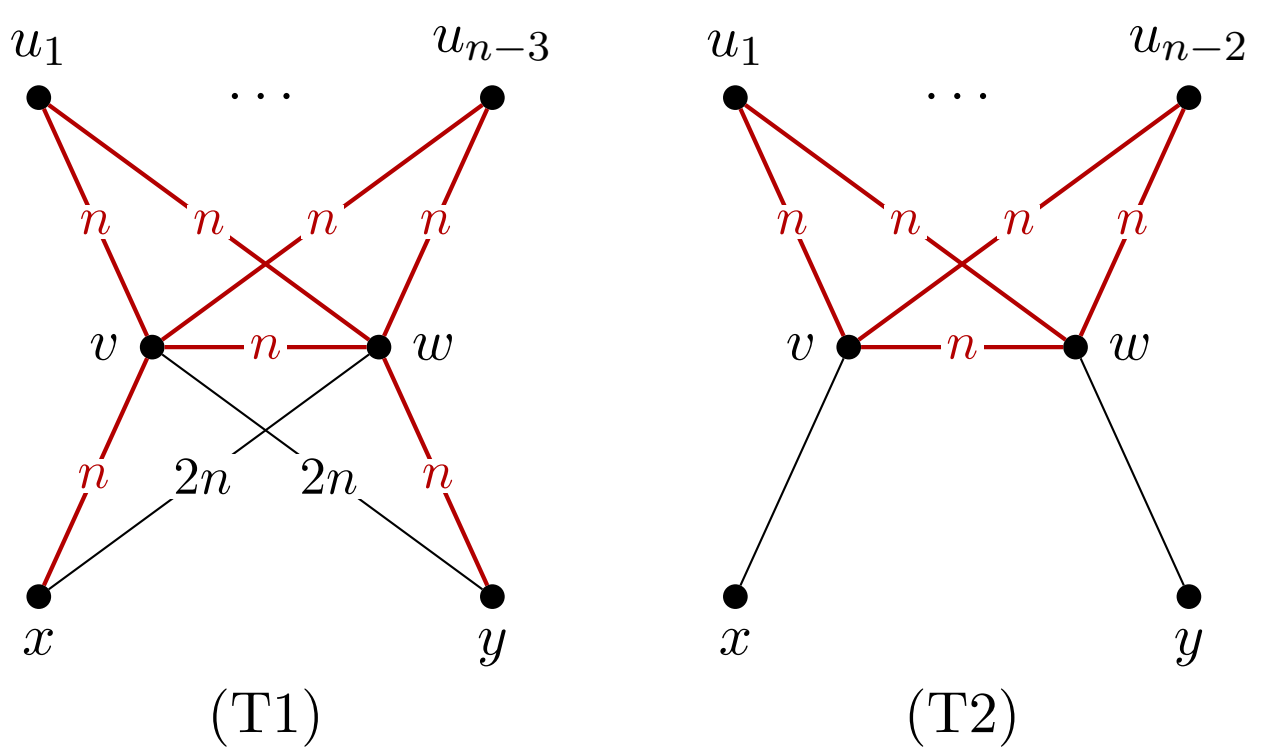

Next, we show that there are exactly two possible types of irreducible edges, called (T1) and (T2), as illustrated in Figure~\ref{fig: quadrics T1 T2}.

\begin{lemma}\label{lem : quadrics T1 T2}
    Let $e=\{v,w\}\in E(\G)$ with $\m(e)=n$.
    Then precisely one of the following two cases holds.
    \begin{itemize}
        \item[(T1)] There exist distinct $x,y,u_1,\dots,u_{n-3}\in V(\G)\setminus\{v,w\}$ such that
        \begin{enumerate}
            \item each of $x,y,u_1,\dots,u_{n-3}$ is adjacent to $v$ and $w$ in $\G$,
            \item $\m(vu_i)=\m(wu_i)=n$ for all $i\in\{1,\dots,n-3\}$,
            \item $\m(vx)=\m(wy)=n$, and
            \item $\m(vy)=\m(wx)=2n$.
        \end{enumerate}
        \item[(T2)] There exist distinct $x,y,u_1,\dots,u_{n-2}\in V(\G)\setminus\{v,w\}$ such that
        \begin{enumerate}
            \item each of $u_1,\dots,u_{n-2}$ is adjacent to $v$ and $w$ in $\Gfat$,
            \item $\m(vu_i)=\m(wu_i)=n$ for all $i\in\{1,\dots,n-2\}$, and
            \item $\{v,x\},\{w,y\}\in E(\G)$.
        \end{enumerate}
    \end{itemize}
\end{lemma}
\begin{proof}
    By Lemma~\ref{lem : quadric irred}, $e$ is irreducible and hence $[C_e]$ is minimal among edge curve classes at $v$ and $w$.
    By Theorem~\ref{thm: local Mukai}, there exist $v_1,\dots,v_{n-2},w_1,\dots,w_{n-2}\in V(\G)\setminus\{v,w\}$ such that for all $i\in\{1,\dots,n-2\}$, we have
    \begin{itemize}
        \item $vv_i,vw_i,ww_i,wv_i\in E(\G)^\pm$,
        \item $[C_{vv_i}]=[C_e]=[C_{ww_i}]$, and hence
        \item $\m(vv_i)=\m(ww_i)=n$.
    \end{itemize}
    Since $\G$ is $n$-valent, we have $\left|\{v_1,\dots,v_{n-2}\}\cap\{w_1,\dots,w_{n-2}\}\right|\in\{n-3,n-2\}$, which corresponds precisely to  cases $(T1)$ and $(T2)$, respectively, where the $u_i$ are the elements in the intersection.

    It remains to show that in case (T1), we have $\m(wx)=\m(vy)=2n$, so let us assume to be in case (T1) from now on.
    Each edge curve class is a sum of classes of irreducible edges, which all have magnitude $n$ by Lemma~\ref{lem : quadric irred}.
    Thus, there exist $a,b\in\ZZ$ such that $\m(vy)=an$ and $\m(wx)=bn$.
    Since $b_2(X)=1$, there is $\lambda\in\RR$ such that $c_1(TX)=\lambda[\omega]$.
    Since $\m(e)=n>0$, we have $\lambda>0$ and hence $a,b>0$ (since $\omega(C_{e'})>0$ for all $e'\in E(\G)$).

    By $b_2(X)=1$ and Proposition~\ref{prop : curve class relations}, we have
    \begin{align}\label{eqn : quadric triangle 1}
        \alpha(wv) + \alpha(vx) &= b\alpha(wx),\\\label{eqn : quadric triangle 2}
        \alpha(vw) + \alpha(wy) &= a\alpha(vy),\\\label{eqn : quadric triangle 3}
        \alpha(vu_i) + \alpha(u_iw) &= \alpha(vw)\,\,\forall i\in\{1,\dots,n-3\}.
    \end{align}
    Combining equation~\eqref{eqn : quadric triangle 3} with equation~\eqref{eqn : magnitude from G} for $\m(e)=n$ gives
    \begin{equation}\label{eqn : quadric con xy}
    \alpha(vy) + \alpha(vx) - \alpha(wy)-\alpha(wx) = \alpha(vw).
    \end{equation}
    Substituting equations~\eqref{eqn : quadric triangle 1} and \eqref{eqn : quadric triangle 2} into equation~\eqref{eqn : quadric con xy} gives
    \[
    \alpha(vw) + (1-a)\alpha(wy)+a(b-1)\alpha(wx) = 0.
    \]
    Since $a=1$ or $b=1$ would contradict 2-independence of $\Gfat$ at $w$, we must have $a\ge 2$ and $b\ge 2$.
    The existence of a compatible connection (see Definition~\ref{def : connection}) for $\Gfat$ at $vw$ together with equation~\eqref{eqn : quadric triangle 3} gives that
    \[
    \{\alpha(vx), \alpha(vy)\} \equiv \{\alpha(wx), \alpha(wy)\} \hspace{3mm}\text{mod } \ZZ\alpha(vw).
    \]
    But by equation~\eqref{eqn : quadric triangle 1} and $b\ge 2$ we have $\alpha(vx)\not\equiv\alpha(wx)$ and hence $\alpha(vx)\equiv \alpha(wy)$ and $\alpha(vy)\equiv \alpha(wx)$ mod $\ZZ\alpha(vw)$.
    Equations~\eqref{eqn : quadric triangle 1} and \eqref{eqn : quadric triangle 2} mod $\ZZ\alpha(vw)$ therefore give $a=b$.
    Substituting equations~\eqref{eqn : quadric triangle 1}, \eqref{eqn : quadric triangle 2} into equation~\eqref{eqn : quadric con xy} once more gives
    \[
        (a-1)(\alpha(vy)-\alpha(wx)) = \alpha(vw).
    \]
    Since $\alpha(vy)-\alpha(wx)$ is an integral multiple of $\alpha(vw)$ and $\alpha(vw)\neq 0$, we must have $a\in\{0,2\}$.
    Since $a\ge 2$, we finally obtain $a=b=2$, as desired.
\end{proof}

\begin{lemma}\label{lem : quadric vx || wy}
    Let $e=\{v,w\}\in E(\G)$ with $\m(e)=n$.
    Then in the notation of Lemma~\ref{lem : quadrics T1 T2}, we have $\alpha(vx) = \alpha(wy)$.
\end{lemma}
\begin{proof}
    In case (T1), the proof of Lemma~\ref{lem : quadrics T1 T2} shows that equations~\eqref{eqn : quadric triangle 1},\eqref{eqn : quadric triangle 2} hold with $a=b=2$.
    Substituting this into equation~\eqref{eqn : quadric con xy} gives the desired $\alpha(vx)-\alpha(wy)=0$.

    In case (T2), Lemma~\ref{lem : quadrics T1 T2}, Proposition~\ref{prop : curve class relations}, and $b_2(X)=1$ similarly give
    \[
    \alpha(vu_i) + \alpha(u_iw) = \alpha(vw)\,\,\forall i\in\{1,\dots,n-2\}
    \]
    and hence equation~\eqref{eqn : magnitude from G} for $\m(vw)=n$ gives $\alpha(vx)-\alpha(wy)=0$.
\end{proof}

\begin{lemma}\label{lem : quadrics all T1 or all T2}
    Either all irreducible edges have type (T1), or all irreducible edges have type (T2) (in the notation of Lemma~\ref{lem : quadrics T1 T2}). 
\end{lemma}
\begin{proof}
    Recall Lemma~\ref{lem : quadric irred} that $e\in E(\G)$ is irreducible if and only if $\m(e)=n$.
    By Lemma~\ref{lem : quadrics T1 T2}, every irreducible edge has type (T1) or (T2), so it suffices to prove that if one irreducible edge has type (T1) then all irreducible edges have type (T1).
    By Lemma~\ref{lem : quadric irred}, any two vertices are connected by a path of irreducible edges, so it suffices to prove that if $e=\{v,w\}\in E(\G)$ is an irreducible edge of type (T1) then all irreducible edges at $w$ have type (T1).
    To do this, observe the following in the notation of Lemma~\ref{lem : quadrics T1 T2} (cf. Figure~\ref{fig: quadrics T1 T2}).

    Since $\m(wy)=n$ there are at least $n-3$ common neighbours $z$ of $w$ and $y$ with $\m(wz)=\m(yz)=n$.
        Since $\m(wx)=\m(vy)=2n\neq n$, the only possible such $z$ are $u_1,\dots,u_{n-3}$.
        Since in type (T2) there would be $n-2$ such $z$, it follows that $wy$ has type (T1) and $yu_1,\dots,yu_{n-3}\in E(\G)^\pm$ with $\m(yu_i)=n$ for all $i\in\{1,\dots,n-3\}$.

        Analogously, $vx\in E(\G)^\pm$ is irreducible of type (T1) and $xu_i\in E(\G)^\pm$ with $\m(xu_i)=n$ for all $i\in \{1,\dots,n-3\}$.

        It remains to show that $wu_i$ has type (T1) for all $i\in\{1,\dots,n-3\}$, so assume for contradiction that $wu_{i}$ has type (T2) for some $i$.
        Then there are $n-2$ common neighbours $z$ of $w$ and $u_i$ with $\m(wz)=\m(u_iz)=n$.
        Since $\m(wx)=2n\neq n$, the only possible such $z$ are $v,y,u_1,\dots,\widehat{u_i},\dots,u_{n-3}$, where $\widehat{u_i}$ means that $u_i$ is omitted.
        Thus, all of these are neighbours of $w$ and $u_i$.
        Moreover, $xu_i\in E(\G)^\pm$ by the previous paragraph.
        In total, all neighbours of $w$ except $u_i$ are also neighbours of $u_i$, so $wu_i$ cannot be of type (T2) and must be of type (T1), as desired.
\end{proof}

From now on, we will say to be \emph{in case (T1)} if all irreducible edges of $\G$ are of type (T1), and \emph{in case (T2)} if all irreducible edges of $\G$ are of type (T2).
By Lemma~\ref{lem : quadrics all T1 or all T2}, precisely one of these two cases applies.

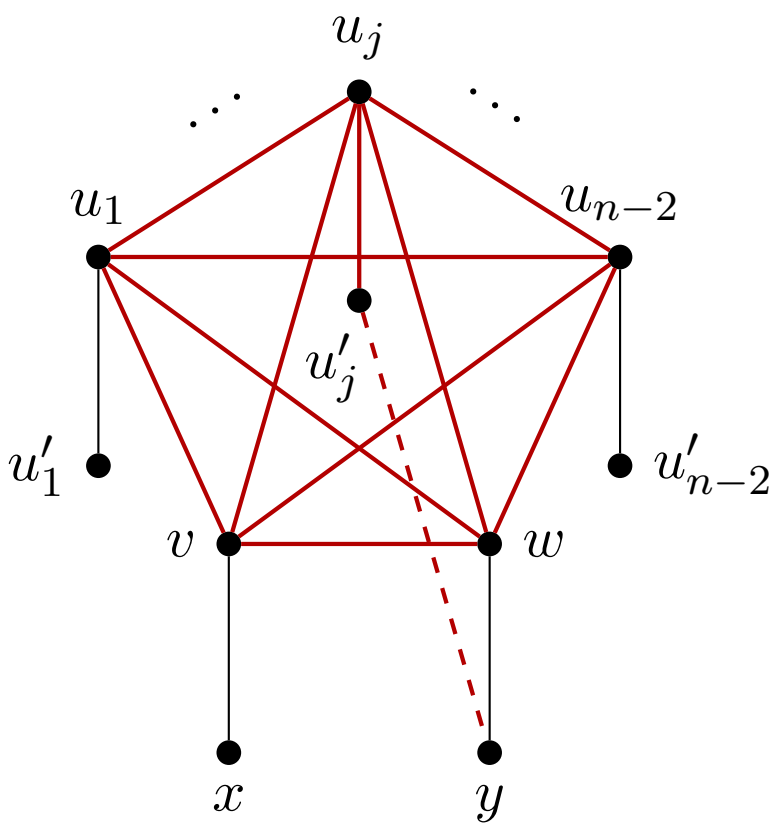
\begin{figure}[h]
    \centering
    \begin{tikzpicture}[scale=1.15]
\node[draw, circle, fill=black, inner sep=1.5pt, label=left:$v$] (v) at (-0.75,0) {};
\node[draw, circle, fill=black, inner sep=1.5pt, label=right:$w$] (w) at (0.75,0) {};
\node[draw, circle, fill=black, inner sep=1.5pt, label=above:$u_1$] (uOne) at (-1.5,1.65) {};
\node[draw, circle, fill=black, inner sep=1.5pt, label=above:$u_{n-2}$] (uLast) at (1.5,1.65) {};
\node[draw, circle, fill=black, inner sep=1.5pt, label=above:$u_j$] (uJ) at (0,2.6) {};
\node[rotate=32.3] at (-0.8,2.5) {$\cdots$};
\node[rotate=-32.3] at (0.8,2.5) {$\cdots$};
\node[draw, circle, fill=black, inner sep=1.5pt, label=below:$x$] (x) at (-0.75,-1.2) {};
\node[draw, circle, fill=black, inner sep=1.5pt, label=below:$y$] (y) at (0.75,-1.2) {};
\node[draw, circle, fill=black, inner sep=1.5pt, label=left:$u_1'$] (uOneP) at (-1.5,0.45) {};
\node[draw, circle, fill=black, inner sep=1.5pt, label=right:$u_{n-2}'$] (uLastP) at (1.5,0.45) {};
\node[draw, circle, fill=black, inner sep=1.5pt, label={[xshift=-5pt]below:$u_j'$}] (uJP) at (0,1.4) {};

\draw[red!70!black, thick] (v) -- (w);

\foreach \u in {uOne, uLast} {
  \draw[red!70!black, thick] (v) -- (\u);
  \draw[red!70!black, thick] (w) -- (\u);
}
\draw[red!70!black, thick] (v) -- (uJ);
\draw[red!70!black, thick] (w) -- (uJ);

\draw[red!70!black, thick] (uOne) -- (uLast);
\draw[red!70!black, thick] (uOne) -- (uJ);
\draw[red!70!black, thick] (uJ) -- (uLast);

\draw[red!70!black, thick] (uJ) -- (uJP);

\draw[red!70!black, thick, dashed] (uJP) -- (y);

\draw (v) -- (x);
\draw (w) -- (y);
\draw (uOne) -- (uOneP);
\draw (uLast) -- (uLastP);
\end{tikzpicture}
    \caption{Illustration of the proof of Lemma~\ref{lem : quadric T2 all irred}. The red (thick) edges, including the dashed one, all have magnitude $n$.}
    \label{fig : quadric T2 all irred}
\end{figure}

\begin{lemma}\label{lem : quadric T2 all irred}
    In case (T2), all edges are irreducible.
\end{lemma}
\begin{proof}
    By Lemma~\ref{lem : quadric irred},
    every vertex is adjacent to an irreducible edge.
    Since we are in case (T2), it therefore suffices (in the notation of Lemma~\ref{lem : quadrics T1 T2} and Figure~\ref{fig: quadrics T1 T2}) to prove that $\m(wy)=n$ for any irreducible $e=\{v,w\}\in E(\G)$.

    Assume for contradiction that $\m(wy)\neq n$, so in particular $\m(wy)\ge 2n$ by Lemma~\ref{lem : quadric irred}.
    Let $i\in\{1,\dots,n-2\}$.
    Since $wu_i$ is irreducible and of type (T2), there are $n-2$ common neighbours $z$ of $w$ and $u_i$ such that $\m(wz)=\m(u_iz)=n$.
    Since the only options for $z$ are the $n-2$ vertices $u_1,\dots,\widehat{u_i},\dots,u_{n-2},v$, we have $u_iu_j\in E(\G)^\pm$ and $\m(u_iu_j)=n$ for all $j\neq i$.
    Thus, the subgraph induced by $S:=\{v,w,u_1,\dots,u_{n-2}\}$ is a complete graph of valency $n-1$ in which every edge has magnitude $n$.

    Since each $wu_i$ is of type (T2), there exists $u_i'\in V(\G)$ adjacent to $u_i$ but no other vertex of $u_1,\dots,u_{n-2},v,w$ with $\alpha(u_iu_i')=\alpha(wy)=\alpha(vx)$ (cf. Lemma~\ref{lem : quadric vx || wy}).
    By 2-independence, we further get $u_i'\notin\{x,y\}$ and all $u_i'$ are pairwise distinct as $i\in\{1,\dots,n-2\}$ varies.

    By Lemma~\ref{lem : quadric irred}, there exists a path from $w$ to $y$ consisting only of irreducible edges.
    Therefore, there either exists $j\in\{1,\dots,n-2\}$ such that $\m(u_ju_j')=n$, or $\m(vx)=n$.
    We assume the former holds, and note that the latter case can be treated by exactly the same argument after swapping $(u_j,u_j')$ and $(v,x)$ in the following paragraph.
    
    Since $\m(u_ju_j')=n$, the path $P=wu_ju_j'$ has $\m(P)=2n$.
    Since $n\ge 3$, we have $\m(P)>n+1$.
    By Theorem~\ref{thm: unique min path restated}, there exists another path $Q$ from $w$ to $u_j'$ with $\m(Q)\le 2n$.
    In particular, $\len(Q)\le 2$ and hence $\len(Q)=2$ since $wu_j'\notin E(\G)^\pm$. 
    So there exists $z\in V(\G)\setminus\{u_j\}$ such that $Q=wzu_j'$ is a path on $\G$ with $\m(wz)=\m(zu_j')=n$.
    Since neither $v$ nor any $u_i$ with $i\neq j$ is adjacent to $u_j'$, we must have $z=y$.
    Thus, $\m(wy)=n$, contradicting our assumption and finishing the proof by contradiction.
\end{proof}

\begin{lemma}\label{lem : quadrics T2 rectangle}
    In case $(T2)$ and in the notation of Lemma~\ref{lem : quadrics T1 T2} (see Figure~\ref{fig: quadrics T1 T2}), each edge $e=\{v,w\}\in E(\G)$ satisfies $\{x,y\}\in E(\G)$ and $\alpha(xy)=\alpha(vw)$.
    Moreover, each common neighbour of $v$ and $w$ is adjacent to $x$ and $y$.
\end{lemma}
\begin{proof}
    Recall that by Lemma~\ref{lem : quadrics all T1 or all T2} and Lemma~\ref{lem : quadric T2 all irred}, all edges $e=\{v,w\}$ are irreducible and of type (T2), so the local model of $e$ is as in the right of Figure~\ref{fig: quadrics T1 T2}.
    Note also that $y$ is characterised as the unique neighbour of $w$ not adjacent to $v$ and distinct from $v$.
    Since $wy$ is again of type (T2), we can construct an infinite path $P=v_0v_1v_2\dots$ as follows:
    \begin{itemize}
        \item $v_0=v$,
        \item $v_1=w$, and
        \item for any $i\ge0$, $v_{i+2}$ is the unique neighbour of $v_{i+1}$ not adjacent to $v_i$ and distinct from $v_i$.
    \end{itemize}
    This is well-defined since all edges are irreducible and of type (T2).

    \textit{Claim 1: Each $u_j$ is adjacent to all vertices of $P$.}

    \textit{Proof of Claim 1:}
    We prove that $v_i$ is adjacent to each $u_j$ by induction on $i\ge 0$.
    The claim clearly holds for $i\in\{0,1\}$ by definition of $u_1,\dots,u_{n-2}$ in Lemma~\ref{lem : quadrics T1 T2}.
    For the induction step, assume $u_j$ is incident to $v_i$ and $v_{i+1}$.
    Since $v_{i+1}v_{i+2}$ is irreducible of type (T2), and since $v_{i+2}v_i\notin E(\G)^\pm$, $v_i$ is the unique neighbour of $v_{i+1}$ not adjacent to $v_{i+2}$.
    Since $u_j\neq v_i,v_{i+1}$ and $u_j$ is adjacent to $v_{i+1}$, it is therefore also adjacent to $v_{i+2}$, so the claim is proved.

    \textit{Claim 2: We have $\alpha(v_{i}v_{i+1})=-\alpha(v_{i+2}v_{i+3})$ for all $i\ge 0$.}

    \textit{Proof of Claim 2:}
    The edge $v_{i+1}v_{i+2}$ is irreducible of type (T2) and $v_iv_{i+2},v_{i+1}v_{i+3}\notin E(\G)^\pm$ by construction.
    Thus, by Lemma~\ref{lem : quadric vx || wy} we have $\alpha(v_{i+1}v_i)=\alpha(v_{i+2}v_{i+3})$ as desired.

    \textit{Claim 3: We have $v_{i+4}=v_i$ for all $i\ge 0$.}

    \textit{Proof of Claim 3:}
    Since $n\ge 3$, the set $\{u_1,\dots,u_{n-2}\}$ is non-empty.
    Since $u_1$ is adjacent to $v_i,v_{i+1},\dots,v_{i+4}$ (by Claim 1) and all edges have magnitude $n$ (by Lemma~\ref{lem : quadric irred}) and hence the same curve class (as $b_2(X)=1$), applying Corollary~\ref{cor : curve class triangle} to any triangle $v_ju_1v_{j+1}$ gives
    \[
        \alpha(u_1v_{j+1}) = \alpha(u_1v_j)+\alpha(v_{j}v_{j+1})
    \]
    and hence
    \[
        \alpha(u_1v_{i+4}) = \alpha(u_1v_i) + \sum_{j=0}^3 \alpha(v_{i+j}v_{i+j+1}).
    \]
    Since the sum on the right is zero by Claim 2, we obtain $\alpha(u_1v_{i+4})=\alpha(u_1v_i)$ and therefore $v_i=v_{i+4}$ by 2-independence of $\Gfat$, proving the Claim.

    \textit{Conclusion:}
    By Claim 2 and 3, we have
    \[
    \alpha(vx)=\alpha(wy)=\alpha(v_1v_2)=\alpha(v_4v_3)=\alpha(vv_3).
    \]
    By 2-independence of $\Gfat$, we therefore obtain $v_3=x$, so $xy=v_3v_2\in E(\G)^\pm$ and, by Claim 2, we have $\alpha(xy)=\alpha(v_3v_2)=\alpha(v_0v_1)=\alpha(vw)$ as desired.
    The \emph{moreover} part of the statement follows from Claim 1 since the common neighbours of $v,w$ are precisely $u_1,\dots,u_{n-2}$.
\end{proof}

\begin{lemma}\label{lem : quadrics T2 unique opposites}
    In case (T2), for each $v\in V(\G)$ there exists a unique vertex at distance two from $v$.
\end{lemma}
\begin{proof}
    Since each edge is irreducible of type (T2) by Lemma~\ref{lem : quadric T2 all irred}, we may use the local model of (T2) from Lemma~\ref{lem : quadrics T1 T2} to order the neighbours of $v$ as $x_1,\dots,x_{n/2},w_1,\dots,w_{n/2}$ such that for each $i$,
    \begin{itemize}
        \item $x_i$ and $w_i$ are not adjacent, and
        \item there exists $y_i\in V(\G)\setminus\{v\}$ adjacent to $x_i$ and $w_i$ but not to $v$
        (by Lemma~\ref{lem : quadrics T2 rectangle}).
    \end{itemize}
    In particular, $2\mid n$.

    \textit{Claim 1: If $y\in V(\G)$ has distance two from $v$ then $y\in\{y_1,\dots,y_{n/2}\}$.}
    
    \textit{Proof of Claim 1:}    
    By assumption, there exists some $y'\in V(\G)$ adjacent to $v$ and $y$.
    Thus, $y'\in \{x_1,\dots,x_{n/2},w_1,\dots,w_{n/2}\}$.
    However, by definition of type (T2) in Lemma~\ref{lem : quadrics T1 T2}, each of those has a unique neighbour which is distinct from and not adjacent to $v$.
    Those unique neighbours are precisely, $y_1,\dots,y_{n/2}$, respectively, proving the Claim.

    \textit{Claim 2: We have $y_1=\cdots=y_{n/2}$.}

    \textit{Proof of Claim 2:}
    Assume $n/2\ge 2$.
    For any $i\in\{1,\dots,n/2\}$, 
    Lemma~\ref{lem : quadrics T2 rectangle} applied to $e=\{v,w_i\}$ shows that all neighbours of $v$ other than $x_i,w_i$ are also neighbours of $x_i,w_i$, and $y_i$.
    Thus,
    \begin{itemize}
        \item the neighbours of $x_i$ are $v,y_i,x_1,\dots,\widehat{x_i},\dots,x_{n/2},w_1,\dots,\widehat{w_i},\dots,w_{n/2}$, and
        \item $y_i$ is adjacent to all of $x_1,\dots,x_{n/2},w_1,\dots,w_{n/2}$.
    \end{itemize}
    In particular, $y_1$ is adjacent to $x_2$, so by the list of neighbours of $x_2$ and using that $y_1$ has distance two from $v$, we must have $y_1=y_2$.
    By symmetry, $y_1=\cdots=y_{n/2}$, proving the Claim.

    Combining Claim 1 and Claim 2, the unique vertex at distance two from $v$ is $y_1$.
\end{proof}

\begin{proposition}\label{prop : quadrics shape of G}
    Let $\Gfat$ be the GKM graph of a Hamiltonian GKM space $(X,\omega,T)$ of dimension $2n$ with $\psind(X,\omega)=n$ and $n\ge 3$.
    Then $\Gfat$ satisfies the condition of Proposition~\ref{prop : quadric criterion}.
\end{proposition}
\begin{proof}
    Recall from Lemma~\ref{lem : quadrics all T1 or all T2} that we are either in case (T1) or (T2).

    For case (T1), let $v\in V(\G)$.
    By Lemma~\ref{lem : quadric irred} there exists an irreducible edge $\{v,w\}\in E(\G)$.
    By Lemma~\ref{lem : quadrics T1 T2}, all neighbours of $w$ other than $v$ are also neighbours of $v$.
    Since this holds for all $v$, $\G$ is a complete graph and $|V(\G)|=n+1$, since $\G$ is $n$-valent.

    Moreover, for each $v\in V(\G)$ we can define $\overline{v}$ to be the unique neighbour of $v$ such that $\m(v\overline{v})=2n$.
    (That is, $\overline{v}=y$ on the left of Figure~\ref{fig: quadrics T1 T2}.)
    The assignment $v\mapsto\overline{v}$ is then automatically a fixed-point free involution, since each vertex has a unique edge of magnitude $2n$ in case (T1), and $\G$ is simple, so has no self-loops at vertices.
    In particular, $|V(\G)|=n+1$ is even, so $2\nmid n$.

    Let $e_0\in E(\G)$ be irreducible and let $C:=[C_{e_0}]\in H_2(X;\ZZ)$.
    Then since $b_2(X)=1$ and each edge $e=\{v,w\}\in E(\G)$ has magnitude $n$ or $2n$, we have
    \[
    [C_e] = \begin{cases}
        C & \text{if }\m(e)=n,\\
        2C & \text{if }\m(e)=2n.
    \end{cases}
    \]
    The second and third statement on $\alpha$ in the condition of Proposition~\ref{prop : quadric criterion} now follow immediately by applying Proposition~\ref{prop : curve class relations} to any triangle in $\G$, using the above description of $[C_e]$.
    The statement that $\alpha(vw)=\alpha(\overline w\,\overline v)$ follows since either $w=\overline{v}$ (which implies $\overline w\,\overline v=vw$), or $w\neq \overline{v}$ and so $\m(vw)=\m(w\overline{v})=n$.
    In the latter case, Lemma~\ref{lem : quadric vx || wy} applied to the irreducible edge $w\overline{v}$ implies the desired statement. (Note that in this application of the lemma, we have $x=v$ and $y=\overline{w}$.)
    
    So far, we have shown that case (T1) implies $2\nmid n$ and the proposition holds in this case.

    Let us now deal with case (T2).
    By Lemma~\ref{lem : quadrics T2 unique opposites}, each $v\in V(\G)$ has a unique vertex at distance two from $v$, which we call $\overline{v}$.
    By \cite[Theorem~3.5]{GaarKrenn2023}, since $\G$ is regular, $\G$ is a complete graph on $n+2$ vertices with only the edges $\{v,\overline{v}\}$ removed, and $2\mid n$.

    By Lemma~\ref{lem : quadric irred} and \ref{lem : quadric T2 all irred}, all edges are irreducible of magnitude $n$, and hence of the same curve class (as $b_2(X)=1$).
    Corollary~\ref{cor : curve class triangle} therefore gives the desired
    \[
    \alpha(uv)+\alpha(vw)=\alpha(uw)
    \]
    for any distinct $u,v,w\in V(\G)$ forming a triangle.
    Finally, each edge $\{v,w\}\in E(\G)$ is of type (T2), so in the terminology of Lemma~\ref{lem : quadrics T2 rectangle} we have $\overline{v}=y$ and $\overline{w}=x$, and hence $\alpha(vw)=\alpha(\overline{w}\,\overline{v})$ by the same lemma.
    
    In total, we have shown that (T1) implies $2\nmid n$ and (T2) implies $2\mid n$, and the proposition holds in both cases.
\end{proof}

Finally, Theorem~\ref{thm : quadrics} follows immediately from Proposition~\ref{prop : quadrics shape of G} and Proposition~\ref{prop : quadric criterion} for $n\ge 3$, and Lemma~\ref{lem : quadric n ge 3, b2 is 1} for $n\in \{1,2\}$.

\section{Reflexive Delzant polytopes}
\label{sec : polytopes}

In this section, we translate Corollary~\ref{cor : connectedness} and Theorem~\ref{thm: local Mukai} to statements about \emph{reflexive Delzant polytopes}.
We recall the following background from the exposition given in \cite[\S2--3]{MR3695881}.
By a \emph{polytope} we mean the convex hull of finitely many points in $\RR^n$.
By an \emph{edge} of a polytope we mean a 1-dimensional face of it. By a \emph{lattice polytope} we mean a polytope in $\RR^n$ whose vertices (i.e., 0-dimensional faces) all lie in $\ZZ^n$.

\begin{definition}\label{def : Delzant polytope}
    A polytope $\Delta\subset\RR^n$ is \emph{Delzant} if the following conditions hold.
\begin{enumerate}
    \item $\Delta$ is \emph{simple}: each vertex is incident to precisely $n$ edges.
    \item $\Delta$ is \emph{rational}: for any vertices $v,w$ of $\Delta$ connected by an edge, the vector $w-v$ is proportional over $\RR\setminus\{0\}$ to an element of $\ZZ^n$.
    
    \emph{Notation:} let $(w-v)_\text{prim}$ denote the \emph{primitive integer vector} associated to $w-v$.
    That is, $(w-v)_\text{prim}$ is the unique element of $\ZZ^n$ with coprime components such that $(w-v)$ is a positive real multiple of $(w-v)_\text{prim}$.
    \item $\Delta$ is \emph{smooth}: for each vertex $v$ of $\Delta$, let $w_1,\dots,w_n$ be the vertices of $\Delta$ connected to $v$ by an edge.
    Then
    \[
    \left((w_1-v)_\text{prim},\dots,(w_n-v)_\text{prim}\right)
    \]
    forms a basis of the lattice $\ZZ^n$.
\end{enumerate}
\end{definition}

By the convexity theorem \cite{Ati82, MR664117}, the moment map image of any closed symplectic toric manifold is a Delzant polytope.
By \cite{MR984900}, this assignment gives rise to a bijection between Delzant polytopes $\Delta\subset\RR^n$ (up to translation in $\mathbb{R}^n$) and closed symplectic toric manifolds $(X,\omega,T)$ of dimension $2n$ (up to equivariant symplectomorphism).

Under this bijection, the GKM graph of $(X,\omega,T)$ is $\Gfat=(\G,\alpha)$ where $\G$ is the 1-skeleton of $\Delta$ and the axial function $\alpha$ assigns to each edge $vw$ of $\Delta$ the vector $(w-v)_\text{prim}$ (cf.\ \cite[Lemma~3.12]{MR3695881}).
By the same lemma, if $w-v = \lambda_{vw}(w-v)_{\text{prim}}$ then 
\begin{equation}\label{eqn : affine length is volume}
\int_{C_{\{v,w\}}}\omega = \lambda_{vw}.
\end{equation}

\begin{definition}
    Let $\Delta\subset\RR^n$ be a lattice polytope and let $e$ be an edge of $\Delta$.
    Let $v,w$ be vertices of $\Delta$ connected by $e$. The \emph{affine length of $e$}, denoted by $L(e)$, is the greatest common divisor of the components of $w-v\in\ZZ^n$.
\end{definition}

Specialising equation~\eqref{eqn : affine length is volume} to the case when $\Delta$ is a lattice polytope, each $e\in E(\G)$ satisfies
\begin{equation}\label{eqn : affine length is volume 2}
\int_{C_e}\omega = L(e).
\end{equation}

Next, we specialise to \emph{reflexive Delzant polytopes} since this allows us to access our results on magnitudes.
Originally, reflexive polytopes were introduced by Batyrev in the context of mirror symmetry \cite{MR1269718}.

\begin{definition}\label{def : reflexive polytope}
    A polytope $\Delta\subset \RR^n$ is \emph{reflexive} if the following conditions hold.
    \begin{enumerate}
        \item $\Delta$ contains the origin $0\in\RR^n$ in its interior.
        \item $\Delta$ is a lattice polytope.
        \item The \emph{dual polytope}
        \[
        \Delta^* := \{y\in\RR^n : \langle x, y\rangle\ge -1\,\,\,\forall\,x\in\Delta\}
        \]
        is also a lattice polytope.
    \end{enumerate}
\end{definition}

Given a Delzant polytope $\Delta\subset\RR^n$, its associated closed symplectic toric manifold $(X,\omega,T)$ satisfies $c_1(TX)=[\omega]$ if and only if $\Delta$ is reflexive up to translation in $\RR^n$ (cf.\ \cite[Proposition~3.10]{MR3695881}).
We are now ready to state and prove Corollary~\ref{cor : intro polytopes}.

\begin{corollary}\label{cor : polytopes}
    Let $\Delta\subset\RR^n$ be a reflexive Delzant polytope. Then
    \begin{enumerate}
        \item The $1$-skeleton of $\Delta$ remains connected after removing all edges of affine length greater than $n+1$.
        \item Let $v$ be any vertex of $\Delta$ and let $L_{\min}(v)$ be the minimum affine length of the edges of $\Delta$ incident to $v$.
        Then there exist at least $L_{\min}(v) - 1$ edges of affine length $L_{\min}(v)$ at $v$.
    \end{enumerate}
\end{corollary}
\begin{proof}
    Since $\Delta$ is Delzant, let $(X,\omega,T)$ be the closed toric symplectic manifold associated to $\Delta$, and let $\Gfat=(\G,\alpha)$ be its GKM graph, which is described above.
    Since $\Delta$ is reflexive, we have $c_1(TX)=[\omega]$ and equation~\eqref{eqn : affine length is volume 2} becomes
    \[
    \m(e) = \int_{C_e} c_1(TX) = \int_{C_e}\omega = L(e)
    \]
    for each edge $e$ of $\Delta$.
    Since $\G$ is precisely the 1-skeleton of $\Delta$, statements (1) and (2) follow immediately from Corollary~\ref{cor : connectedness} and Theorem~\ref{thm: local Mukai}, respectively, both applied to the Hamiltonian GKM space $(X,\omega,T)$.
    To apply Theorem~\ref{thm: local Mukai}, we use here that an edge $e$ of $\G$ that minimises $\m(e)$ among edges at $v$ automatically has $[C_e]$ minimal among the classes of edges at $v$ because $c_1(TX)=[\omega]$ is positive on all edge curve classes (cf.\ Definition~\ref{def : partial order on H2}).
\end{proof}

\bibliographystyle{amsalpha}
\bibliography{refs}

@misc{Holmes_Muratore_2025a, title={\textsf{GKMtools.jl}, {A} {J}ulia package for {GKM} spaces and their equivariant {G}romov--{W}itten theory}, url={https://mgemath.github.io/GKMtools.jl/}, author={Holmes, Daniel and Muratore, Giosuè}, year={2025}, language={en},
    note={Documentation and instructions available at \href{https://mgemath.github.io/GKMtools.jl/}{https://mgemath.github.io/GKMtools.jl/}},
}

@misc{Holmes_Muratore_2025b, type={Preprint}, title={{C}omputations in equivariant {G}romov--{W}itten theory of {GKM} spaces}, author={Holmes, Daniel and Muratore, Giosuè}, year={2025}, month=sep, language={en}, note="Preprint (submitted), \href{https://arxiv.org/abs/2509.07562}{arXiv:2509.07562 [math.AG]}" }

@article{Atiyah_Bott_1984,
    AUTHOR = {Atiyah, M. F. and Bott, R.},
     TITLE = {The moment map and equivariant cohomology},
   JOURNAL = {Topology},
  FJOURNAL = {Topology. An International Journal of Mathematics},
    VOLUME = {23},
      YEAR = {1984},
    NUMBER = {1},
     PAGES = {1--28},
      ISSN = {0040-9383},
   MRCLASS = {58F05 (55N91 57R20 57R22)},
  MRNUMBER = {721448},
MRREVIEWER = {G.\ J.\ Heckman},
       DOI = {10.1016/0040-9383(84)90021-1},
       URL = {https://doi.org/10.1016/0040-9383(84)90021-1},
}

@book {Cox_Katz_1999,
    AUTHOR = {Cox, David A. and Katz, Sheldon},
     TITLE = {Mirror symmetry and algebraic geometry},
    SERIES = {Mathematical Surveys and Monographs},
    VOLUME = {68},
 PUBLISHER = {American Mathematical Society, Providence, RI},
      YEAR = {1999},
     PAGES = {xxii+469},
      ISBN = {0-8218-1059-6},
   MRCLASS = {14J32 (14-02 14M25 14N10 14N35 32G81 32J81 32Q25)},
  MRNUMBER = {1677117},
MRREVIEWER = {Andreas\ Gathmann},
       DOI = {10.1090/surv/068},
       URL = {https://doi.org/10.1090/surv/068},
}

@book {MR766741,
    AUTHOR = {Kirwan, Frances Clare},
     TITLE = {Cohomology of quotients in symplectic and algebraic geometry},
    SERIES = {Mathematical Notes},
    VOLUME = {31},
 PUBLISHER = {Princeton University Press, Princeton, NJ},
      YEAR = {1984},
     PAGES = {i+211},
      ISBN = {0-691-08370-3},
   MRCLASS = {58F05 (14D25 14L30)},
  MRNUMBER = {766741},
MRREVIEWER = {Martin\ A.\ Guest},
       DOI = {10.2307/j.ctv10vm2m8},
       URL = {https://doi.org/10.2307/j.ctv10vm2m8},
}

@article{Guillemin_Zara_2002,
    AUTHOR = {Guillemin, Victor and Zara, Catalin},
     TITLE = {Equivariant de {R}ham theory and graphs},
   JOURNAL = {Asian J. Math.},
  FJOURNAL = {Asian Journal of Mathematics},
 PUBLISHER = {Int. Press, Somerville, MA},
    VOLUME = {3},
      YEAR = {1999},
    NUMBER = {1},
     PAGES = {49--76},
      ISSN = {1093-6106,1945-0036},
   MRCLASS = {58J20 (53D20 55N91)},
  MRNUMBER = {1701922},
MRREVIEWER = {Johannes\ Huebschmann},
       DOI = {10.4310/AJM.1999.v3.n1.a3},
       URL = {https://doi.org/10.4310/AJM.1999.v3.n1.a3},
}

@book{Hsiang_1975,
address={Berlin, Heidelberg},
title={Cohomology Theory of Topological Transformation Groups},
ISBN={978-3-642-66054-2},
url={http://link.springer.com/10.1007/978-3-642-66052-8},
DOI={10.1007/978-3-642-66052-8},
publisher={Springer},
author={Hsiang, Wu Yi},
year={1975}
}

@article {MR87176,
    AUTHOR = {Grothendieck, Alexander},
     TITLE = {Sur la classification des fibr\'es holomorphes sur la sph\`ere
              de {R}iemann},
   JOURNAL = {Amer. J. Math.},
  FJOURNAL = {American Journal of Mathematics},
    VOLUME = {79},
      YEAR = {1957},
     PAGES = {121--138},
      ISSN = {0002-9327,1080-6377},
   MRCLASS = {53.3X},
  MRNUMBER = {87176},
MRREVIEWER = {C.\ Chevalley},
       DOI = {10.2307/2372388},
       URL = {https://doi.org/10.2307/2372388},
}

@article {Li_2017,
    AUTHOR = {Li, Ping},
     TITLE = {The {H}irzebruch {$\chi_y$}-genus and {P}oincar{\'e} polynomial revisited},
   JOURNAL = {Commun. Contemp. Math.},
  FJOURNAL = {Communications in Contemporary Mathematics},
    VOLUME = {19},
      YEAR = {2017},
    NUMBER = {5},
     PAGES = {1650048, 19},
      ISSN = {0219-1997,1793-6683},
   MRCLASS = {32Q55 (32G15 37J10 53D20)},
  MRNUMBER = {3670789},
MRREVIEWER = {V.\ V.\ Chueshev},
       DOI = {10.1142/S0219199716500486},
       URL = {https://doi.org/10.1142/S0219199716500486},
}

@article {Godinho_Sabatini_2014,
    AUTHOR = {Godinho, Leonor and Sabatini, Silvia},
     TITLE = {New tools for classifying {H}amiltonian circle actions with
              isolated fixed points},
   JOURNAL = {Found. Comput. Math.},
  FJOURNAL = {Foundations of Computational Mathematics. The Journal of the
              Society for the Foundations of Computational Mathematics},
    VOLUME = {14},
      YEAR = {2014},
    NUMBER = {4},
     PAGES = {791--860},
      ISSN = {1615-3375,1615-3383},
   MRCLASS = {53D20 (19J35 37B05)},
  MRNUMBER = {3230015},
MRREVIEWER = {Manuel\ de Le\'on},
       DOI = {10.1007/s10208-014-9204-1},
       URL = {https://doi.org/10.1007/s10208-014-9204-1},
}

@article {Ati82,
    AUTHOR = {Atiyah, M. F.},
     TITLE = {Convexity and commuting {H}amiltonians},
   JOURNAL = {Bull. London Math. Soc.},
  FJOURNAL = {The Bulletin of the London Mathematical Society},
    VOLUME = {14},
      YEAR = {1982},
    NUMBER = {1},
     PAGES = {1--15},
      ISSN = {0024-6093,1469-2120},
   MRCLASS = {53C15 (22E30 53C55 58F05)},
  MRNUMBER = {642416},
MRREVIEWER = {A.\ Morimoto},
       DOI = {10.1112/blms/14.1.1},
       URL = {https://doi.org/10.1112/blms/14.1.1},
}

@incollection {GKZ22,
    AUTHOR = {Goertsches, Oliver and Konstantis, Panagiotis and Zoller,
              Leopold},
     TITLE = {Low-dimensional {GKM} theory},
 BOOKTITLE = {Group actions and equivariant cohomology},
    SERIES = {Contemp. Math.},
    VOLUME = {808},
     PAGES = {113--146},
 PUBLISHER = {Amer. Math. Soc., [Providence], RI},
      YEAR = {2024},
      ISBN = {978-1-4704-7180-4; [9781470477233]},
   MRCLASS = {57R91 (55N91 57S12)},
  MRNUMBER = {4840842},
       DOI = {10.1090/conm/808/16183},
       URL = {https://doi.org/10.1090/conm/808/16183},
}

@article {GKM98,
    AUTHOR = {Goresky, Mark and Kottwitz, Robert and MacPherson, Robert},
     TITLE = {Equivariant cohomology, {K}oszul duality, and the localization
              theorem},
   JOURNAL = {Invent. Math.},
  FJOURNAL = {Inventiones Mathematicae},
    VOLUME = {131},
      YEAR = {1998},
    NUMBER = {1},
     PAGES = {25--83},
      ISSN = {0020-9910,1432-1297},
   MRCLASS = {55N91 (14F25 14F32 16E99 18G10 55N33)},
  MRNUMBER = {1489894},
MRREVIEWER = {Roy\ Joshua},
       DOI = {10.1007/s002220050197},
}

@article {Graber_Pandharipande_1999,
    AUTHOR = {Graber, Tom and Pandharipande, Rahul},
     TITLE = {Localization of virtual classes},
   JOURNAL = {Invent. Math.},
  FJOURNAL = {Inventiones Mathematicae},
    VOLUME = {135},
      YEAR = {1999},
    NUMBER = {2},
     PAGES = {487--518},
      ISSN = {0020-9910,1432-1297},
   MRCLASS = {14C17 (14D20 14N10 14N35)},
  MRNUMBER = {1666787},
MRREVIEWER = {Paolo\ Aluffi},
       DOI = {10.1007/s002220050293},
       URL = {https://doi.org/10.1007/s002220050293},
}

@book{Guillemin_Guillemin_Ohsawa_Ginzburg_Karshon_2014,
    AUTHOR = {Guillemin, Victor and Ginzburg, Viktor and Karshon, Yael},
     TITLE = {Moment maps, cobordisms, and {H}amiltonian group actions},
    SERIES = {Mathematical Surveys and Monographs},
    VOLUME = {98},
      NOTE = {Appendix J by Maxim Braverman},
 PUBLISHER = {American Mathematical Society, Providence, RI},
      YEAR = {2002},
     PAGES = {viii+350},
      ISBN = {0-8218-0502-9},
   MRCLASS = {53D20 (53-02 53D50 57R90 57S15)},
  MRNUMBER = {1929136},
MRREVIEWER = {Tara\ S.\ Holm},
       DOI = {10.1090/surv/098},
       URL = {https://doi.org/10.1090/surv/098},
}

@article {GZ01,
    AUTHOR = {Guillemin, Victor and Zara, Catalin},
     TITLE = {1-skeleta, {B}etti numbers, and equivariant cohomology},
   JOURNAL = {Duke Math. J.},
  FJOURNAL = {Duke Mathematical Journal},
    VOLUME = {107},
      YEAR = {2001},
    NUMBER = {2},
     PAGES = {283--349},
      ISSN = {0012-7094,1547-7398},
   MRCLASS = {53D20 (55N91 57S15)},
  MRNUMBER = {1823050},
MRREVIEWER = {Thomas\ Delzant},
       DOI = {10.1215/S0012-7094-01-10724-2},
       URL = {https://doi.org/10.1215/S0012-7094-01-10724-2},
}

@article {Guillemin_Holm_Zara_2006,
    AUTHOR = {Guillemin, Victor and Holm, Tara and Zara, Catalin},
     TITLE = {A {GKM} description of the equivariant cohomology ring of a
              homogeneous space},
   JOURNAL = {J. Algebraic Combin.},
  FJOURNAL = {Journal of Algebraic Combinatorics. An International Journal},
    VOLUME = {23},
      YEAR = {2006},
    NUMBER = {1},
     PAGES = {21--41},
      ISSN = {0925-9899,1572-9192},
   MRCLASS = {55N91 (22F30 57S10)},
  MRNUMBER = {2218848},
MRREVIEWER = {Tomasz\ Rybicki},
       DOI = {10.1007/s10801-006-6027-4},
       URL = {https://doi.org/10.1007/s10801-006-6027-4},
}

@article {LS17,
    AUTHOR = {Liu, Chiu-Chu Melissa and Sheshmani, Artan},
     TITLE = {Equivariant {G}romov--{W}itten invariants of algebraic {GKM}
              manifolds},
   JOURNAL = {SIGMA Symmetry Integrability Geom. Methods Appl.},
  FJOURNAL = {SIGMA. Symmetry, Integrability and Geometry. Methods and
              Applications},
    VOLUME = {13},
      YEAR = {2017},
     PAGES = {Paper No. 048, 21},
      ISSN = {1815-0659},
   MRCLASS = {14N35 (14D20 14H10)},
  MRNUMBER = {3667222},
MRREVIEWER = {Jie\ Zhou},
       DOI = {10.3842/SIGMA.2017.048},
}

@incollection {Manin_1995,
    AUTHOR = {Manin, Yu. I.},
     TITLE = {Generating functions in algebraic geometry and sums over
              trees},
 BOOKTITLE = {The moduli space of curves ({T}exel {I}sland, 1994)},
    SERIES = {Progr. Math.},
    VOLUME = {129},
     PAGES = {401--417},
 PUBLISHER = {Birkh\"auser Boston, Boston, MA},
      YEAR = {1995},
      ISBN = {0-8176-3784-2},
   MRCLASS = {14N10 (05C05 14H10)},
  MRNUMBER = {1363064},
MRREVIEWER = {Bernd\ Siebert},
       DOI = {10.1007/978-1-4612-4264-2\_14},
       URL = {https://doi.org/10.1007/978-1-4612-4264-2_14},
}

@article {Hirschi_2023,
    AUTHOR = {Hirschi, Amanda},
     TITLE = {Properties of {G}romov-{W}itten invariants defined via global
              {K}uranishi charts},
   JOURNAL = {Ann. H. Lebesgue},
  FJOURNAL = {Annales Henri Lebesgue},
    VOLUME = {8},
      YEAR = {2025},
     PAGES = {811--871},
      ISSN = {2644-9463},
   MRCLASS = {53D45 (32Q65)},
  MRNUMBER = {5005386},
       DOI = {10.5802/ahl.249},
       URL = {https://doi.org/10.5802/ahl.249},
}

@article{Charton_Kessler_2023,
    AUTHOR = {Charton, Isabelle and Kessler, Liat},
     TITLE = {Monotone symplectic six-manifolds that admit a {H}amiltonian
              {GKM} action are diffeomorphic to smooth {F}ano threefolds},
   JOURNAL = {Transform. Groups},
  FJOURNAL = {Transformation Groups},
    VOLUME = {31},
      YEAR = {2026},
    NUMBER = {1},
     PAGES = {263--326},
      ISSN = {1083-4362,1531-586X},
   MRCLASS = {14J45 (53D20 57R91 57S25)},
  MRNUMBER = {5005925},
       DOI = {10.1007/s00031-025-09914-2},
       URL = {https://doi.org/10.1007/s00031-025-09914-2},
}

@article {Goldin_Tymoczko_2023,
    AUTHOR = {Goldin, Rebecca and Tymoczko, Julianna},
     TITLE = {Which {H}essenberg varieties are {GKM}?},
   JOURNAL = {Pure Appl. Math. Q.},
  FJOURNAL = {Pure and Applied Mathematics Quarterly},
    VOLUME = {19},
      YEAR = {2023},
    NUMBER = {4},
     PAGES = {1899--1942},
      ISSN = {1558-8599,1558-8602},
   MRCLASS = {14M15},
  MRNUMBER = {4671386},
MRREVIEWER = {Praise\ Adeyemo},
       DOI = {10.4310/pamq.2023.v19.n4.a8},
       URL = {https://doi.org/10.4310/pamq.2023.v19.n4.a8},
}

@article {MR2064318_Guillemin_Holm_04,
    AUTHOR = {Guillemin, Victor and Holm, Tara S.},
     TITLE = {G{KM} theory for torus actions with nonisolated fixed points},
   JOURNAL = {Int. Math. Res. Not.},
  FJOURNAL = {International Mathematics Research Notices},
      YEAR = {2004},
    NUMBER = {40},
     PAGES = {2105--2124},
      ISSN = {1073-7928,1687-0247},
   MRCLASS = {53D20},
  MRNUMBER = {2064318},
MRREVIEWER = {Catalin\ Zara},
       DOI = {10.1155/S1073792804132285},
       URL = {https://doi.org/10.1155/S1073792804132285},
}

@article {Tolman_1998,
    AUTHOR = {Tolman, Susan},
     TITLE = {Examples of non-{K}\"ahler {H}amiltonian torus actions},
   JOURNAL = {Invent. Math.},
  FJOURNAL = {Inventiones Mathematicae},
    VOLUME = {131},
      YEAR = {1998},
    NUMBER = {2},
     PAGES = {299--310},
      ISSN = {0020-9910,1432-1297},
   MRCLASS = {53D20},
  MRNUMBER = {1608575},
MRREVIEWER = {Stephanie\ F.\ Singer},
       DOI = {10.1007/s002220050205},
       URL = {https://doi.org/10.1007/s002220050205},
}

@article {Woodward_1998,
    AUTHOR = {Woodward, Chris},
     TITLE = {Multiplicity-free {H}amiltonian actions need not be
              {K}\"ahler},
   JOURNAL = {Invent. Math.},
  FJOURNAL = {Inventiones Mathematicae},
    VOLUME = {131},
      YEAR = {1998},
    NUMBER = {2},
     PAGES = {311--319},
      ISSN = {0020-9910,1432-1297},
   MRCLASS = {53D20},
  MRNUMBER = {1608579},
MRREVIEWER = {Stephanie\ F.\ Singer},
       DOI = {10.1007/s002220050206},
       URL = {https://doi.org/10.1007/s002220050206},
}

@article {Rosu_03,
    AUTHOR = {Rosu, Ioanid},
     TITLE = {Equivariant {$K$}-theory and equivariant cohomology},
      NOTE = {With an appendix by Allen Knutson and Rosu},
   JOURNAL = {Math. Z.},
  FJOURNAL = {Mathematische Zeitschrift},
    VOLUME = {243},
      YEAR = {2003},
    NUMBER = {3},
     PAGES = {423--448},
      ISSN = {0025-5874,1432-1823},
   MRCLASS = {19L47 (14F43)},
  MRNUMBER = {1970011},
MRREVIEWER = {Gabriele\ Vezzosi},
       DOI = {10.1007/s00209-002-0447-1},
       URL = {https://doi.org/10.1007/s00209-002-0447-1},
}

@incollection {Kontsevich_torus_enumeration,
    AUTHOR = {Kontsevich, Maxim},
     TITLE = {Enumeration of rational curves via torus actions},
 BOOKTITLE = {The moduli space of curves ({T}exel {I}sland, 1994)},
    SERIES = {Progr. Math.},
    VOLUME = {129},
     PAGES = {335--368},
 PUBLISHER = {Birkh\"auser Boston, Boston, MA},
      YEAR = {1995},
      ISBN = {0-8176-3784-2},
   MRCLASS = {14N10 (14D22 14L30)},
  MRNUMBER = {1363062},
MRREVIEWER = {Anatoly\ Libgober},
       DOI = {10.1007/978-1-4612-4264-2\_12},
       URL = {https://doi.org/10.1007/978-1-4612-4264-2_12},
}

@article {MR3665359_Sabatini_2017,
    AUTHOR = {Sabatini, Silvia},
     TITLE = {On the {C}hern numbers and the {H}ilbert polynomial of an
              almost complex manifold with a circle action},
   JOURNAL = {Commun. Contemp. Math.},
  FJOURNAL = {Communications in Contemporary Mathematics},
    VOLUME = {19},
      YEAR = {2017},
    NUMBER = {4},
     PAGES = {1750043, 51},
      ISSN = {0219-1997,1793-6683},
   MRCLASS = {57R20 (37J10 57S15)},
  MRNUMBER = {3665359},
MRREVIEWER = {Zhi\ L\"u},
       DOI = {10.1142/S0219199717500432},
       URL = {https://doi.org/10.1142/S0219199717500432},
}

@article {MR1189503,
    AUTHOR = {Koll\'ar, J\'anos and Miyaoka, Yoichi and Mori, Shigefumi},
     TITLE = {Rational connectedness and boundedness of {F}ano manifolds},
   JOURNAL = {J. Differential Geom.},
  FJOURNAL = {Journal of Differential Geometry},
    VOLUME = {36},
      YEAR = {1992},
    NUMBER = {3},
     PAGES = {765--779},
      ISSN = {0022-040X,1945-743X},
   MRCLASS = {14J45},
  MRNUMBER = {1189503},
MRREVIEWER = {Yuri\ G.\ Prokhorov},
       URL = {http://projecteuclid.org/euclid.jdg/1214453188},
}

@book {MR4655919,
    AUTHOR = {Anderson, David and Fulton, William},
     TITLE = {Equivariant cohomology in algebraic geometry},
    SERIES = {Cambridge Studies in Advanced Mathematics},
    VOLUME = {210},
 PUBLISHER = {Cambridge University Press, Cambridge},
      YEAR = {2024},
     PAGES = {xv+446},
      ISBN = {978-1-00-934998-7},
   MRCLASS = {14L30 (05E14 14-02 14F43 14Mxx 20G05 55N91)},
  MRNUMBER = {4655919},
MRREVIEWER = {Michael\ Orin\ Joyce},
}

@article {MR3695881,
    AUTHOR = {Godinho, Leonor and von Heymann, Frederik and Sabatini,
              Silvia},
     TITLE = {12, 24 and beyond},
   JOURNAL = {Adv. Math.},
  FJOURNAL = {Advances in Mathematics},
    VOLUME = {319},
      YEAR = {2017},
     PAGES = {472--521},
      ISSN = {0001-8708,1090-2082},
   MRCLASS = {52B20 (14J45 14M25 53D20)},
  MRNUMBER = {3695881},
MRREVIEWER = {Mateusz\ Micha\l ek},
       DOI = {10.1016/j.aim.2017.08.023},
       URL = {https://doi.org/10.1016/j.aim.2017.08.023},
}

@article {MR776689,
    AUTHOR = {Hattori, Akio},
     TITLE = {{$S^1$}-actions on unitary manifolds and quasi-ample line
              bundles},
   JOURNAL = {J. Fac. Sci. Univ. Tokyo Sect. IA Math.},
  FJOURNAL = {Journal of the Faculty of Science. University of Tokyo.
              Section IA. Mathematics},
    VOLUME = {31},
      YEAR = {1985},
    NUMBER = {3},
     PAGES = {433--486},
      ISSN = {0040-8980},
   MRCLASS = {57S15 (57R15)},
  MRNUMBER = {776689},
MRREVIEWER = {S\"oren\ Illman},
}

@article {MR595007,
    AUTHOR = {Michelsohn, M. L.},
     TITLE = {Clifford and spinor cohomology of {K}\"ahler manifolds},
   JOURNAL = {Amer. J. Math.},
  FJOURNAL = {American Journal of Mathematics},
    VOLUME = {102},
      YEAR = {1980},
    NUMBER = {6},
     PAGES = {1083--1146},
      ISSN = {0002-9327,1080-6377},
   MRCLASS = {53C55 (32L20 58A14 58G05)},
  MRNUMBER = {595007},
MRREVIEWER = {N.\ J.\ Hitchin},
       DOI = {10.2307/2374181},
       URL = {https://doi.org/10.2307/2374181},
}

@article {MR984900,
    AUTHOR = {Delzant, Thomas},
     TITLE = {Hamiltoniens p\'eriodiques et images convexes de l'application
              moment},
   JOURNAL = {Bull. Soc. Math. France},
  FJOURNAL = {Bulletin de la Soci\'et\'e{} Math\'ematique de France},
    VOLUME = {116},
      YEAR = {1988},
    NUMBER = {3},
     PAGES = {315--339},
      ISSN = {0037-9484},
   MRCLASS = {58F05},
  MRNUMBER = {984900},
MRREVIEWER = {J.\ J.\ Duistermaat},
       URL = {http://www.numdam.org/item?id=BSMF_1988__116_3_315_0},
}

@article {MR1063222,
    AUTHOR = {Wi{\'s}niewski, Jaros{\l }aw A.},
     TITLE = {On a conjecture of {M}ukai},
   JOURNAL = {Manuscripta Math.},
  FJOURNAL = {Manuscripta Mathematica},
    VOLUME = {68},
      YEAR = {1990},
    NUMBER = {2},
     PAGES = {135--141},
      ISSN = {0025-2611,1432-1785},
   MRCLASS = {14J45 (14J05)},
  MRNUMBER = {1063222},
MRREVIEWER = {Peter\ Nielsen},
       DOI = {10.1007/BF02568756},
       URL = {https://doi.org/10.1007/BF02568756},
}

@article {MR554387,
    AUTHOR = {Mori, Shigefumi},
     TITLE = {Projective manifolds with ample tangent bundles},
   JOURNAL = {Ann. of Math. (2)},
  FJOURNAL = {Annals of Mathematics. Second Series},
    VOLUME = {110},
      YEAR = {1979},
    NUMBER = {3},
     PAGES = {593--606},
      ISSN = {0003-486X},
   MRCLASS = {14E05 (14M20)},
  MRNUMBER = {554387},
MRREVIEWER = {Luciana\ Picco Botta},
       DOI = {10.2307/1971241},
       URL = {https://doi.org/10.2307/1971241},
}

@article {MR1998396,
    AUTHOR = {Bonavero, Laurent and Casagrande, Cinzia and Debarre, Olivier
              and Druel, St\'ephane},
     TITLE = {Sur une conjecture de {M}ukai},
   JOURNAL = {Comment. Math. Helv.},
  FJOURNAL = {Commentarii Mathematici Helvetici},
    VOLUME = {78},
      YEAR = {2003},
    NUMBER = {3},
     PAGES = {601--626},
      ISSN = {0010-2571,1420-8946},
   MRCLASS = {14J45 (14E30 14M25)},
  MRNUMBER = {1998396},
MRREVIEWER = {Yuri\ G.\ Prokhorov},
       DOI = {10.1007/s00014-003-0765-x},
       URL = {https://doi.org/10.1007/s00014-003-0765-x},
}

@incollection {Mukai1988Problems,
AUTHOR = {Mukai, Shigeru},
TITLE = {Problems on characterization of the complex projective space},
BOOKTITLE = {Birational geometry of algebraic varieties, open problems
({K}atata, 1988)},
PAGES = {57--60},
PUBLISHER = {The Taniguchi Foundation, Katata},
YEAR = {1988},
}

@article {MR4671382,
    AUTHOR = {Caviedes Castro, Alexander and Pabiniak, Milena and Sabatini,
              Silvia},
     TITLE = {Generalizing the {M}ukai conjecture to the symplectic category
              and the {K}ostant game},
   JOURNAL = {Pure Appl. Math. Q.},
  FJOURNAL = {Pure and Applied Mathematics Quarterly},
    VOLUME = {19},
      YEAR = {2023},
    NUMBER = {4},
     PAGES = {1803--1837},
      ISSN = {1558-8599,1558-8602},
   MRCLASS = {53D05 (05C31)},
  MRNUMBER = {4671382},
MRREVIEWER = {Philippe\ Rukimbira},
       DOI = {10.4310/pamq.2023.v19.n4.a4},
       URL = {https://doi.org/10.4310/pamq.2023.v19.n4.a4},
}

@article {MR316745,
    AUTHOR = {Kobayashi, Shoshichi and Ochiai, Takushiro},
     TITLE = {Characterizations of complex projective spaces and
              hyperquadrics},
   JOURNAL = {J. Math. Kyoto Univ.},
  FJOURNAL = {Journal of Mathematics of Kyoto University},
    VOLUME = {13},
      YEAR = {1973},
     PAGES = {31--47},
      ISSN = {0023-608X},
   MRCLASS = {32C10},
  MRNUMBER = {316745},
MRREVIEWER = {J.\ A.\ Morrow},
       DOI = {10.1215/kjm/1250523432},
       URL = {https://doi.org/10.1215/kjm/1250523432},
}

@article {MR4088417,
    AUTHOR = {Goertsches, Oliver and Konstantis, Panagiotis and Zoller,
              Leopold},
     TITLE = {G{KM} theory and {H}amiltonian non-{K}\"ahler actions in
              dimension 6},
   JOURNAL = {Adv. Math.},
  FJOURNAL = {Advances in Mathematics},
    VOLUME = {368},
      YEAR = {2020},
     PAGES = {107141, 17},
      ISSN = {0001-8708,1090-2082},
   MRCLASS = {57S15 (55N91 57R91)},
  MRNUMBER = {4088417},
MRREVIEWER = {Andrea\ Spiro},
       DOI = {10.1016/j.aim.2020.107141},
       URL = {https://doi.org/10.1016/j.aim.2020.107141},
}

@article {MR4545925,
    AUTHOR = {Goertsches, Oliver and Konstantis, Panagiotis and Zoller,
              Leopold},
     TITLE = {G{KM} manifolds are not rigid},
   JOURNAL = {Algebr. Geom. Topol.},
  FJOURNAL = {Algebraic \& Geometric Topology},
    VOLUME = {22},
      YEAR = {2022},
    NUMBER = {7},
     PAGES = {3511--3532},
      ISSN = {1472-2747,1472-2739},
   MRCLASS = {57R91 (53D20 55N91)},
  MRNUMBER = {4545925},
       DOI = {10.2140/agt.2022.22.3511},
       URL = {https://doi.org/10.2140/agt.2022.22.3511},
}

@article {MR1269718,
    AUTHOR = {Batyrev, Victor V.},
     TITLE = {Dual polyhedra and mirror symmetry for {C}alabi-{Y}au
              hypersurfaces in toric varieties},
   JOURNAL = {J. Algebraic Geom.},
  FJOURNAL = {Journal of Algebraic Geometry},
    VOLUME = {3},
      YEAR = {1994},
    NUMBER = {3},
     PAGES = {493--535},
      ISSN = {1056-3911,1534-7486},
   MRCLASS = {14J32 (14M25)},
  MRNUMBER = {1269718},
MRREVIEWER = {I.\ Dolgachev},
}

@book {MR3674984,
    AUTHOR = {McDuff, Dusa and Salamon, Dietmar},
     TITLE = {Introduction to symplectic topology},
    SERIES = {Oxford Graduate Texts in Mathematics},
   EDITION = {Third},
 PUBLISHER = {Oxford University Press, Oxford},
      YEAR = {2017},
     PAGES = {xi+623},
      ISBN = {978-0-19-879490-5; 978-0-19-879489-9},
   MRCLASS = {53D35 (53D40 57R17 57R57 57R58)},
  MRNUMBER = {3674984},
MRREVIEWER = {Hansj\"org\ Geiges},
       DOI = {10.1093/oso/9780198794899.001.0001},
       URL = {https://doi.org/10.1093/oso/9780198794899.001.0001},
}

@article {MR664117,
    AUTHOR = {Guillemin, V. and Sternberg, S.},
     TITLE = {Convexity properties of the moment mapping},
   JOURNAL = {Invent. Math.},
  FJOURNAL = {Inventiones Mathematicae},
    VOLUME = {67},
      YEAR = {1982},
    NUMBER = {3},
     PAGES = {491--513},
      ISSN = {0020-9910,1432-1297},
   MRCLASS = {58F05 (57S25 70H05)},
  MRNUMBER = {664117},
MRREVIEWER = {I.\ Vaisman},
       DOI = {10.1007/BF01398933},
       URL = {https://doi.org/10.1007/BF01398933},
}

@article {MR2072765,
    AUTHOR = {Fulton, W. and Woodward, C.},
     TITLE = {On the quantum product of {S}chubert classes},
   JOURNAL = {J. Algebraic Geom.},
  FJOURNAL = {Journal of Algebraic Geometry},
    VOLUME = {13},
      YEAR = {2004},
    NUMBER = {4},
     PAGES = {641--661},
      ISSN = {1056-3911,1534-7486},
   MRCLASS = {14N35 (14M15 14N15)},
  MRNUMBER = {2072765},
MRREVIEWER = {Harry\ Tamvakis},
       DOI = {10.1090/S1056-3911-04-00365-0},
       URL = {https://doi.org/10.1090/S1056-3911-04-00365-0},
}

@article {MR2228683,
    AUTHOR = {Casagrande, Cinzia},
     TITLE = {The number of vertices of a {F}ano polytope},
   JOURNAL = {Ann. Inst. Fourier (Grenoble)},
  FJOURNAL = {Universit\'e{} de Grenoble. Annales de l'Institut Fourier},
    VOLUME = {56},
      YEAR = {2006},
    NUMBER = {1},
     PAGES = {121--130},
      ISSN = {0373-0956,1777-5310},
   MRCLASS = {14M25 (14J45 52B20)},
  MRNUMBER = {2228683},
MRREVIEWER = {Boris\ J.\ Kazarnovski\u i},
       DOI = {10.5802/aif.2175},
       URL = {https://doi.org/10.5802/aif.2175},
}

@incollection {MR2761927,
    AUTHOR = {Chen, Jiun-Cheng},
     TITLE = {On {F}ano varieties with large pseudo-index},
 BOOKTITLE = {Algebraic geometry in {E}ast {A}sia---{S}eoul 2008},
    SERIES = {Adv. Stud. Pure Math.},
    VOLUME = {60},
     PAGES = {195--213},
 PUBLISHER = {Math. Soc. Japan, Tokyo},
      YEAR = {2010},
      ISBN = {978-4-931469-63-1},
   MRCLASS = {14J45 (14E30)},
  MRNUMBER = {2761927},
MRREVIEWER = {Anne-Sophie\ Kaloghiros},
       DOI = {10.2969/aspm/06010195},
       URL = {https://doi.org/10.2969/aspm/06010195},
}

@article {MR2660578,
    AUTHOR = {Novelli, Carla and Occhetta, Gianluca},
     TITLE = {Rational curves and bounds on the {P}icard number of {F}ano
              manifolds},
   JOURNAL = {Geom. Dedicata},
  FJOURNAL = {Geometriae Dedicata},
    VOLUME = {147},
      YEAR = {2010},
     PAGES = {207--217},
      ISSN = {0046-5755,1572-9168},
   MRCLASS = {14J45 (14E30)},
  MRNUMBER = {2660578},
MRREVIEWER = {Jaros\l aw\ A.\ Wi\'sniewski},
       DOI = {10.1007/s10711-009-9452-4},
       URL = {https://doi.org/10.1007/s10711-009-9452-4},
}

@article {MR2725271,
    AUTHOR = {Pasquier, Boris},
     TITLE = {The pseudo-index of horospherical {F}ano varieties},
   JOURNAL = {Internat. J. Math.},
  FJOURNAL = {International Journal of Mathematics},
    VOLUME = {21},
      YEAR = {2010},
    NUMBER = {9},
     PAGES = {1147--1156},
      ISSN = {0129-167X,1793-6519},
   MRCLASS = {14J45 (14L30 52B20)},
  MRNUMBER = {2725271},
MRREVIEWER = {Gianluca\ Occhetta},
       DOI = {10.1142/S0129167X10006422},
       URL = {https://doi.org/10.1142/S0129167X10006422},
}

@article {MR3592522,
    AUTHOR = {Gagliardi, Giuliano and Hofscheier, Johannes},
     TITLE = {The generalized {M}ukai conjecture for symmetric varieties},
   JOURNAL = {Trans. Amer. Math. Soc.},
  FJOURNAL = {Transactions of the American Mathematical Society},
    VOLUME = {369},
      YEAR = {2017},
    NUMBER = {4},
     PAGES = {2615--2649},
      ISSN = {0002-9947,1088-6850},
   MRCLASS = {14M27 (14J45 14L30 52B20)},
  MRNUMBER = {3592522},
MRREVIEWER = {Pierre-Emmanuel\ Chaput},
       DOI = {10.1090/tran/6738},
       URL = {https://doi.org/10.1090/tran/6738},
}

@article {MR4790927,
    AUTHOR = {Reineke, Markus},
     TITLE = {The {M}ukai conjecture for {F}ano quiver moduli},
   JOURNAL = {Algebr. Represent. Theory},
  FJOURNAL = {Algebras and Representation Theory},
    VOLUME = {27},
      YEAR = {2024},
    NUMBER = {4},
     PAGES = {1641--1644},
      ISSN = {1386-923X,1572-9079},
   MRCLASS = {16G20 (14D20 14J45)},
  MRNUMBER = {4790927},
MRREVIEWER = {Izzet\ Coskun},
       DOI = {10.1007/s10468-024-10268-8},
       URL = {https://doi.org/10.1007/s10468-024-10268-8},
}

@article {MR5048523,
    AUTHOR = {Bravi, P. and Pezzini, G.},
     TITLE = {The generalized {M}ukai conjecture for spherical varieties
              with a reductive general isotropy group},
   JOURNAL = {Manuscripta Math.},
  FJOURNAL = {Manuscripta Mathematica},
    VOLUME = {177},
      YEAR = {2026},
    NUMBER = {2},
     PAGES = {Paper No. 28, 36},
      ISSN = {0025-2611,1432-1785},
   MRCLASS = {14M27 (14J45)},
  MRNUMBER = {5048523},
       DOI = {10.1007/s00229-026-01703-9},
       URL = {https://doi.org/10.1007/s00229-026-01703-9},
}

@misc{Gagliardi_Hofscheier_Pearson_2025,
  type={Preprint},
  title={{T}he generalised {M}ukai conjecture for spherical varieties},
  author={Gagliardi, Giuliano and Hofscheier, Johannes and Pearson, Heath},
  year={2025},
  month=feb,
  language={en},
  note="Preprint, \href{https://arxiv.org/abs/2502.21155}{arXiv:2502.21155 [math.AG]}"
}

@article {MR2415341,
    AUTHOR = {Pasquier, Boris},
     TITLE = {Vari\'et\'es horosph\'eriques de {F}ano},
   JOURNAL = {Bull. Soc. Math. France},
  FJOURNAL = {Bulletin de la Soci\'et\'e{} Math\'ematique de France},
    VOLUME = {136},
      YEAR = {2008},
    NUMBER = {2},
     PAGES = {195--225},
      ISSN = {0037-9484,2102-622X},
   MRCLASS = {14L30 (14J45)},
  MRNUMBER = {2415341},
MRREVIEWER = {Yuri\ Tschinkel},
       DOI = {10.24033/bsmf.2554},
       URL = {https://doi.org/10.24033/bsmf.2554},
}

@incollection {MR1929792,
    AUTHOR = {Cho, Koji and Miyaoka, Yoichi and Shepherd-Barron, N. I.},
     TITLE = {Characterizations of projective space and applications to
              complex symplectic manifolds},
 BOOKTITLE = {Higher dimensional birational geometry ({K}yoto, 1997)},
    SERIES = {Adv. Stud. Pure Math.},
    VOLUME = {35},
     PAGES = {1--88},
 PUBLISHER = {Math. Soc. Japan, Tokyo},
      YEAR = {2002},
      ISBN = {4-931469-19-1},
   MRCLASS = {14M20 (14E05)},
  MRNUMBER = {1929792},
MRREVIEWER = {S\'andor\ J.\ Kov\'acs},
       DOI = {10.2969/aspm/03510001},
       URL = {https://doi.org/10.2969/aspm/03510001},
}

@article {MR4029085,
    AUTHOR = {Charton, Isabelle},
     TITLE = {Hamiltonian {$S^1$}-spaces with large equivariant
              pseudo-index},
   JOURNAL = {J. Geom. Phys.},
  FJOURNAL = {Journal of Geometry and Physics},
    VOLUME = {147},
      YEAR = {2020},
     PAGES = {103521, 10},
      ISSN = {0393-0440,1879-1662},
   MRCLASS = {53D35 (32Q99)},
  MRNUMBER = {4029085},
       DOI = {10.1016/j.geomphys.2019.103521},
       URL = {https://doi.org/10.1016/j.geomphys.2019.103521},
}

@article{Klyachko1990,
  author  = {Klyachko, A. A.},
  title   = {Equivariant bundles on toral varieties},
  journal = {Math. USSR-Izv.},
  volume  = {35},
  year    = {1990},
  number  = {2},
  pages   = {337--375},
  doi     = {10.1070/IM1990v035n02ABEH000707},
  mrnumber = {1024452}
}

@misc{KSZ_moment_maps_in_a_box, 
    type={Preprint},
    title={{F}initeness and boundedness of positive monotone {H}amiltonian $\text{GKM}_3$ spaces},
    author={Panagiotis Konstantis and Silvia Sabatini and Leopold Zoller},
    year={2026},
    language={en},
    note="Preprint (submitted), \href{https://arxiv.org/abs/2606.28768}{arXiv:2606.28768 [math.SG]}"
}

@article{BialynickiBirula1973,
  author  = {Bia{\l}ynicki-Birula, Andrzej},
  title   = {Some theorems on actions of algebraic groups},
  journal = {Annals of Mathematics},
  series  = {2},
  volume  = {98},
  number  = {3},
  year    = {1973},
  pages   = {480--497},
  doi     = {10.2307/1970915},
  mrnumber = {0366940}
}

@article{GaarKrenn2023,
  author  = {Elisabeth Gaar and Daniel Krenn},
  title   = {A characterization of graphs with regular distance-2 graphs},
  journal = {Discrete Applied Mathematics},
  volume  = {324},
  pages   = {181--218},
  year    = {2023},
  doi     = {10.1016/j.dam.2022.09.020}
}

@article{Kuroki2025,
  author  = {Kuroki, Shintar{\^o}},
  title   = {Equivariant cohomology of even-dimensional complex quadrics
             from a combinatorial point of view},
  journal = {Osaka Journal of Mathematics},
  volume  = {62},
  number  = {4},
  pages   = {539--563},
  year    = {2025}
}

@article{KurokiPaul2026,
  author  = {Kuroki, Shintar{\^o} and Paul, Bidhan},
  title   = {Equivariant cohomology of odd-dimensional complex quadrics
             from a combinatorial point of view},
  journal = {International Journal of Mathematics},
  volume  = {37},
  number  = {3},
  pages   = {2650021},
  year    = {2026},
  doi     = {10.1142/S0129167X26500217},
  note    = {corrected version available at \href{https://arxiv.org/abs/2407.17921v3}{arXiv:2407.17921v3 [math.AT]}}
}

@misc{GKZ2026,
      title={On orientability, {P}oincar\'e duality, and connectivity of {GKM} graphs}, 
      author={Oliver Goertsches and Panagiotis Konstantis and Leopold Zoller},
      year={2026},
      eprint={2608.26999},
      archivePrefix={arXiv},
      primaryClass={math.AT},
      url={https://arxiv.org/abs/2608.26999}, 
      note="Preprint, \href{https://arxiv.org/abs/2608.26999}{arXiv:2608.26999 [math.AT]}"
}

@misc{APC_upcoming, type={Preprint}, title={{C}lassification of almost positive {H}amiltonian {GKM} spaces}, author={Holmes, Daniel}, year={2027}, language={en}, note="In preparation" }

@misc{Pearson_2026,
  title={The {Mukai} conjecture via {Cox} rings for smooth canonical ambient embeddings},
  url={https://arxiv.org/abs/2604.25023},
  DOI={10.48550/arXiv.2604.25023},
  number={arXiv:2604.25023},
  publisher={arXiv},
  author={Pearson, Heath},
  year={2026},
  month=apr,
  note={Preprint, \href{https://arxiv.org/abs/2604.25023}{arXiv:2604.25023 [math.AG]}},
}

@article {MR3592468,
    AUTHOR = {Dedieu, Thomas and H\"oring, Andreas},
     TITLE = {Numerical characterisation of quadrics},
   JOURNAL = {Algebr. Geom.},
  FJOURNAL = {Algebraic Geometry},
    VOLUME = {4},
      YEAR = {2017},
    NUMBER = {1},
     PAGES = {120--135},
      ISSN = {2313-1691,2214-2584},
   MRCLASS = {14J45 (14E30 14M22)},
  MRNUMBER = {3592468},
MRREVIEWER = {James\ McKernan},
       DOI = {10.14231/AG-2017-006},
       URL = {https://doi.org/10.14231/AG-2017-006},
}

@article {MR2997473,
    AUTHOR = {Novelli, Carla},
     TITLE = {On {F}ano manifolds with an unsplit dominating family of
              rational curves},
   JOURNAL = {Kodai Math. J.},
  FJOURNAL = {Kodai Mathematical Journal},
    VOLUME = {35},
      YEAR = {2012},
    NUMBER = {3},
     PAGES = {425--438},
      ISSN = {0386-5991,1881-5472},
   MRCLASS = {14J45 (14E30)},
  MRNUMBER = {2997473},
MRREVIEWER = {James\ McKernan},
       DOI = {10.2996/kmj/1352985447},
       URL = {https://doi.org/10.2996/kmj/1352985447},
}

@misc{OAI:The-Generalized-Mukai-Conjecture-September-24-2026,
  author = {{OpenAI}},
  title = {{The generalized Mukai conjecture}},
  howpublished = {OpenAI Math Release preprint
                  \href{https://github.com/openai/math/blob/main/preprints/The-Generalized-Mukai-Conjecture-September-24-2026/article.pdf}{OAI:The-Generalized-Mukai-Conjecture-September-24-2026}},
  year = {2026}
}

@incollection {MR1922103,
    AUTHOR = {Kebekus, Stefan},
     TITLE = {Characterizing the projective space after {C}ho, {M}iyaoka and
              {S}hepherd-{B}arron},
 BOOKTITLE = {Complex geometry ({G}\"ottingen, 2000)},
     PAGES = {147--155},
 PUBLISHER = {Springer, Berlin},
      YEAR = {2002},
      ISBN = {3-540-43259-0},
   MRCLASS = {14M20 (14J45)},
  MRNUMBER = {1922103},
MRREVIEWER = {S\'andor\ J.\ Kov\'acs},
}

@article {MR3098425,
    AUTHOR = {Hwang, Jun-Muk},
     TITLE = {Varieties of minimal rational tangents of codimension 1},
   JOURNAL = {Ann. Sci. \'Ec. Norm. Sup\'er. (4)},
  FJOURNAL = {Annales Scientifiques de l'\'Ecole Normale Sup\'erieure.
              Quatri\`eme S\'erie},
    VOLUME = {46},
      YEAR = {2013},
    NUMBER = {4},
     PAGES = {629--649},
      ISSN = {0012-9593,1873-2151},
   MRCLASS = {14E05 (14J45 14M20)},
  MRNUMBER = {3098425},
MRREVIEWER = {Karl\ Schwede},
       DOI = {10.24033/asens.2197},
       URL = {https://doi.org/10.24033/asens.2197},
}

@article {MR662120,
    AUTHOR = {Mori, Shigefumi},
     TITLE = {Threefolds whose canonical bundles are not numerically
              effective},
   JOURNAL = {Ann. of Math. (2)},
  FJOURNAL = {Annals of Mathematics. Second Series},
    VOLUME = {116},
      YEAR = {1982},
    NUMBER = {1},
     PAGES = {133--176},
      ISSN = {0003-486X},
   MRCLASS = {14J30 (14E35)},
  MRNUMBER = {662120},
MRREVIEWER = {Ulf\ Persson},
       DOI = {10.2307/2007050},
       URL = {https://doi.org/10.2307/2007050},
}

@book {MR1440180,
    AUTHOR = {Koll\'ar, J\'anos},
     TITLE = {Rational curves on algebraic varieties},
    SERIES = {Ergebnisse der Mathematik und ihrer Grenzgebiete. 3. Folge. A
              Series of Modern Surveys in Mathematics [Results in
              Mathematics and Related Areas. 3rd Series. A Series of Modern
              Surveys in Mathematics]},
    VOLUME = {32},
 PUBLISHER = {Springer-Verlag, Berlin},
      YEAR = {1996},
     PAGES = {viii+320},
      ISBN = {3-540-60168-6},
   MRCLASS = {14-02 (14C05 14E05 14F17 14J45)},
  MRNUMBER = {1440180},
MRREVIEWER = {Yuri\ G.\ Prokhorov},
       DOI = {10.1007/978-3-662-03276-3},
       URL = {https://doi.org/10.1007/978-3-662-03276-3},
}

@article {MR1191735,
    AUTHOR = {Campana, F.},
     TITLE = {Connexit\'e{} rationnelle des vari\'et\'es de {F}ano},
   JOURNAL = {Ann. Sci. \'Ecole Norm. Sup. (4)},
  FJOURNAL = {Annales Scientifiques de l'\'Ecole Normale Sup\'erieure.
              Quatri\`eme S\'erie},
    VOLUME = {25},
      YEAR = {1992},
    NUMBER = {5},
     PAGES = {539--545},
      ISSN = {0012-9593},
   MRCLASS = {14J45},
  MRNUMBER = {1191735},
MRREVIEWER = {Luciana\ Picco Botta},
       URL = {http://www.numdam.org/item?id=ASENS_1992_4_25_5_539_0},
}

\end{document}